\documentclass[10pt]{article}
   \usepackage{amsmath,amsfonts,amsthm,amssymb,amscd,graphicx}
\theoremstyle{plain}
\newtheorem{thm}{Theorem}[section]
\newtheorem{prop}[thm]{Proposition}
\newtheorem{lemma}[thm]{Lemma}
\newtheorem{cor}[thm]{Corollary}

\theoremstyle{definition}
\newtheorem{defn}[thm]{Definition}

\theoremstyle{remark}
\newtheorem{rmk}[thm]{Remark}

\newcommand{\vol}{\ensuremath{\mathsf{vol}}}
\newcommand{\hol}{\ensuremath{\mathrm{Hol}_g (M)}}

\newcommand{\G}{\ensuremath{\operatorname{G_2}}}
\newcommand{\SP}{\ensuremath{\operatorname{Spin}(7)}}
\newcommand{\SUtw}{\ensuremath{\mathrm{SU}(2)}}
\newcommand{\SUth}{\ensuremath{\mathrm{SU}(3)}}
\newcommand{\Gs}{\ensuremath{\operatorname{G_2}}{-structure}}
\newcommand{\SPs}{\ensuremath{\operatorname{Spin}(7)}{-structure}}

\newcommand{\ph}{\ensuremath{\varphi}}
\newcommand{\ps}{\ensuremath{\psi}}
\newcommand{\wzeo}{\ensuremath{\Omega^0_1}}
\newcommand{\wons}{\ensuremath{\Omega^1_7}}
\newcommand{\wtws}{\ensuremath{\Omega^2_7}}
\newcommand{\wtwf}{\ensuremath{\Omega^2_{14}}}
\newcommand{\wtho}{\ensuremath{\Omega^3_1}}
\newcommand{\wths}{\ensuremath{\Omega^3_7}}
\newcommand{\wtht}{\ensuremath{\Omega^3_{27}}}
\newcommand{\wfoo}{\ensuremath{\Omega^4_1}}
\newcommand{\wfos}{\ensuremath{\Omega^4_7}}
\newcommand{\wfot}{\ensuremath{\Omega^4_{27}}}
\newcommand{\wfis}{\ensuremath{\Omega^5_7}}
\newcommand{\wfif}{\ensuremath{\Omega^5_{14}}}

\newcommand{\wthpos}{\ensuremath{\Omega^3_{\text{pos}}}}
\newcommand{\st}{\ensuremath{\ast}}
\newcommand{\hk}{\mathbin{\! \hbox{\vrule height0.3pt width5pt 
depth 0.2pt \vrule height5pt width0.4pt depth 0.2pt}}}
\newcommand{\stph}{\ensuremath{\ast \varphi}}
\newcommand{\lieg}{\ensuremath{\operatorname{\mathfrak {g}_2}}}
\newcommand{\tr}{\ensuremath{\operatorname{Tr}}}
\newcommand{\ws}{\ensuremath{w^{\flat}}}
\newcommand{\us}{\ensuremath{u^{\flat}}}
\newcommand{\vs}{\ensuremath{v^{\flat}}}
\newcommand{\as}{\ensuremath{a^{\flat}}}
\newcommand{\bs}{\ensuremath{b^{\flat}}}
\newcommand{\cs}{\ensuremath{c^{\flat}}}
\newcommand{\ds}{\ensuremath{d^{\flat}}}
\newcommand{\Xs}{\ensuremath{X^{\flat}}}
\newcommand{\ddx}[1]{\ensuremath{\frac{\del}{\del\:\!\! x^{#1}}}}
\newcommand{\dx}[1]{\ensuremath{d\:\!\! x^{#1}}}
\newcommand{\nab}[1]{\ensuremath{\nabla_{\! \! #1 \,}}}

\newcommand{\nph}{\ensuremath{\widetilde \varphi}}

\newcommand{\nstph}{\ensuremath{\tilde \st \widetilde \varphi}}

\newcommand{\del}{\ensuremath{\partial}}

\newcommand{\ddt}{\ensuremath{\frac{\partial}{\partial t}}}
\newcommand{\sgn}{\ensuremath{\mathrm{sgn}}}
\newcommand{\dl}{\ensuremath{\delta}}

\numberwithin{equation}{section}
\numberwithin{table}{section}
\numberwithin{figure}{section}

\begin{document}

\title{Flows of \G-Structures, I.}
\author{Spiro Karigiannis \\ {\it Mathematical Institute, University
of Oxford} \\ {\tt karigiannis@maths.ox.ac.uk} }

\maketitle

\begin{abstract}
  This is a foundational paper on flows of \Gs s.  We use local
  coordinates to describe the four torsion forms of a \Gs\ and derive
  the evolution equations for a general flow of a \Gs\ $\ph$ on a
  $7$-manifold $M$.  Specifically, we compute the evolution of the
  metric $g$, the dual $4$-form $\ps$, and the four independent
  torsion forms. In the process we obtain a simple new proof of a
  theorem of Fern\'andez-Gray.

  As an application of our evolution equations, we derive an analogue
  of the second Bianchi identity in \G-geometry which appears to be
  new, at least in this form. We use this result to derive explicit
  formulas for the Ricci tensor and part of the Riemann curvature
  tensor in terms of the torsion. These in turn lead to new proofs of
  several known results in \G-geometry.
\end{abstract}

\newpage

\tableofcontents

\newpage

\section{Introduction} \label{introsec}

\subsection{Overview} \label{overviewsec}

This work can be considered as a foundational paper on
geometric flows on manifolds with \Gs. This article also serves a
second purpose, which is to collect together many useful identities,
relations, and computational techniques for manifolds with \Gs, some
of which have never before appeared in the literature, although they
may be known to experts in the field.

Most research in geometric flows, such as Ricci flow or mean curvature
flow, is done using local coordinates, and this paper describes such
an approach for flows of \Gs s. One advantage to such an approach is
that we can appeal to known results on the evolution of quantities
which depend only on the variation of the metric, such as the
Christoffel symbols and the Riemann, Ricci, and scalar curvatures, as
can be found, for example, in~\cite{CK}.  The author hopes that the
foundational results collected here will be useful to others,
particularly researchers working on geometric flows who may not yet be
familiar with \Gs s.

A general evolution of a \Gs\ is described by a symmetric $2$-tensor
$h$ and a vector field $X$, and it is only $h$ which affects the
evolution of the associated Riemannian metric. On the other hand, for
the various geometric quantities unique to \Gs s, such as the four
torsion forms, their evolution is determined by both $h$ and $X$.

To date the only attempt at considering a flow of \Gs s has been the
Laplacian flow discussed in~\cite{Br3}. When restricted to {\em
closed} \Gs s, this flow is actually the gradient flow (with respect
to an appropriately chosen inner product) for the Hitchin functional
introduced in~\cite{Hi1} (arXiv version.) This flow is still not very
well understood.

It is not clear if the Laplacian flow is really the `natural' flow to
consider for \Gs s, or indeed if there even exists a natural flow at
all in this context. The advantage of developing the general theory of
flows of \Gs s is that it allows us to examine the evolution equations
for the torsion under a general flow in the hopes that one can find
obvious choices for specific flows which might have nice properties.
This is clearly more efficient than trying to guess the right flow and
computing all the associated evolution equations from scratch each
time.

The main results of this paper are Theorem~\ref{fulltorsionevolutionthm}
on the evolution of the full torsion tensor under a general flow of
\Gs s, and Theorem~\ref{secondbianchithm} which proves a Bianchi-like
identity relating the Riemann curvature and intrinsic torsion of a \Gs.

In the rest of Section~\ref{introsec}, we review our notation and
conventions. In Section~\ref{g2structuressec}, we review \Gs s, the
decomposition of the space of forms, and the intrinsic torsion forms
of a \Gs. Along the way we provide a simple new computational proof of
a Theorem of Fern\'andez-Gray~\cite{FG}. Section~\ref{generalflowssec}
is the heart of the paper, where we compute the evolution equations
for the metric, dual $4$-form, and torsion forms for a general flow of
\Gs s. In Section~\ref{bianchisec}, we apply our evolution equations
to derive Bianchi-type identities in \G-geometry, and use these to
produce new proofs of several known results. We also obtain an
explicit formula for the Ricci tensor of a general \Gs\ in terms of
the torsion. The paper closes with two appendices; the first of which
discusses various identities in \G-geometry, some of which are new and
which should be useful in other contexts as well, and the second
collecting some standard facts about flows of metrics.

The author is currently preparing a sequel~\cite{KY} to this paper in
collaboration with S.T. Yau, in which we present a detailed analysis
of several specific flows, including the Laplacian flow and the Ricci
flow in the context of \G-geometry. We also discuss short-time existence and soliton solutions.

\subsection{Notation and Conventions} \label{notationsec}

In this section we set up our notation and conventions. Throughout
this paper, $M$ is a (not necessarily compact) smooth manifold of
dimension $7$ which admits a \Gs. (See Section~\ref{g2reviewsec} for a
review of \Gs s.) The Einstein summation convention is employed
throughout: an index which appears both as a subscript and as a
superscript in the same term is summed from $1$ to $7$.

We use $S_7$ to denote the group of permutations of seven letters, and
$\sgn(\sigma)$ denotes the sign ($\pm 1$) of an element $\sigma$ of
$S_7$.

The space of $k$-forms on $M$ will be denoted by $\Omega^k$. It is the
space of sections of the bundle $\Lambda^k (T^* M)$. A differential
$k$-form $\alpha$ on $M$ will be written as
\begin{equation*}
\alpha = \frac{1}{k!} \alpha_{i_1 i_2 \cdots i_k} \, \dx{i_1} \wedge
\dx{i_2} \wedge \cdots \wedge \dx{i_k}
\end{equation*}
in local coordinates $(x^1, \ldots, x^7)$, where the sums are all from
$1$ to $7$, and $\alpha_{i_1 i_2 \cdots i_k}$ is completely
skew-symmetric in its indices. With this convention $\alpha$ can also
be written as
\begin{equation*}
\alpha = \sum_{i_1 < i_2 < \cdots < i_k} \, \alpha_{i_1 i_2 \cdots i_k}
\dx{i_1} \wedge \dx{i_2} \wedge \cdots \wedge \dx{i_k}
\end{equation*}
but we will not have need to do so. The advantage of this approach is
that if we take the interior product $\ddx{m} \hk \alpha$ of the
$k$-form $\alpha$ with a vector field $\ddx{m}$, we obtain the
$(k-1)$-form
\begin{equation*}
\ddx{m} \hk \alpha = \frac{1}{(k-1)!} \alpha_{m i_1 i_2 \cdots
i_{k-1}} \, \dx{i_1} \wedge \dx{i_2} \wedge \cdots \wedge \dx{i_{k-1}}
\end{equation*}

Given a Riemannian metric $g$ on $M$, it
induces a metric on $k$-forms which is defined on decomposable
elements to be
\begin{eqnarray*}
g( \dx{i_1} \wedge \cdots \wedge \dx{i_k} , \dx{j_1} \wedge
\cdots \wedge \dx{j_k} ) \, & = & \, \det_{a,b = 1, \ldots , k}(
g( \dx{i_a}, \dx{j_b}) ) \, = \, \det(g^{i_a j_b} ) \\ & =
& \sum_{\sigma \in S_7} \sgn(\sigma) g^{i_1 j_{\sigma(1)}} g^{i_2
j_{\sigma(2)}} \ldots g^{i_k j_{\sigma(k)}} 
\end{eqnarray*}
where $g^{ij} = g( \dx{i}, \dx{j} )$ is the induced metric
on the cotangent bundle and $g^{ij}$ is the inverse matrix of the
matrix $g_{ij}$. With this convention, one can check that the
inner product of two $k$-forms $\alpha = \frac{1}{k!} \alpha_{i_1
\cdots i_k} \, \dx{i_1} \wedge \cdots \wedge \dx{i_k}$ and $\beta =
\frac{1}{k!}  \beta_{j_1 \cdots j_k} \, \dx{j_1} \wedge \cdots \wedge
\dx{j_k}$ is
\begin{equation} \label{formsmetriceq}
g(\alpha , \beta) = \frac{1}{k!} \, \alpha_{i_1 \cdots
i_k} \beta_{j_1 \cdots j_k} g^{i_1 j_1} \ldots g^{i_k j_k}
\end{equation}
which differs from other conventions by the factor of $k!$.

The Levi-Civita covariant derivative associated to $g$ is denoted by
$\nabla$, and its associated Christoffel symbols by $\Gamma^k_{ij}$,
where $\Gamma^k_{ij} = \Gamma^k_{ji}$. We write $\nab{i}$ for
covariant differentiation in the $\ddx{i}$ direction. In local
coordinates, we have $\nab{i} (\dx{k}) = - \Gamma^{k}_{ij} \dx{j}$.
If $T_{i_1 \cdots i_k}$ is a tensor of type $(0,k)$, then $\nab{m}
T_{i_1 \cdots i_k}$ always means $(\nab{m} T)_{i_1 \cdots i_k}$, which
is 
\begin{equation} \label{coveq}
\nab{m} T_{i_1 \cdots i_k} = \ddx{m} T_{i_1 \cdots i_k} -
\Gamma^l_{mi_1} T_{l i_2 \cdots i_k} - \Gamma^l_{mi_2} T_{i_1 l \cdots
i_k} - \cdots - \Gamma^l_{mi_k} T_{i_1 i_2 \cdots l}
\end{equation}
and is a tensor of type $(0,k+1)$. Because the metric $g$ is parallel
with respect to $\nabla$, covariant differentiation commutes with
contractions. Explicitly, we have
\begin{equation*}
\nab{m} \left( T_{ai_1 \cdots i_k} S_{bj_1 \cdots j_l} g^{ab} \right)
= \left( \nab{m} T_{ai_1 \cdots i_k} \right) S_{bj_1 \cdots j_l}
g^{ab} + T_{ai_1 \cdots i_k} \left( \nab{m} S_{bj_1 \cdots j_l}
\right) g^{ab}
\end{equation*}
where $T$ is a $(0,k+1)$ tensor and $S$ is a $(0,l+1)$ tensor.

The exterior derivative $d\alpha$ of a $k$-form $\alpha$ can be
written in terms of the covariant derivative as
\begin{eqnarray} \nonumber
d\alpha & = & \frac{1}{k!} \, \left( \nab{m} \alpha_{i_1 \cdots i_k} +
\Gamma^l_{mi_1} \alpha_{l i_2 \cdots i_k} + \cdots + \Gamma^l_{mi_k}
\alpha_{i_1 i_2 \cdots l} \right) \dx{m} \wedge \dx{i_1} \wedge \cdots
\wedge \dx{i_k} \\ \label{dformeq} & = & \frac{1}{k!} (\nab{m}
\alpha_{i_1 \cdots i_k}) \, \dx{m} \wedge \dx{i_1} \cdots \wedge \dx{i_k}
\end{eqnarray}
since the $\Gamma^k_{ij}$'s are symmetric in $i,j$. We will always
write exterior derivatives of forms in this way.

The metric $g$ determines a `musical' isomorphism between the
tangent and cotangent bundles of $M$. If $v$ is a vector field, then
the {\em metric dual} $1$-form $\vs$ is defined by $\vs(w) = g(v, w)$
for all vector fields $w$. In coordinates, $\left( \ddx{i}
\right)^{\flat} = g_{ik} \dx{k}$. Similarly a $1$-form $\alpha$ has a
metric dual vector field $\alpha^{\sharp}$ defined by $\beta
(\alpha^{\sharp}) = g(\alpha, \beta)$ for all $1$-forms $\beta$, and
$\left( \dx{i}\right)^{\sharp} = g^{ik} \ddx{k}$. This isomorphism is
an isometry: $g(\vs, \ws) = g(v,w)$.

We use `$\vol$' to denote the volume form on $M$ associated to a
metric $g$ and an orientation, rather than something like $d \vol$, to
avoid confusion, since the volume form is never exact on a compact
manifold. In local coordinates the volume form is
\begin{equation*}
\vol = \sqrt{\det(g)} \, \dx{1} \wedge \cdots \wedge \dx{7}
\end{equation*}
where $\det(g)$ is the determinant of the matrix $g_{ij} = 
g( \ddx{i}, \ddx{j} )$.

The metric and orientation together determine the Hodge star
operator $\st$ taking $k$-forms to $(7-k)$-forms, characterized by the
relation
\begin{equation*}
\alpha \wedge \st \beta = g( \alpha, \beta ) \, \vol
\end{equation*}
on two $k$-forms $\alpha$ and $\beta$. We also have $\st^2 = 1$. Suppose
$v$ is a vector field and $\alpha$ is a $k$-form. The interior
product, wedge product, and Hodge star operator are all related by the
following identities (the signs here are all specific to the odd dimension
$7$):
\begin{eqnarray} \label{iosrelationseq}
& & \st ( v \hk \alpha) = (-1)^{k+1} (\vs \wedge \st \alpha) \\
\nonumber & & \st (v \hk \st \alpha) = (-1)^k (\vs \wedge \alpha) \\
\nonumber & & v \hk \vol = \st \vs
\end{eqnarray}
More details can be found, for example, in~\cite{K1}.

We also have need to consider the adjoint $\dl = d^{\st} :
\Omega^k \to \Omega^{k-1}$ of the exterior derivative, with respect to
the metric. This operator is called the {\em coderivative} and it
satisfies
\begin{equation*}
\int_M g(d\alpha, \beta) \, \vol = \int_M g(\alpha, \dl \beta) \, \vol
\end{equation*}
whenever $M$ is compact without boundary. The operator $\dl$ can be
written in terms of $d$ and $\st$ as
\begin{equation} \label{deldefneq}
\dl = {(-1)}^k \st d \st \qquad \text{ on \, } \Omega^k
\end{equation}
where again, the signs are specific to the odd dimension $7$ which we
consider exclusively throughout this paper. The coderivative $\dl$ can
be written in terms of the metric $g$ and the covariant derivative
$\nabla$ as follows:
\begin{eqnarray}
\nonumber \dl \alpha & = & \frac{1}{(k-1)!} (\dl \alpha)_{i_1 i_2 \cdots i_{k-1}}
\, \dx{i_1} \wedge \dx{i_2} \wedge \cdots \wedge \dx{i_{k-1}} \\
\label{delcoordinateseq} \text{where \, } (\dl \alpha)_{i_1 i_2 \cdots i_{k-1}}
& = & - g^{lm} \, \nab{l} \alpha_{m i_1 \cdots i_{k-1}}
\end{eqnarray}

We need the expression for the Lie derivative $\mathcal{L}_Y (\alpha)$
of a tensor $\alpha$ in the direction of a vector field $Y$ in terms
of the covariant derivative. Suppose $\alpha = \alpha_{i_1 \cdots i_k}$
is a tensor, where $\alpha_{i_1 \cdots i_k} = \alpha( \ddx{i_1}, \ldots, 
\ddx{i_k})$. Then by the Liebnitz rule for both $\mathcal{L}_Y$ and $\nab{Y}$,
we have
\begin{eqnarray} \nonumber
(\mathcal{L}_Y (\alpha))_{i_1 \cdots i_k} & = & Y (\alpha_{i_1 \cdots i_k}) 
- \sum_{j=1}^k \alpha(\ddx{i_1}, \ldots, \mathcal{L}_Y (\ddx{i_j}), \ldots, 
\ddx{i_k}) \\ \nonumber & = &  Y (\alpha_{i_1 \cdots i_k}) 
- \sum_{j=1}^k \alpha(\ddx{i_1}, \ldots, \nab{Y} (\ddx{i_j}) - \nab{i_j} Y,
\ldots, \ddx{i_k}) \\ \label{liederivativeeq} & = &
(\nab{Y} \alpha)_{i_1 \cdots i_k} + (\nab{i_1} Y^l) \alpha_{l i_2 \cdots i_k}
+ \cdots + (\nab{i_k} Y^l) \alpha_{i_1 \cdots i_{k-1} l}
\end{eqnarray}

Finally, we discuss our conventions for labelling the Riemann curvature tensor.
We define
\begin{equation*}
R^m_{ijk} \ddx{m} = (\nab{i} \nab{j} - \nab{j} \nab{i})(\ddx{k})
\end{equation*}
in terms of coordinate vector fields (so the usual Lie bracket term vanishes), and
we choose to lower indices by
\begin{equation*}
R_{ijkl} = R^m_{ijk} g_{ml}
\end{equation*}
With this convention, the Ricci tensor must be defined as $R_{jk} = R_{ijkl} g^{il}$
to ensure that the round sphere has positive Ricci curvature. Recall that we have
\begin{equation*}
R_{ijkl} = -R_{jikl} = -R_{ijlk} = R_{klij}
\end{equation*}
and the first Bianchi identity
\begin{equation} \label{riemannbianchieq}
R_{ijkl} + R_{iklj} + R_{iljk} = 0
\end{equation}
We will also need the Ricci identity
\begin{equation} \label{ricciidentityeq}
\nab{k} \nab{i} X_l - \nab{i} \nab{k} X_l = - R_{kilm} X^m
\end{equation}

Two possible references for this section (although the conventions do not
always agree with ours) are~\cite{CLN, J}.

\subsection{Acknowledgements} \label{acksec}

The bulk of this paper was completed while the author was a Postdoctoral
Fellow at the Mathematical Sciences Research Institute in 2006-2007 as part
of the program on `Geometric Evolution Equations and Related Topics.'
The inspiration for this work are the papers~\cite{Br3} and~\cite{Hi1}
by Robert Bryant and Nigel Hitchin, respectively. These were the first
papers (as far as the author knows) to discuss flows in \G-geometry,
and served to motivate the author to consider the general situation of
\G-flows.

The author is indebted to Robert Bryant, Hsiao-Bing Cheng, Ben Chow,
Dominic Joyce, Naichung Conan Leung, Lei Ni, and Andrejs Treibergs for
useful discussions. The author also benefited from a brief talk with
Richard Hamilton, from whom the author learned valuable
intuition. Finally, the author would like to acknowledge his former
thesis advisor Shing-Tung Yau, without whose constant encouragement and
advice, none of this work would have been possible.

\section{Manifolds with \Gs} \label{g2structuressec}

In this section we review the concept of a \Gs\ on a manifold $M$ and
the associated decompositions of the space of forms. More details
about \Gs s can be found, for example, in~\cite{Br3, J4, K1, Sa}. We
also describe explicitly the four torsion tensors associated to a \Gs\
and compute their expressions in local coordinates. These results are
needed to determine the evolution equations of the torsion tensors in
Section~\ref{torsionevolutionsec}.

\subsection{Review of \Gs s} \label{g2reviewsec}

Consider a $7$-manifold $M$ with a \G\ structure
$\ph$. Such a structure exists if and only if $M$ is orientable
and spin, which is equivalent to the vanishing of the first and
second Stiefel-Whitney classes $w_1(M) = w_2(M) = 0$. In fact the
space of $3$-forms $\ph$ on $M$ which determine a \Gs\ is an open
subbundle $\wthpos$ of the bundle $\Omega^3$ of $3$-forms on $M$,
sometimes called the bundle of {\em positive} $3$-forms. Such a
structure determines a Riemannian metric and an orientation in a
non-linear fashion which we now describe. Given local coordinates
$x^1, \ldots, x^7$, we define
\begin{equation} \label{Bijdefn}
B_{ij} \, \dx{1} \wedge \ldots \wedge \dx{7} = (\ddx{i} \hk \ph) \wedge
(\ddx{j} \hk \ph) \wedge \ph
\end{equation}
It is clear that $B_{ij} = B_{ji}$. This top form can be shown to
be equal to
\begin{equation} \label{Bijdefn2}
B_{ij} \, \dx{1} \wedge \ldots \wedge \dx{7} = -6 g_{ij} \vol = -6 g_{ij}
\sqrt{\det(g)} \, \dx{1} \wedge \ldots \wedge \dx{7}
\end{equation}
In other words, the $3$-form $\ph$ naturally determines the tensor
product of the metric $g_{ij}$ with the volume form $\vol$. From
this we can extract the metric $g_{ij}$ and subsequently the
volume form $\vol$ as follows.
\begin{eqnarray*}
B_{ij} & = & \left( (\ddx{i} \hk \ph) \wedge (\ddx{j} \hk \ph)
\wedge \ph \right) \, \left( \ddx{1}, \ldots, \ddx{7} \right) \\
B_{ij} & = & -6 g_{ij} \sqrt{\det(g)} \\ \det(B) & = & (-6)^7 \,
\det(g) \det(g)^{\frac{7}{2}} = -6^7 \, \det(g)^{\frac{9}{2}} \\
\sqrt{\det(g)} & = & - \frac{1}{6^{\frac{7}{9}}} \,
\det(B)^{\frac{1}{9}} \\ g_{ij} & = & - \frac{1}{6} \,
\frac{B_{ij}}{\sqrt{\det(g)}} = \frac{1}{6^{\frac{2}{9}}} \,
\frac{B_{ij}}{\det(B)^{\frac{1}{9}}}
\end{eqnarray*}
Some authors chose a different convention in which the $(-6)$
factor above is actually a $(+6)$. Since this cancels out in the
above calculation, the equation
\begin{equation} \label{gijeq}
\boxed{\, \, g_{ij} = \frac{1}{6^{\frac{2}{9}}} \,
\frac{B_{ij}}{\det(B)^{\frac{1}{9}}} \, \,}
\end{equation}
is true for both sign conventions (see~\cite{K2}.)

We stress that the tensor $B$ is actually a section of
$\mathrm{Sym}^2 (T^* M) \otimes \wedge^7 (T^* M)$. Therefore, if
we change to a new basis $\dx{i} = P^i_l d\tilde x^l$, then
$\tilde B_{ij} = P^k_i P^l_j \det(P) B_{kl}$ from which it follows
from~\eqref{gijeq} that $\tilde g_{ij} = P^k_i P^l_j g_{kl}$ as
expected.

The metric $g$ and orientation (determined by the volume form)
determine a Hodge star operator $\st$, and we therefore have the
associated dual $4$-form $\ps = \st \ph$.  The metric also determines
the Levi-Civita connection $\nabla$, and the manifold $(M, \ph)$ is
called a \G\ manifold if $\nab{} \ph = 0$. Note that this is a
nonlinear partial differential equation for $\ph$, since $\nabla$
depends on $g$ which depends non-linearly on $\ph$. Such manifolds
(where $\ph$ is parallel) have Riemannian holonomy $\hol$ contained in
the exceptional Lie group $\G \subset \mathrm{SO}(7)$.

The \Gs\ corresponding to the $3$-form $\ph$ is also called {\em
torsion-free} if $\ph$ is parallel with respect to the metric
$g_{\ph}$ induced by $\ph$. Torsion-free \Gs s have sometimes been
called `integrable' by some authors, but since that term has also been
used in different contexts as well, we prefer to stick to the
unambigous `torsion-free' at all times. See also~\cite{Br3} for more
discussion about the sometimes confusing terminology about \Gs s.

In~\cite{FG}, the following theorem is proved.

\smallskip

\begin{thm}[Fern\'andez-Gray, 1982] \label{FGthm}
The \Gs\ corresponding to $\ph$ is torsion-free if and only if $\ph$
is both closed and co-closed:
\begin{equation*}
d \ph = 0 \qquad \qquad d \st_{\ph} \ph = d \ps = 0
\end{equation*}
\end{thm}
More recent proofs of this theorem can be found in~\cite{BS, J4,
Br3}.

\smallskip

\begin{rmk} \label{FGrmk}
A differential form $\ph$ is harmonic if $(dd^* + d^*d)\ph = 0$. On a compact manifold,
this is equivalent to $d\ph = 0$ and $d \st \ph = 0$. 
Since a parallel differential form is always closed and co-closed, 
Theorem~\ref{FGthm} says that for a compact manifold with $\G$~structure $\ph$, the $3$-form  $\ph$ being parallel is equivalent to it being harmonic (with respect to its induced metric.)
\end{rmk}

\smallskip

\begin{rmk} \label{holonomyrmk}
It is well known (see~\cite{J4}, for example), that for a \G-manifold
the holonomy must be one of the following four possibilities:
\begin{eqnarray*}
\hol = \{ 1 \} \quad \, \! & \Leftrightarrow & b_1(M) = 7 \\ \hol = \SUtw &
\Leftrightarrow & b_1(M) = 3 \\ \hol = \SUth & \Leftrightarrow &
b_1(M) = 1 \\ \hol = \G \quad \, \, & \Leftrightarrow & b_1(M) = 0
\end{eqnarray*}
We are interested in constructing
manifolds with full holonomy $\G$, and not a strictly smaller
subgroup. Therefore we can assume that the fundamental group $\pi_1 (M)$
is finite.
\end{rmk}

\smallskip

\begin{rmk} \label{topobstructionsrmk}
There are some topological obstructions (in the compact case)
to the existence of a torsion-free \Gs\ which are known. First,
we need $b_3(M) \geq 1$, since the $3$-form $\ph$ is a
non-zero harmonic $3$-form, and hence represents a non-trivial cohomology
class by the Hodge theorem. Additionally, if we insist on full holonomy
\G\, rather than a strictly smaller subgroup, then we must have $b_1(M) = 0$,
which was already mentioned in Remark~\ref{holonomyrmk}, and also that the
first Pontryagin class $p_1 (M)$ of the manifold must be non-zero.
There are also some conditions on the cohomology ring structure. A detailed
discussion can be found in~\cite{J4}. Of course, sufficient conditions
for the existence of a torsion-free \Gs\
(analogous to the Calabi conjecture in K\"ahler geometry) are far 
from being known.
\end{rmk}

\subsection{Decomposition of the space of forms} \label{formssec}

The existence of a \Gs\ $\ph$ on $M$ (with no condition on $\nab{}
\ph$) determines a decomposition of the spaces of differential
forms on $M$ into irreducible \G\ representations. This is
analogous to the decomposition of complex-valued differential
forms on an almost complex manifold into forms of type $(p,q)$. We
will see explicitly that the spaces $\Omega^2$ and 
$\Omega^3$ of $2$-forms and $3$-forms decompose as
\begin{eqnarray*}
\Omega^2 & = & \wtws \oplus \wtwf \\
\Omega^3 & = & \wtho \oplus \wths \oplus \wtht
\end{eqnarray*}
where $\Omega^k_l$ has (pointwise) dimension $l$ and this
decomposition is orthogonal with respect to the metric $g$.
The spaces $\wtws$ and $\wths$ are both isomorphic to the
cotangent bundle $\wons = T^*M$ (and hence also to the tangent
bundle $TM$.) We show below that the space $\wtwf$ is isomorphic
to the Lie algebra $\lieg$, and $\wtht$ is isomorphic to the traceless
symmetric $2$-tensors $\mathrm{Sym}^2_0 (T^*M)$ on $M$. The explicit
descriptions are as follows:
\begin{eqnarray}
\label{wtwsdesc} \wtws & = & \{ X \hk \ph;\, X \in \Gamma(TM) \}
\quad = \quad \{ \beta \in \Omega^2; \, \st ( \ph \wedge \beta ) =
-2 \beta \} \\ \label{wtwfdesc} \wtwf & = & \{ \beta \in \Omega^2;
\, \beta \wedge \ps = 0 \} \, \, \, \quad = \quad \{ \beta \in
\Omega^2; \, \st ( \ph \wedge \beta) = \beta \} \\
\label{wthodesc} \wtho & = & \{ f \ph; \, f \in C^{\infty}(M) \}
\\ \label{wthsdesc} \wths & = & \{ X \hk \ps; \, X \in \Gamma(TM)
\} \\ \label{wthtdesc} \wtht & = & \{ h_{ij} g^{jl} \, \dx{i}
\wedge \left( \ddx{l} \hk \ph \right) \, ; \, h_{ij} =
h_{ji} \, , \, \tr_g(h_{ij}) = g^{ij} h_{ij} = 0 \}
\end{eqnarray}
The decompositions $\Omega^4 = \wfoo \oplus \wfos \oplus \wfot$
and $\Omega^5 = \wfis \oplus \wfif$ are obtained by taking the
Hodge star of the decompositions of $\Omega^3$ and $\Omega^2$,
respectively.

\smallskip

\begin{rmk} \label{eigenvaluessignrmk}
There is another orientation convention for \Gs s which differs from this
one. In the other convention, the eigenvalues of the operator $\beta
\mapsto \st (\ph \wedge \beta)$ are $+2$ and $-1$ instead of $-2$ and
$+1$, respectively. See~\cite{K2} for more on sign and orientation
conventions in \G\ geometry.
\end{rmk}

\smallskip

We now establish the decompositions of the spaces $\Omega^2$ and $\Omega^3$
in detail. Let $\beta = \frac{1}{2} \beta_{ij} \, \dx{i} \wedge \dx{j}$ be an
arbitrary $2$-form. Then we have
\begin{eqnarray*}
\st ( \ph \wedge \beta) & = & \frac{1}{2} \, \beta_{ij} \st (\dx{i}
\wedge \dx{j} \wedge \ph) = \frac{1}{2} \, \beta_{ij} g^{il} \ddx{l}
\hk \st (\dx{j} \wedge \ph) \\ & = & -\frac{1}{2} \, \beta_{ij} g^{il}
g^{jm} \ddx{l} \hk \ddx{m} \hk \st \ph = -\frac{1}{2} \, \beta_{ij}
g^{il} g^{jm} \left(\frac{1}{2} \ps_{mlab} \, \dx{a} \wedge \dx{b}
\right) \\ & = & \frac{1}{4} \, \beta_{ij} \ps_{lmab} g^{il} g^{jm} \,
\dx{a} \wedge \dx{b}
\end{eqnarray*}
where we have used~\eqref{iosrelationseq} twice. This computation
allows us to express the projection operators $\pi_{7}$ and $\pi_{14}$
from $\Omega^2$ to $\wtws$ and $\wtwf$, respectively, as
follows. From~\eqref{wtwsdesc} and~\eqref{wtwfdesc}, we see that
\begin{equation*}
\pi_7 (\beta) = \frac{1}{3} \beta - \frac{1}{3} \st (\ph \wedge \beta)
\qquad \qquad \qquad \pi_{14} (\beta) = \frac{2}{3} \beta +
\frac{1}{3} \st (\ph \wedge \beta)
\end{equation*}
Hence for $\beta = \frac{1}{2} \, \beta_{ij} \dx{i} \wedge
\dx{j}$, we have
\begin{eqnarray} \label{wtwsprojeq}
\pi_7 (\beta) & = & \frac{1}{2} \, \left( \frac{1}{3} \beta_{ab} -
\frac{1}{6} \beta_{ij} g^{il} g^{jm} \ps_{lmab} \right) \dx{a} \wedge \dx{b} \\
\label{wtwfprojeq} \pi_{14} (\beta) & = & \frac{1}{2} \, \left( \frac{2}{3}
\beta_{ab} + \frac{1}{6} \beta_{ij} g^{il} g^{jm} \ps_{lmab} \right)
\dx{a} \wedge \dx{b}
\end{eqnarray}
The alternative characterizations of $\wtws$ and $\wtwf$ given
in~\eqref{wtwsdesc} and~\eqref{wtwfdesc} can also be expressed in
local coordinates: $\beta \in \wtws$ if and only if $\beta = X \hk
\ph$ for some vector field $X$, which is equivalent to $\beta_{ij} =
X^k \ph_{ijk}$, and $\beta \in \wtwf$ if and only if $\psi \wedge
\beta = 0$, which, using~\eqref{iosrelationseq}, becomes $\beta_{ij}
g^{il} g^{jm} \ph_{lmk} = 0$. We summarize the descriptions of
$\Omega^2$ in the following:

\smallskip

\begin{prop} \label{wtwprop}
Let $\beta = \frac{1}{2} \, \beta_{ij} \, \dx{i} \wedge \dx{j}$ be in
$\Omega^2$. Then
\begin{eqnarray*}
& & \beta \in \wtws \quad \Leftrightarrow \quad \beta_{ij} g^{il}
g^{jm} \ps_{lmab} = -4 \, \beta_{ab} \quad \Leftrightarrow \quad
\beta_{ij} = X^k \ph_{ijk} \\ & & \beta \in \wtwf \quad
\Leftrightarrow \quad \beta_{ij} g^{il} g^{jm} \ps_{lmab} = 2 \,
\beta_{ab} \quad \Leftrightarrow \quad \beta_{ij} g^{il} g^{jm}
\ph_{lmk} = 0
\end{eqnarray*}
\end{prop}

\smallskip

\begin{rmk}
It is an enlightening exercise to establish the above equivalences (in
index notation) directly using the identities of
Lemmas~\ref{g2identities1lemma}, \ref{g2identities2lemma},
and~\ref{g2identities3lemma}.
\end{rmk}

\smallskip

\begin{rmk} \label{wtwvecsrmk}
Because the space $\wtws$ is isomorphic to the space of vector fields
(and also to the space of $1$-forms), it is convenient to be able to
go back and forth between the two descriptions. If $X = X^k \ddx{k}$
is a vector field, the associated $\wtws$ form, by
Proposition~\ref{wtwprop}, is $\frac{1}{2} X_{ab} \, \dx{a}
\wedge \dx{b}$, where $X_{ab} = X^l \ph_{lab}$. It is easy to check using
Lemma~\ref{g2identities1lemma} that we can solve for $X^k$ in this expression
as $X^k = \frac{1}{6} X_{ab} \ph_{mpq} g^{ap} g^{bq} g^{mk}$.
In summary,
\begin{equation} \label{wtwvecseq}
X_{ab} = X^l \ph_{lab} \quad \Leftrightarrow \quad 
X^k = \frac{1}{6} \, X_{ab} \ph_{mpq} g^{ap} g^{bq} g^{mk}
\end{equation}
which we will have occasion to use frequently.
\end{rmk}

\smallskip

We need a few more relations involving $\wtws$ and $\wtwf$ which will be used
to simplify the evolution equations for the torsion forms in
Section~\ref{torsionevolutionsec}.

\smallskip

\begin{lemma} \label{twoformsidentitieslemma}
Suppose $\beta_{ij}$ is a $2$-form. Then if $\beta_{ij} \in \wtwf$,
\begin{equation} \label{wtwfidentityeq}
\beta_{ab} g^{bl} \ph_{lpq} = \beta_{pl} g^{lm} \ph_{maq} - \beta_{ql} g^{lm}
\ph_{map}
\end{equation}
whereas if $\beta \in \wtws$, then
\begin{equation} \label{wtwsidentityeq}
\beta_{ab} g^{bl} \ph_{lpq} = -\frac{1}{2} \, \beta_{pl} g^{lm} \ph_{maq} + 
\frac{1}{2} \, \beta_{ql} g^{lm} \ph_{map} - \frac{3}{2} \, g_{pa} \beta_q +
\frac{3}{2} \, g_{qa} \beta_p
\end{equation}
where in this case $\beta_k = \frac{1}{6} \, \beta_{ij} g^{ia} g^{jb} \ph_{abk}$
as given by~\eqref{wtwvecseq}.
\end{lemma}
\begin{proof}
We prove both statements at once. By Proposition~\ref{wtwprop}, we can write
$\beta_{ab} = \lambda \beta_{ij} g^{im} g^{jn} \ps_{mnab}$ where $\lambda =
\frac{1}{2}$ if $\beta \in \wtwf$ and $\lambda = -\frac{1}{4}$ if $\beta \in
\wtws$. Now we use Lemma~\ref{g2identities2lemma} and compute
\begin{eqnarray*}
& & \beta_{ab} g^{bl} \ph_{lpq} = \lambda \beta_{ij} g^{im} g^{jn} (\ph_{pql}
\ps_{mnab} g^{lb}) \\ & = & \lambda \beta_{ij} g^{im} g^{jn} ( g_{pm} \ph_{qna}
+ g_{pn} \ph_{mqa} + g_{pa} \ph_{mnq} - g_{qm} \ph_{pna} - g_{qn} \ph_{mpa}
- g_{qa} \ph_{mnp}) \\ & = & \lambda( \beta_{pj} g^{jn} \ph_{qna} + \beta_{ip}
g^{im} \ph_{mqa} + 6 g_{pa} \beta_q - \beta_{qj} g^{jn} \ph_{pna} -
\beta_{iq} g^{im} \ph_{mpa} - 6 g_{qa} \beta_p) 
\end{eqnarray*}
where $\beta_k = 0$ if $\beta \in \wtwf$ by Proposition~\ref{wtwprop}. This can
be simplified to
\begin{equation*}
\beta_{ab} g^{bl} \ph_{lpq} = \lambda( \, 2 \beta_{pl} g^{lm} \ph_{maq} - 
2 \beta_{ql} g^{lm} \ph_{map} + 6 g_{pa} \beta_q - 6 g_{qa} \beta_p \, )
\end{equation*}
and the statements now follow by substituting the value of $\lambda$ in each case.
\end{proof}

\smallskip

\begin{cor} \label{g2liealgebracor}
The space $\wtwf$ is a Lie algebra with respect to the commutator of matrices:
\begin{equation*}
[ \beta, \mu]_{ij} = \beta_{il} g^{lm} \mu_{mj} - \mu_{il} g^{lm} \beta_{mj}
\end{equation*}
\end{cor}
\begin{proof}
We know that $\Omega^2 \cong \mathfrak{so}(7)$, and that the matrix commutator
is a Lie algebra bracket on $\Omega^2$ which satisfies the Jacobi identity. We
need to show is that the bracket of two elements of $\wtwf$ is again in $\wtwf$.
Suppose $\beta$ and $\mu$ are in $\wtwf$. Then $[\beta, \mu] \in \wtwf$ if
and only if
\begin{equation*}
[\beta, \mu]_{ij} g^{ia} g^{jb} \ph_{abc} = 0
\end{equation*}
We have
\begin{eqnarray*}
[\beta, \mu]_{ij} g^{ia} g^{jb} \ph_{abc} & = & \beta_{il} g^{lm} \mu_{mj} g^{ia}
g^{jb} \ph_{abc} - \mu_{il} g^{lm} \beta_{mj} g^{ia} g^{jb} \ph_{abc} \\
& = & \beta_{il} g^{lm} \mu_{mj} g^{ia} g^{jb} \ph_{abc} + \mu_{il} g^{lm}
g^{ia} ( \beta_{ak} g^{kn} \ph_{nmc} - \beta_{ck} g^{kn} \ph_{nma})
\end{eqnarray*}
where we have used~\eqref{wtwfidentityeq} in the last line above since $\beta \in
\wtwf$. The first two terms cancel each other, and the last term vanishes
by Proposition~\ref{wtwprop}, since $\mu \in \wtwf$.
\end{proof}

\smallskip

\begin{rmk} \label{g2liealgebrarmk}
Of course, $\wtwf \cong \lieg$, the Lie algebra of $\G$.
\end{rmk}

\smallskip

We can regard a vector field $X^k$ as a $\wtws$ form $X_{ab}$
by~\eqref{wtwvecseq}.  Given a $2$-form $\beta_{ij}$, it acts on $X^k$
by matrix multiplication: $(\beta(X))^m = \beta_{ij} X^j g^{im}$ to
give another vector field $\beta(X)$, which we can turn into a
$\wtws$-form by $(\beta(X))_{ab} = (\beta(X))^m \ph_{mab}$.

\smallskip

\begin{prop} \label{wtwactionvecsprop}
If $\beta \in \wtwf$, then
\begin{equation} \label{wtwfactioneq}
(\beta(X))_{ab} = [ \beta, X ]_{ab}
\end{equation}
whereas if $\beta \in \wtws$, then
\begin{equation} \label{wtwsactioneq}
(\beta(X))_{ab} = -\frac{1}{2} [\beta, X]_{ab} - \frac{3}{2} \, \beta_a X_b
+ \frac{3}{2} \, \beta_b X_a 
\end{equation}
where $[ \beta, X]$ is the matrix commutator of $\beta$ and $X$ regarded
as elements of $\Omega^2$.
\end{prop}
\begin{proof}
We use~\eqref{wtwfidentityeq} and compute, if $\beta \in \wtwf$:
\begin{eqnarray*}
(\beta(X))_{ab} & = & \beta_{ij} X^j g^{il} \ph_{lab} = - X^j (\beta_{al} g^{lm}
\ph_{mjb} - \beta_{bl} g^{lm} \ph_{mja}) \\ & = & \beta_{al} g^{lm} X_{mb} -
X_{am} g^{lm} \beta_{lb} = [\beta, X]_{ab} 
\end{eqnarray*}
and similarly using~\eqref{wtwsidentityeq} in the case $\beta \in \wtws$ to
establish~\eqref{wtwsactioneq}.
\end{proof}

\smallskip

To conclude our discussion of $\Omega^2$, we note that if $\beta \in
\wtws$, we can also write $(\beta(X))^m = \beta_{ij} X^j g^{im} =
\beta^k X^j \ph_{kij} g^{im}$, and thus $\beta(X) = - (X \hk
\beta \hk \ph)^{\sharp} = - (\beta \times X)$, in
terms of the cross product, where we have used~\eqref{crosseq}.
Therefore equation~\eqref{wtwsactioneq} can be restated as
\begin{equation} \label{wtwsactioneq2}
\beta_a X_b - \beta_b X_a = -\frac{1}{3} [ \beta, X ]_{ab} + \frac{2}{3} \, 
(\beta^{\flat} \times X)_{ab}
\end{equation}
However, we see also that
\begin{equation} \label{wtwsactioneq3}
\frac{1}{6} (\beta_a X_b - \beta_b X_a) g^{ai} g^{bj} \ph_{ijk} = \frac{1}{3}
 \, \beta^i X^j \ph_{ijk} = \frac{1}{3} (\beta \times X)_k
\end{equation}
and therefore, considering the $\wtws$ parts of~\eqref{wtwsactioneq} as vectors
via~\eqref{wtwvecseq}, and using~\eqref{wtwsactioneq3} gives
\begin{equation*}
-( \beta \times X) = -\frac{1}{2} \, \pi_7 ([ \beta, X]) - \frac{3}{2} (\frac{1}{3}
\, \beta \times X)
\end{equation*}
from which it follows that
\begin{equation} \label{pi7commutatorcrosseq}
\pi_7 ([ \beta, X]) = \beta \times X
\end{equation}
for $\beta \in \wtws$, which will be used in Section~\ref{torsionevolutionsec}.

\smallskip

The decomposition of the space $\Omega^3$ of $3$-forms can be
understood by considering the infinitesmal action of $\mathrm{GL}(7,
\mathbb R)$ on $\ph$. Let $A = A^i_l \in \mathfrak{gl}(7, \mathbb R)$.
Hence $e^{At} \in \mathrm{GL}(7, \mathbb R)$, and we have
\begin{equation*}
e^{A t} \cdot \ph = \frac{1}{6} \, \ph_{ijk} \, (e^{At} \dx{i}) \wedge 
(e^{At} \dx{j}) \wedge (e^{At} \dx{k}) 
\end{equation*}
Differentiating with respect to $t$ and setting $t = 0$, we obtain:
\begin{eqnarray*}
\left. \frac{d}{dt} \right|_{t=0} (e^{At} \cdot \ph) & = & \frac{1}{6}
\, \left( A^l_i \ph_{ljk} + A^l_j \ph_{ilk} + A^l_k \ph_{ijl} \right) \dx{i} \wedge
\dx{j} \wedge \dx{k} \\ & = & \frac{1}{2} \, A^l_i \ph_{ljk} \, \dx{i} \wedge
\dx{j} \wedge \dx{k} = A^l_i \, \dx{i} \wedge \left( \ddx{l} \hk \ph \right)
\end{eqnarray*}
Now let $A^l_i = g^{lj} A_{ij}$, and decompose $A_{ij} = S_{ij} +
C_{ij}$ into symmetric and skew-symmetric parts, where $S_{ij} =
\frac{1}{2} (A_{ij} + A_{ji})$ and $C_{ij} = \frac{1}{2} (A_{ij} -
A_{ji})$. We have a map
\begin{eqnarray*}
D & : &  \mathrm{GL}(7, \mathbb R) \to \Omega^3 \\
D & : & A \mapsto \left. \frac{d}{dt} \right|_{t=0} (e^{At} \cdot \ph) \\
& & = S_{ij} g^{jl} \, \dx{i} \wedge \left( \ddx{l} \hk \ph \right) +
C_{ij} g^{jl} \, \dx{i} \wedge \left( \ddx{l} \hk \ph \right)  
\end{eqnarray*}
\begin{prop} \label{wtwfdescprop} The kernel of $D$ is isomorphic to
the subspace $\wtwf$. It is also isomorphic to the Lie algebra
$\lieg$ of the Lie group \G\ which is the subgroup of
$\mathrm{GL}(7, \mathbb R)$ which preserves $\ph$.
\end{prop}
\begin{proof}
Since we are defining \G\ to be the group preserving $\ph$, the
kernel of $D$ is isomorphic to $\lieg$ by definition. To
show explicitly that this is isomorphic to $\wtwf$, decompose
$C_{ij} = (C_7)_{ij} + (C_{14})_{ij}$, where $C_7 \in \wtws$ and
$C_{14} \in \wtwf$. We have
\begin{equation*}
(C_{14})_{ij} g^{jl} \, \dx{i} \wedge \left( \ddx{l} \hk \ph \right) =
\frac{1}{6} \,  \left( (C_{14})^l_i \ph_{ljk} + (C_{14})^l_j \ph_{ilk} +
(C_{14})^l_k \ph_{ijl} \right) \dx{i} \wedge \dx{j} \wedge \dx{k} 
\end{equation*}
From Proposition~\ref{wtwprop}, we have $(C_{14})_{ij} = \frac{1}{2}
\, (C_{14})_{ab} g^{ap} g^{bq} \ps_{pqij}$. Using this together with the
final equation of Lemma~\ref{g2identities2lemma}, one can compute easily that
\begin{equation*}
(C_{14})^l_i \ph_{ljk} + (C_{14})^l_j \ph_{ilk} + (C_{14})^l_k \ph_{ijl}
= -2 \left( (C_{14})^l_i \ph_{ljk} + (C_{14})^l_j \ph_{ilk} +
(C_{14})^l_k \ph_{ijl} \right)
\end{equation*}
and hence that $(C_{14})_{ij} g^{jl} \, \dx{i} \wedge \left( \ddx{l}
\hk \ph \right) = 0$. Therefore $\wtwf$ is in the kernel of $D$. We will
see below that $D$ restricted to $\wtws$ or the symmetric tensors $S^2(T)$
is injective. This completes the proof.
\end{proof}

By counting dimensions, we must have $\wths = \{(C_{7})_{ij} g^{jl} \,
\dx{i} \wedge \left( \ddx{l} \hk \ph \right) \}$ and also $\wtho
\oplus \wtht = \{S_{ij} g^{jl} \, \dx{i} \wedge \left( \ddx{l} \hk \ph
\right)\}$. We now proceed to establish these explicitly.

To show $(C_{7})_{ij} g^{jl} \, \dx{i} \wedge \left(
\ddx{l} \hk \ph \right)$ is $X \hk \ps$ for some vector field $X$, we
use Proposition~\ref{wtwprop} to write $(C_7)_{ij} = (C_7)^k \ph_{kij}$, where
$(C_7)^k = \frac{1}{6} \, (C_7)_{ij} g^{ia} g^{jb} \ph_{abc} g^{kc}$, 
by~\eqref{wtwvecseq}. In fact since $(C_{14})_{ij} g^{ia} g^{jb} \ph_{abc} = 0$,
we actually have $(C_7)^k = \frac{1}{6} \, C_{ij} g^{ia} g^{jb} \ph_{abc} g^{kc}$.
Now we use Lemma~\ref{g2identities1lemma} and compute:
\begin{eqnarray*}
& & \frac{1}{6} \,  \left( (C_{7})_{il} g^{lm} \ph_{mjk} + (C_{7})_{jl} g^{lm} \ph_{imk} +
(C_{7})_{kl} g^{lm} \ph_{ijm} \right) \\ & = & \frac{1}{6} \, \left( (C_7)^n \ph_{nil}
g^{lm} \ph_{mjk} + (C_7)^n \ph_{njl} g^{lm} \ph_{imk} + (C_7)^n \ph_{nkl} g^{lm} \ph_{ijm}
\right) \\ & = & \frac{1}{6} \, (C_7)^n \left( g_{nj} g_{ik} - g_{nk} g_{ij} - \ps_{nijk}
+ g_{nk} g_{ji} - g_{ni} g_{jk} - \ps_{njki} + g_{ni} g_{kj} - g_{nj} g_{ki} - \ps_{nkij}
\right) \\ & = & \frac{1}{6} \, (-3 (C_7)^n ) \ps_{nijk} = \frac{1}{6} \, X^n \ps_{nijk}  
\end{eqnarray*}
and hence we have shown that for $C_{ij}$ skew-symmetric,
\begin{eqnarray} \label{skewthreeformeq}
& & C_{il} g^{lm} \ph_{mjk} + C_{jl} g^{lm} \ph_{imk} +
C_{kl} g^{lm} \ph_{ijm} = X^n \ps_{nijk} \\ \nonumber & & \text{ where } \quad
X^n = -\frac{1}{2} \, C_{ij} g^{ia} g^{jb} \ph_{abc} g^{cn}
\end{eqnarray}

Following the notation of Bryant~\cite{Br3}, we define maps
$i : S^2(T) \to \Omega^3$ and $j : \Omega^3 \to S^2(T)$ as follows
(our definition of the map $i$ differs from Bryant's by a factor of $2$):
\begin{eqnarray} \label{idefeq}
i(h_{ij}) & = & h_{ij} g^{jl} \, \dx{i} \wedge \left( \ddx{l} \hk \ph \right) =
\frac{1}{2} \, h^l_i \ph_{ljk} \, \dx{i} \wedge \dx{j} \wedge \dx{k} \\
\label{jdefeq} (j(\eta))(v,w) & = & \st\left( ( v \hk \ph) \wedge (w \hk \ph)
\wedge \eta \right) 
\end{eqnarray}

\smallskip

We will have several occasions to use the following.
\begin{prop} \label{symmstarprop}
Suppose that $h_{ij}$ is a symmetric tensor. It corresponds to the
form $\eta = i(h_{ij})$ in $\Omega^3$ given by
\begin{equation*}
\eta = h_{ij} g^{jl} \, \dx{i} \wedge \left( \ddx{l} \hk \ph \right) =
\frac{1}{2} \, h^l_i \ph_{ljk} \, \dx{i} \wedge \dx{j} \wedge \dx{k}
\end{equation*}
Then the Hodge star $\st \eta$ of $\eta$ is
\begin{equation*}
\st \eta =  \left( \frac{1}{4} \tr_g(h) \, g_{ij} - h_{ij} \right) g^{jl} \dx{i}
\wedge \left( \ddx{l} \hk \ps \right)
\end{equation*}
where $\tr_g(h) = g^{ij} h_{ij}$.
\end{prop}
\begin{proof}
We compute
\begin{eqnarray*}
\st \eta & = & h^l_i \st \left( \dx{i} \wedge \left( \ddx{l} \hk \ph
\right) \right) = h^l_i g^{im} \ddx{m} \hk \st \left( \ddx{l}
\hk \ph \right) \\ & = & h^{lm} \ddx{m} \hk \left( g_{lk} \, \dx{k}
\wedge \ps \right) = h^m_k \left( \delta^k_m \ps - \dx{k}
\wedge \left( \ddx{m} \hk \ps \right) \right) \\ & = & \frac{1}{4} \tr_g(h) \,
\delta^l_i \, \dx{i} \wedge \left( \ddx{l} \hk \ps \right)
- h^l_i \, \dx{i} \wedge \left( \ddx{l} \hk \ps \right)
\end{eqnarray*}
where we have used~\eqref{iosrelationseq} several times.
\end{proof}

\smallskip

\begin{prop} \label{symmmetricprop}
Suppose $f_{ij}$ and $h_{ij}$ are two symmetric tensors. Let $i(f)$ and
$i(h)$ be their corresponding forms in $\Omega^3$. Then we have
\begin{equation*}
i(f) \wedge \st (i(h)) = g(i(f), i(h)) \, \vol = \left( \tr_g(f) \tr_g(h) + 
2 f^k_l h^l_k \right) \vol
\end{equation*}
\end{prop}
\begin{proof}
Using Proposition~\ref{symmstarprop}, we compute $i(f) \wedge \st (i(h)) =$
\begin{eqnarray*}
& & f^l_i \, \dx{i} \wedge \left( \ddx{l}
\hk \ph \right) \wedge \left( \frac{1}{4} \tr_g(h) \delta^m_k - h^m_k \right) \,
\dx{k} \wedge \left( \ddx{m} \hk \ps \right) \\ & = & f^{la}
\left( \frac{1}{4} \tr_g(h) g^{mb} - h^{mb} \right) \, \left(\ddx{a}\right)^{\flat}
\wedge \left(\ddx{b} \right)^{\flat} \wedge \left( \ddx{l} \hk \ph
\right) \wedge \left( \ddx{m} \hk \ps \right) \\ & = & f^{la}
\left( \frac{1}{4} \tr_g(h) g^{mb} - h^{mb} \right)
(2 \, g_{al} g_{bm} - 2 \, g_{am} g_{bl} + \ps_{ablm} ) \, \vol \\ & =
& \left( \frac{14}{4} \tr_g(f) \tr_g(h) - 2 \tr_g(f) \tr_g(h) -
\frac{2}{4} \tr_g(f) \tr_g(h) + 2 \, f^l_m h^m_l \right) \vol + 0 \\ & = & 
\left( \tr_g(f) \tr_g(h) + 2 f^k_l h^l_k \right) \vol
\end{eqnarray*}
using Proposition~\ref{vecvecphipsiprop} and the symmetry of $f_{ij}$ and $h_{ij}$.
\end{proof}

\smallskip

\begin{cor} \label{imapinjectivecor}
The map $i : S^2(T) \to \Omega^3$ is injective. It is therefore an isomorphism onto
its image, $\wtho \oplus \wtht$.
\end{cor}
\begin{proof}
Suppose $i(h) = 0$. Then Proposition~\ref{symmmetricprop} gives
$\tr_g(h)^2 + g^{ia} g^{jb} h_{ij} h_{ab} = 0$. The second term is just $g(h,h)$,
the natural inner product on $S^2(T)$. Thus both terms are non-negative and
hence vanish. Therefore $g(h,h) = 0$, so $h_{ij} = 0$.
\end{proof}

\smallskip

\begin{prop} \label{threeformsisomorphismsprop} The map $j : \Omega^3 \to S^2(T)$
is an isomorphism between $\wtho \oplus \wtht$ and $S^2(T)$.
Consequently, $\wths$ is the kernel of $j$. Explicitly, we have
\begin{eqnarray*}
& & \text{If \, } \eta = h_{ij} g^{jl} \, \dx{i} \wedge \left( \ddx{l} \hk
\ph \right) + (X \hk \ps) = i(h) + (X \hk \ps) \\ & & \text{then \, } j(\eta) = 
- 2 \tr_g(h) g_{ij} - 4 h_{ij}
\end{eqnarray*}
\end{prop}
\begin{proof}
The map $j$ is linear. First let $\eta = X \hk \ps$. Then
\begin{equation*}
j(\eta) = \st\left( ( v \hk \ph) \wedge (w \hk \ph) \wedge (X \hk \ps) \right) = 0
\end{equation*}
by Theorem 2.4.7 of~\cite{K1}. Thus $\wths$ is in the kernel of the map $j$.
Now suppose that $\eta = h_{ij} g^{jl} \, \dx{i} \wedge \left( \ddx{l} \hk
\ph \right) = \frac{1}{2} h^l_i \ph_{ljk} \, \dx{l} \wedge \dx{j} \wedge \dx{k}$.
Then we have
\begin{eqnarray*}
(j(\eta))_{ab} & = & \frac{1}{2} h^l_i \ph_{ljk} \st \left( (\ddx{a} \hk \ph)
\wedge (\ddx{b} \hk \ph) \wedge \dx{i} \wedge \dx{j} \wedge \dx{k} \right) \\ & = & 
\frac{1}{8} h^l_i \ph_{ljk} \ph_{apq} \ph_{bmn} \, \st (\dx{p} \wedge \dx{q} \wedge
\dx{m} \wedge \dx{n} \wedge \dx{i} \wedge \dx{j} \wedge \dx{k}) \\ & = & \frac{1}{8}
\, \sum_{\sigma \in S_7} \sgn(\sigma) \ph_{a\sigma(1)\sigma(2)}
\ph_{b\sigma(3)\sigma(4)} \ph_{l\sigma(5)\sigma(6)} h^l_{\sigma(7)} \st (\dx{1}
\wedge \cdots \wedge \dx{7}) \\ & = & \frac{1}{8} \, \st \left( \frac{8}{3} (B_{ij}
h^l_l + B_{il} h^l_j + B_{jl} h^l_i ) \, \dx{1} \wedge \cdots \wedge \dx{7} \right)
\\ & = & \frac{1}{3} \st \left( -6 \tr_g(h) \, g_{ij} \vol - 6 h_{ij} \vol -
6 h_{ji}) \vol \right) \\ & = & -2 \tr_g(h) \, g_{ij} -4 h_{ij}
\end{eqnarray*}
where we have used~\eqref{wtwscubedeq} and $B_{ij} = -6 g_{ij} \sqrt{\det(g)}$. From
this it follows immediately that $j$ is injective on $\wtho \oplus \wtht$, for if
$j(i(h)) = -2 \tr_g(h) \, g_{ij} - 4 h_{ij} = 0$, taking the trace gives
$-18 \tr_g(h) = 0$ and hence $h_{ij} = 0$.
\end{proof}

\smallskip

To summarize, we have seen that an arbitrary $3$-form $\eta$ on a manifold $M$
with \Gs\ $\ph$ consists of the data of a vector field $X$ and a symmetric $2$-tensor
$h$. Explicitly, we have
\begin{eqnarray*}
\eta & = & h_{ij} g^{jl} \, \dx{i} \wedge \left( \ddx{l} \hk
\ph \right) + X^l \ddx{l} \hk \ps \\ & = & \frac{1}{2} \, h^l_i
\ph_{ljk} \, \dx{i} \wedge \dx{j} \wedge \dx{k} + \frac{1}{6} \,
X^l \ps_{lijk} \, \dx{i} \wedge \dx{j} \wedge \dx{k} 
\end{eqnarray*}

\smallskip

\begin{rmk} \label{symmrmk}
Note that the symmetric $2$-tensor $h_{ij}$ decomposes as $h_{ij}
= \frac{1}{7} \tr_g(h) \, g_{ij} + h^0_{ij}$ where $h^0_{ij}$ is the
trace-free part of $h_{ij}$. Hence the first term in the above expression
can be written as
\begin{equation*}
\frac{3}{7} \, h \ph + \frac{1}{2}\, (h^0)^l_i \ph_{ljk} \, \dx{i}
\wedge \dx{j} \wedge \dx{k}
\end{equation*}
which shows explicitly the $\wtho$ and $\wtht$ components.
However, it is more convenient to consider the tensor
$h_{ij}$ directly, using the isomorphism $\wtho \oplus \wtht \cong
\mathrm{Sym}^2 (T^*M)$.
\end{rmk}

\subsection{The intrinsic torsion forms of a \Gs}
\label{torsionsec}

Using the decomposition of the spaces of forms on $M$ determined by
$\ph$, given in Section~\ref{formssec}, we can decompose $d\ph$ and $d
\ps$ into types. This defines the {\em torsion forms} of the \Gs.

\smallskip

\begin{defn} \label{torsiondef}
There are four independent torsion forms corresponding to a \Gs\ $\ph$.
Following the notation introduced in~\cite{Br3}, we denote them by
\begin{eqnarray*}
\tau_0 \, \, \in \, \, \wzeo \, \, & \qquad & \tau_1 \, \, \in \, \, \wons
\\ \tau_2 \, \, \in \, \, \wtwf & \qquad & \tau_3 \, \, \in \, \,
\wtht
\end{eqnarray*} 
They are defined via the equations
\begin{eqnarray} \label{torsiondefeq}
d \ph & = & \tau_0 \, \ps + 3 \tau_1 \wedge \ph + \st \tau_3 \\
\nonumber d \ps & = & 4 \tau_1 \wedge \ps + \st \tau_2
\end{eqnarray}
The constants are chosen for convenience. The fact that the $\wfos$
component of $d\ph$ and the $\wfis$ component of $d\ps$ are the same,
up to a constant (when projected back to $\wons$ via the \G-invariant
isomorphisms), is something that needs to be proved. A nice
representation-theoretic proof is described in~\cite{BS, Br3}. Below
we will give a brute force computational proof of this fact, which is
useful as it is an exercise in the type of manipulations that will be
used frequently in Section~\ref{torsionevolutionsec}, when we compute
the evolution equations for the torsion forms.
\end{defn}

\smallskip

\begin{rmk} \label{torsionnamesrmk}
We call $\tau_0$ the {\em scalar torsion}, $\tau_1$ the {\em
vector torsion}, $\tau_2$ the {\em Lie algebra torsion}, and $\tau_3$
the {\em symmetric traceless torsion}. These names are non-standard,
but are clearly reasonable.
\end{rmk}

\smallskip

\begin{rmk} \label{torsiondefrmk1}
In~\cite{Br3}, the second defining equation is given as $d\ps = 4
\tau_1 \wedge \ps + \tau_2 \wedge \ph$. This is the same thing, since
in our chosen orientation convention, $\tau_2 \wedge \ph = \st \tau_2$
for $\tau_2 \in \wtwf$. We prefer this way of writing the equation to
make it more symmetric with the $d\ph$ equation.
\end{rmk}

\smallskip

\begin{rmk} \label{torsiondefrmk2}
In~\cite{K1}, a one-form $\theta$ is discussed, sometimes called the
{\em Lee form}. It is easy to check that $\theta = -12 \tau_1$.
\end{rmk}

\smallskip

\begin{thm} \label{onetauonethm}
The appearance of the same one-form $\tau_1$ in the expressions
for both $d\ph$ and $d\ps$ in~\eqref{torsiondefeq} above is justified.
\end{thm}
\begin{proof}
We begin by not assuming that the two $\tau_1$'s are the same. Let $d
\ph = \tau_0 \, \ps + 3 \tau_1 \wedge \ph + \st \tau_3$ and $d \ps =
4 \tilde \tau_1 \wedge \ps + \st \tau_2$. We must show that $\tilde
\tau_1 = \tau_1$. We manipulate these relations as follows
\begin{eqnarray*}
d \ph = \tau_0 \, \ps + 3 \tau_1 \wedge \ph + \st \tau_3 \quad \, \, &
\qquad & \qquad \quad d \ps = 4 \tilde \tau_1 \wedge \ps + \st \tau_2
\\ \st (d \ph) = \tau_0 \, \ph + 3 \st( \tau_1 \wedge \ph) + \tau_3 \,
\! & \qquad & \quad \, \, \, \, \st (d \ps) = 4 \st (\tilde \tau_1
\wedge \ps) + \tau_2 \\ \ph \wedge \st (d \ph) = 0 - 3 \ph \st ( \ph
\wedge \tau_1) + 0 \, \, \, \, \, & \qquad & \ps \wedge \st (d \ps) =
4 \ps \wedge \st (\ps \wedge \tilde \tau_1) + 0 \\ \ph \wedge \st (d
\ph) = 12 \st \tau_1 \qquad \qquad \qquad \quad & \qquad & \ps \wedge
\st (d \ps) = 12 \st \tilde \tau_1
\end{eqnarray*}
where we have used Proposition~\ref{g2relationsprop}. Therefore we see that
\begin{eqnarray*}
\tau_1 = \tilde \tau_1 & \Leftrightarrow & \ph \wedge \st (d\ph) = \ps
\wedge \st (d\ps) \\ & \Leftrightarrow & \dx{p} \wedge \ph \wedge \st
(d\ph) = \dx{p} \wedge \ps \wedge \st (d \ps) \qquad \text{ for all
$p$} \\ & \Leftrightarrow & g(d \ph , \dx{p} \wedge \ph ) = g(d\ps ,
\dx{p} \wedge \ps) \qquad \text{ for all $p$}
\end{eqnarray*}
and it is the last equality above that we will now establish.

Let $X = X_i \dx{i}$ be an arbitrary one-form. Then we have
\begin{eqnarray*}
X \wedge \ph & = & \frac{1}{6} \, X_q \ph_{ijk} \, \dx{q} \wedge
\dx{i} \wedge \dx{j} \wedge \dx{k} \\ & = & \frac{1}{24} \left( X_q
\ph_{ijk} - X_i \ph_{qjk} - X_j \ph_{iqk} - X_k \ph_{ijq} \right) \,
\dx{q} \wedge \dx{i} \wedge \dx{j} \wedge \dx{k} \\ & = & \frac{1}{24}
\, A_{qijk} \, \dx{q} \wedge \dx{i} \wedge \dx{j} \wedge \dx{k}
\end{eqnarray*}
where we have skew-symmetrized the coefficients. Similarly we have
\begin{eqnarray*}
d \ph & = & \frac{1}{6} \, \left( \nab{m} \ph_{abc} - \nab{a}
\ph_{mbc} - \nab{b} \ph_{amc} - \nab{c} \ph_{abm} \right) \, \dx{m}
\wedge \dx{a} \wedge \dx{b} \wedge \dx{c} \\ & = & \frac{1}{24} \,
B_{mabc} \, \dx{m} \wedge \dx{a} \wedge \dx{b} \wedge \dx{c}
\end{eqnarray*}
Now using~\eqref{formsmetriceq}, we have
\begin{eqnarray*}
g( X \wedge \ph , d\ph) & = & \frac{1}{24} \, A_{qijk} B_{mabc} g^{qm}
g^{ia} g^{jb} g^{kc} \\ & = & \frac{1}{6} \, \left( X_q \ph_{ijk} -
X_i \ph_{qjk} - X_j \ph_{iqk} - X_k \ph_{ijq} \right) \left( \nab{m}
\ph_{abc} \right) g^{qm} g^{ia} g^{jb} g^{kc}
\end{eqnarray*}
Let $X = \dx{p}$, so that $X_i = \delta^p_i$, and this expression
becomes
\begin{eqnarray*}
g( \dx{p} \wedge \ph , d\ph) & = & \frac{1}{6} \, \left( \delta^p_q
\ph_{ijk} - \delta^p_i \ph_{qjk} - \delta^p_j \ph_{iqk} - \delta^p_k
\ph_{ijq} \right) \left( \nab{m} \ph_{abc} \right) g^{qm} g^{ia}
g^{jb} g^{kc} \\ & = & \frac{1}{6} \, \ph_{ijk} \left( \nab{m}
\ph_{abc} \right) g^{pm} g^{ia} g^{jb} g^{kc} - \frac{1}{2} \,
\ph_{qjk} \left( \nab{m} \ph_{abc} \right) g^{qm} g^{pa} g^{jb} g^{kc}
\end{eqnarray*}
By Proposition~\ref{g2derivativeidentitiesprop}, the first term
vanishes, and the second term becomes
\begin{equation} \label{tauoneeq1}
g( \dx{p} \wedge \ph , d\ph) = \frac{1}{2} \, (\nab{m} \ph_{ijk} )
\ph_{abc} g^{im} g^{pa} g^{jb} g^{kc}
\end{equation}
An exactly analogous calculation, which we omit, yields the expression
\begin{equation} \label{tauoneeq2}
g( \dx{p} \wedge \ps , d\ps ) = \frac{1}{6} \, (\nab{m} \ps_{ijkl} )
\ps_{abcd} g^{im} g^{pa} g^{jb} g^{kc} g^{ld}
\end{equation}
Combining the two expressions, $g( \dx{p} \wedge \ph , d\ph) = g(
\dx{p} \wedge \ps , d\ps)$ if and only if
\begin{equation*}
(\nab{m} \ps_{ijkl} ) \ps_{abcd} g^{im} g^{jb} g^{kc} g^{ld} = 3 \,
(\nab{m} \ph_{ijk} ) \ph_{abc} g^{im} g^{jb} g^{kc}
\end{equation*}
But this is precisely the content of Proposition~\ref{nabphirelnabpsiprop},
after contracting with $g^{im}$.
\end{proof}

\smallskip

The reason to consider the torsion forms of a \Gs\ $\ph$ is because
$\ph$ is torsion-free if and only if all four torsion forms vanish,
and these forms are independent. This is clear since the decomposition
of $\Omega^k$ into \G-representations is orthogonal, and because the
maps $\alpha \mapsto \ph \wedge \alpha$ from $\wons \to \wfos$ and
$\alpha \mapsto \ps \wedge \alpha$ from $\wons \to \wfis$ are
isomorphisms.

Having defined the four torsion forms, we need to derive their
expressions in local coordinates, in order to be able to compute their
evolutions under a flow. To do this most efficiently, we will first
determine how to write $\nab{}\ph$ in terms of the four torsion
tensors. In the process we will obtain a new computational proof of
the Fern\'andez-Gray Theorem~\ref{FGthm}. We begin with the following important
observation.

\smallskip

\begin{lemma} \label{torsionsymmetrieslemma}
For any vector field $X$, the $3$-form $\nab{X} \ph$ lies in the subspace
$\wths$ of $\Omega^3$. Therefore, the covariant derivative $\nab{}\ph$ lies in 
the space $\wons \otimes \wths$, a $49$-dimensional space (pointwise.)
\end{lemma}
\begin{proof}
Let $X = \ddx{l}$, and consider the $3$-form $\nab{l} \ph$. An arbitrary
element $\eta$ of $\wtho \oplus \wtht$ can be written (for some symmetric tensor
$h_{ij}$) as
\begin{equation*}
\eta = \frac{1}{2} \, h^m_i \ph_{mjk} \, \dx{i} \wedge \dx{j} \wedge \dx{k} = 
\frac{1}{6} \, \left( h^m_i \ph_{mjk} + h^m_j \ph_{imk} + h^m_k \ph_{ijm} \right)
\dx{i} \wedge \dx{j} \wedge \dx{k}
\end{equation*}
Using~\eqref{formsmetriceq}, the inner product of $\eta$ with $\nab{l}\ph = 
\frac{1}{6} \, \nab{l} \ph_{abc} \, \dx{a} \wedge \dx{b} \wedge \dx{c}$ is
\begin{eqnarray*}
g(\nab{l}\ph, \eta) & = & \frac{1}{6} \, (\nab{l}\ph_{abc}) \left( h^m_i \ph_{mjk} +
h^m_j \ph_{imk} + h^m_k \ph_{ijm} \right) g^{ai} g^{bj} g^{ck} \\ & = & 
\frac{1}{2} \, (\nab{l}\ph_{abc}) h^m_i \ph_{mjk} g^{ai} g^{bj} g^{ck} = 
\frac{1}{2} \, (\nab{l}\ph_{abc}) h^{ma} \ph_{mjk} g^{bj} g^{ck} 
\end{eqnarray*}
which vanishes since the third equation of
Proposition~\ref{g2derivativeidentitiesprop} says
$(\nab{l}\ph_{abc}) \ph_{mjk} g^{bj} g^{ck}$ is skew-symmetric in $a$ and $m$.
Since $g(\nab{l}\ph, \eta) = 0$ for all $\eta \in \wtho \oplus \wtht$, we have
that $\nab{l}\ph \in \wths$ for all $l = 1, \ldots, 7$ as claimed.
\end{proof}
\begin{rmk} \label{torsionsymmetriesrmk1}
In~\cite{FG}, Fern\`andez and Gray study the symmetries of $\nab{}\ph$, and
deduce that $\nab{}\ph \in \wons \otimes \wths$ as we have just shown. Their
arguments involve the complicated analysis of the extension of the cross product
operation to the full exterior bundle. The above proof is, in our opinion, much
more transparent.
\end{rmk}

\smallskip

We pause here to consider Lemma~\ref{torsionsymmetrieslemma} in more detail.
The reason why $\nab{X} \ph \in \wths$ is essentially because of the way that
the $3$-form $\ph$ of a \Gs\ determines a metric $g$. We have seen that
\begin{equation*}
(v \hk \ph) \wedge (w \hk \ph) \wedge \ph = -6 g(v,w) \vol
\end{equation*}
Since the metric and volume form are parallel, the covariant derivative of this
gives
\begin{equation*}
(v \hk \nab{X}\ph) \wedge (w \hk \ph) \wedge \ph + (v \hk \ph) \wedge
(w \hk \nab{X}\ph) \wedge \ph + (v \hk \ph) \wedge (w \hk \ph) \wedge
\nab{X}\ph = 0
\end{equation*}
Using the identity $\ph \wedge (v \hk \ph) = -2 \st (v \hk \ph)$ from
Proposition~\ref{g2relationsprop}, this becomes
\begin{equation*}
-2 \, (v \hk \nab{X}\ph) \wedge \st (w \hk \ph) -2 \, (v \hk \ph) \wedge \st
(w \hk \nab{X}\ph) + (v \hk \ph) \wedge (w \hk \ph) \wedge \nab{X}\ph = 0
\end{equation*}
Finally, the covariant differentiation of the identity $(v \hk \ph) \wedge \st
(w \hk \ph) = 4 \, g(v,w) \vol$ (which is also implicit in
Proposition~\ref{g2relationsprop}) gives
\begin{equation*}
(v \hk \nab{X}\ph) \wedge \st (w \hk \ph) + (v \hk \ph) \wedge \st (w \hk \nab{X}\ph)
= 0
\end{equation*}
Combining this with the previous equation yields the important relation
\begin{equation*}
(v \hk \ph) \wedge (w \hk \ph) \wedge \nab{X}\ph = 0
\end{equation*}
But by the definition of the map $j$ in~\eqref{jdefeq}, this precisely says
that the $\wtho \oplus \wtht$ component of $\nab{X}\ph$ is zero.

\smallskip

\begin{defn} \label{fulltorsiondef}
Since $\nab{l} \ph \in \wths$, by~\eqref{wthsdesc} we can write
$\nab{l}\ph_{abc} = T_{lm} g^{mn} \ps_{nabc}$ for some $2$-tensor
$T_{lm}$, which we shall call the {\em full torsion tensor}.
\end{defn}

\smallskip

In coordinates, 
the four independent torsion forms are the following: $\tau_0$ (a function);
$\tau_1 = (\tau_1)_l \, \dx{l}$, a $1$-form, which by Remark~\ref{wtwvecsrmk}
can also be written as an $\wtws$-form $\tau_1 = \frac{1}{2} \, (\tau_1)_{ab} \,
\dx{a} \wedge \dx{b}$ where $(\tau_1)_{ab} = (\tau_1)_l g^{lk} \ph_{kab}$;
$\tau_2 = \frac{1}{2} (\tau_2)_{ab} \, \dx{a} \wedge \dx{b}$, an $\wtwf$-form; and
$\tau_3 = \frac{1}{2} \, (\tau_3)_{im} g^{ml} \ph_{ljk} \, \dx{i} \wedge \dx{j}
\wedge \dx{k}$, where $(\tau_3)_{im}$ is traceless symmetric. The relationship
between the four torsion forms and the full torsion tensor $T_{lm}$ is given by
the following theorem.

\smallskip

\begin{thm} \label{torsionthm}
The covariant derivative $\nab{}\ph$ of the $3$-form $\ph$ can be written as
\begin{equation*}
\nab{l} \ph_{abc} = T_{lm} g^{mn} \ps_{nabc} 
\end{equation*}
where the full torsion tensor $T_{lm}$ is 
\begin{equation*}
T_{lm} = \frac{\tau_0}{4} g_{lm} - (\tau_3)_{lm} + (\tau_1)_{lm} -
\frac{1}{2} (\tau_2)_{lm}
\end{equation*}
\end{thm}
\begin{proof}
Write the full torsion tensor as $T_{lm} = S_{lm} + C_{lm}$, where
$S_{lm} = \frac{1}{2} (T_{lm} + T_{ml})$ and $C_{lm} = \frac{1}{2}
(T_{lm} - T_{ml})$ are the symmetric and skew-symmetric parts of $T_{lm}$.
Thus we have
\begin{equation} \label{torsionthmtempeq}
\nab{l}\ph_{abc} = (S_{lm} + C_{lm}) g^{mn} \ps_{nabc}
\end{equation}

Since $d\ph = \tau_0 \ps + 3 \, \tau_1 \wedge \ph + \st \tau_3$, the
$\wtho \oplus \wtht$ component of $\st d\ph$ is $\tau_0 \ph + \tau_3$,
which we write as $\frac{3}{7} \left( \frac{7}{3} \tau_0 \right) \ph + \tau_3$.
By Remark~\ref{symmrmk}, this is $f_{ij} g^{jl} \dx{i} \wedge \left( \ddx{l}
\hk \ph \right)$, where $f_{ij} = \frac{1}{7} \left( \frac{7}{3} \tau_0 \right)
g_{ij} + (\tau_3)_{ij}$. Therefore by Proposition~\ref{symmstarprop}, we have
that the $\wfoo \oplus \wfos$ component of $d\ph = \st (\st d\ph)$ is
$\left( \frac{1}{4} \tr_g(f) g_{ij} - f_{ij} \right) g^{jl} \dx{i} \wedge \left(
\ddx{l} \hk \ps \right)$. But $\tr_g(f) = \frac{7}{3} \tau_0$, so
\begin{equation} \label{torsionthmtempeq2}
\frac{1}{4} \tr_g(f) g_{ij} - f_{ij} = \frac{7}{12} \tau_0 \, g_{ij} -
\frac{4}{12} \tau_0 \, g_{ij} - (\tau_3)_{ij} = \frac{1}{4} \tau_0 \, g_{ij} -
(\tau_3)_{ij}
\end{equation}
Now we can also write $d\ph = \frac{1}{6} \, \nab{l}\ph_{abc} \, \dx{l} \wedge
\dx{a} \wedge \dx{b} \wedge \dx{c}$, which by~\eqref{torsionthmtempeq} is
\begin{eqnarray*}
d\ph & = & \frac{1}{6} \, \left( S_{lm} + C_{lm} \right) g^{mn} \ps_{nabc} \,
\dx{l} \wedge \dx{a} \wedge \dx{b} \wedge \dx{c} \\ & = & S_{lm} g^{mn} \dx{l} 
\wedge \left( \ddx{n} \hk \ps \right) + C_{lm} g^{mn} \dx{l} 
\wedge \left( \ddx{n} \hk \ps \right)
\end{eqnarray*}
The second term above is in $\wfos$. Therefore, comparing the $\wfoo \oplus \wfot$
term of $d\ph$ given above and by~\eqref{torsionthmtempeq2}, we see that
\begin{equation*}
S_{lm} = \frac{\tau_0}{4} \, g_{lm} - (\tau_3)_{lm}
\end{equation*}
This proves the first half of the theorem.

To establish the second half, we write $\dl \ph = - \st d \st \ph$ in two ways.
First, since $d\ps = 4 \, \tau_1 \wedge \ps + \st \tau_2$, we have
$\dl \ph = - \st d \ps = - 4 ( {\tau_1}^{\sharp} \hk \ph ) - \tau_2$,
using~\eqref{iosrelationseq}. Therefore
\begin{equation} \label{torsionthmtempeq3} 
\dl \ph = -4 \, \frac{1}{2} (\tau_1)_{ab} \dx{a} \wedge \dx{b}
- \frac{1}{2} (\tau_2)_{ab} \dx{a} \wedge \dx{b}
\end{equation}
From~\eqref{delcoordinateseq}, we also have $\dl \ph = -\frac{1}{2} \, g^{lk} \nab{l}
\ph_{kab} \, \dx{a} \wedge \dx{b}$. Using~\eqref{torsionthmtempeq}, this is
\begin{equation*}
\dl\ph = -\frac{1}{2} \, g^{lk} \left( S_{lm} + C_{lm} \right) g^{mn}
\ps_{nkab} \, \dx{a} \wedge \dx{b}
\end{equation*}
The first term vanished by the symmetry of $S_{lm}$ and skew-symmetry of
$\ps_{nkab}$. Now we decompose $C_{lm} = (C_7)_{lm} + (C_{14})_{lm}$ into 
$\wtws \oplus \wtwf$ components. Using Proposition~\ref{wtwprop}, and
interchanging $n$ and $k$, we see
\begin{eqnarray*}
\dl\ph & = & \frac{1}{2} \, g^{lk} \left( (C_7)_{lm} + (C_{14})_{lm} \right)
g^{mn} \ps_{knab} \, \dx{a} \wedge \dx{b} \\ & = & \frac{1}{2} \, \left(
-4(C_7)_{ab} + 2 (C_{14})_{ab} \right) \, \dx{a} \wedge \dx{b}
\end{eqnarray*}
Comparing this to~\eqref{torsionthmtempeq3}, we see that
$(C_7)_{ab} = (\tau_1)_{ab}$ and $(C_{14})_{ab} = -\frac{1}{2} (\tau_2)_{ab}$, 
hence
\begin{equation*}
C_{lm} = (\tau_1)_{lm} - \frac{1}{2} \, (\tau_2)_{lm}
\end{equation*}
and the proof is complete.
\end{proof}

\smallskip

\begin{cor} \label{torsionthmcor} The $3$-form $\ph$ is parallel if
and only if it is both closed and co-closed. (This is
Theorem~\ref{FGthm}.)
\end{cor}
\begin{proof}
  A parallel form is always automatically closed and co-closed, since
  the exterior derivative $d$ and the coderivative $\dl$ can both be
  written using the covariant derivative $\nabla$. The converse
  follows from Theorem~\ref{torsionthm}, since $d\ph = 0$ and $\dl \ph
  = 0$ are equivalent to the vanishing of all four torsion tensors, so
  $T_{lm} = 0$ and hence $\nab{l} \ph_{abc} = 0$.
\end{proof}

\smallskip

\begin{rmk} \label{nablapsirmk}
Starting from $\nab{l} \ph_{abc} = T_{lm} g^{mn} \ps_{nabc}$ and
using~\eqref{nablapsieq} and Lemma~\ref{g2identities2lemma}, it is 
an easy computation to show that
\begin{equation} \label{nablapsieq2}
\nab{m} \ps_{ijkl} = - T_{mi} \ph_{jkl} + T_{mj} \ph_{ikl}
- T_{mk} \ph_{ijl} + T_{ml} \ph_{ijk}
\end{equation}
which will be used in Section~\ref{bianchisec}.
\end{rmk}

\smallskip

The following lemma gives an explicit formula for $T_{lm}$ in terms
of $\nab{}\ph$. This will be used in Section~\ref{torsionevolutionsec} to
derive the evolution equations for the torsion tensors.

\smallskip

\begin{lemma} \label{fulltorsionlemma}
The full torsion tensor $T_{lm}$ is equal to
\begin{equation} \label{fulltorsioneq}
T_{lm} = \frac{1}{24} \, (\nab{l} \ph_{abc} ) \ps_{mijk} g^{ia} g^{jb} g^{kc}
\end{equation}
\end{lemma}
\begin{proof}
We begin with $\nab{l} \ph_{abc} = T_{lk} g^{kn} \ps_{nabc}$ and use
Lemma~\ref{g2identities3lemma} to compute:
\begin{eqnarray*}
\nab{l} \ph_{abc} \ps_{nijk} g^{ia} g^{jb} g^{kc} & = & T_{lk} g^{kn}
\ps_{nabc} \ps_{mijk} g^{ia} g^{jb} g^{kc} \\ & = & T_{lk} g^{kn}
(24 \, g_{nm}) = 24 \, T_{lm}
\end{eqnarray*}
as claimed.
\end{proof}

\smallskip

Next we present expressions for the four torsion tensors from which we
will calculate their evolution equations.

\smallskip

\begin{prop} \label{torsionformsprop}
The four torsion forms can be written in terms of $T_{pq} = S_{pq} + C_{pq}$
as follows:
\begin{eqnarray*}
\tau_0 & = & \frac{4}{7} g^{pq} S_{pq} \\
(\tau_3)_{pq} & = & \frac{1}{4} \, \tau_0 \, g_{pq} - S_{pq} \\
(\tau_1)_{pq} & = & \frac{1}{3} C_{pq} - \frac{1}{6} C_{ij} g^{ia} g^{jb}
\ps_{abpq} \\ (\tau_2)_{pq} & = & -\frac{4}{3} C_{pq} - \frac{1}{3} C_{ij}
g^{ia}g^{jb} \ps_{abpq}
\end{eqnarray*}
\end{prop}
\begin{proof}
This is immediate from Theorem~\ref{torsionthm} and
equations~\eqref{wtwsprojeq} and~\eqref{wtwfprojeq}.
\end{proof}

\smallskip

We close this section with an important observation about the vector torsion
$\tau_1$. The effect of a conformal scaling on \Gs s is well
understood, and a detailed discussion can be found in Section 3.1
of~\cite{K1}. In particular, we have the following result.

\smallskip

\begin{thm} \label{conformaltorsionthm}
Let $\nph = f^3 \ph$ be a new \G-form, where $f$ is any nowhere
vanishing smooth function, then the metric scales as $\tilde g = f^2
g$ and the $4$-form as $\nstph = f^4 \stph$. Furthermore, the torsion
forms transform as follows:
\begin{eqnarray*}
\tilde \tau_0 = & f^{-1} \tau_0 \qquad \qquad \qquad \tilde \tau_1 = &
\tau_1 + d \log (f) \\ \tilde \tau_2 = & f \, \tau_2 \qquad \qquad
\qquad \, \, \, \, \tilde \tau_3 = & f^2 \, \tau_3
\end{eqnarray*}
\end{thm}
\begin{proof}
This is essentially Theorem 3.1.4 in~\cite{K1}. In that paper, $\theta
= - 12 \, \tau_1$.
\end{proof}

\smallskip

The relevance of this result is evident. If we can construct a \Gs\
$\ph$ for which the three torsion forms $\tau_0$, $\tau_2$ and
$\tau_3$ all vanish, then any conformal scaling of such a structure
remains of such type. It is easy to check that in this case that
vector torsion $\tau_1$ is a {\em closed} $1$-form. If we are
interested in manifolds with full holonomy \G, then by
Remark~\ref{holonomyrmk}, we must necessarily start with a manifold
$M$ with finite fundamental group, and hence $H^1 (M) = 0$, so the
form $\tau_1$ is exact. Then we see from
Theorem~\ref{conformaltorsionthm} that we can always find a conformal
scaling factor $f$ (unique up to a multiplicative constant) to make
the vector torsion $\tau_1$ vanish as well.  Such \Gs s are called
{\em conformally parallel}. Hence we can restrict attention to
constructing \Gs s where $\tau_0$, $\tau_2$, and $\tau_3$ all vanish,
since we can then always make $\tau_1$ vanish as well, provided the
manifold is topologically able to admit a \Gs\ with full holonomy \G.

\section{General flows of \Gs s} \label{generalflowssec}

In this section we derive the evolution equations for a
general flow $\ddt \ph$ of a \Gs\ $\ph$. Let $X$ be a vector
field and $h$ a symmetric $2$-tensor on $M$. Then a general
variation of the \Gs\ $\ph$ can be written as
\begin{equation*}
\ddt \ph = \frac{1}{2} \, h^l_i \ph_{ljk} \, \dx{i} \wedge \dx{j}
\wedge \dx{k} + \frac{1}{6} \, X^l \ps_{lijk} \, \dx{i} \wedge
\dx{j} \wedge \dx{k}
\end{equation*}
We write the left hand side in coordinates, and skew-symmetrize
the first term on the right hand side to obtain
\begin{eqnarray*}
\frac{1}{6} \, \ddt \ph_{ijk} \, \dx{i} \wedge \dx{j} \wedge
\dx{k} = & & \frac{1}{6} \left( h^l_i \ph_{ljk} - h^l_j \ph_{lik} -
h^l_k \ph_{lji} \right) \, \dx{i} \wedge \dx{j} \wedge \dx{k} \\ &
& {} + \frac{1}{6} \, X^l \ps_{lijk} \, \dx{i} \wedge
\dx{j} \wedge \dx{k}
\end{eqnarray*}
We can now equate the (totally skew-symmetric) coefficients on
both sides to obtain the {\em general flow equation}:
\begin{equation} \label{generalfloweq}
\boxed{\, \, \ddt \ph_{ijk} \, = \,  h^l_i \ph_{ljk} + h^l_j
\ph_{ilk} + h^l_k \ph_{ijl} + X^l \ps_{lijk} \, \, }
\end{equation}

\subsection{Evolution of the metric $g_{ij}$ and related objects}
\label{metricevolutionsec}

We now proceed to derive the evolution equations for the metric
$g$ and objects related to the metric, specifically the volume form
$\vol$ and the Christoffel symbols $\Gamma^k_{ij}$.
The first step is to compute the evolution of the tensor $B_{ij}$,
which was defined as
\begin{equation*}
B_{ij} \, \dx{1} \wedge \ldots \wedge \dx{7} =  (\ddx{i} \hk
\ph) \wedge (\ddx{j} \hk \ph) \wedge \ph
\end{equation*}
which in coordinates becomes
\begin{eqnarray*}
B_{ij} \, \dx{1} \wedge \ldots \wedge \dx{7} & = & \frac{1}{24}
\, \ph_{ik_1k_2} \ph_{jk_3k_4} \ph_{k_5k_6k_7} \, \dx{k_1} \wedge
\ldots \wedge \dx{k_7} \\ & = & \frac{1}{24} \, \sum_{\sigma \in
S_7} \sgn(\sigma) \, \ph_{i \sigma(1) \sigma(2)} \ph_{j
\sigma(3) \sigma(4)} \ph_{\sigma(5) \sigma(6) \sigma(7)} \, \dx{1}
\wedge \ldots \wedge \dx{7}
\end{eqnarray*}
and hence
\begin{equation} \label{Bijdeteq}
B_{ij} = \frac{1}{24} \, \sum_{\sigma \in S_7} \sgn(\sigma)
\, \ph_{i \sigma(1) \sigma(2)} \ph_{j \sigma(3) \sigma(4)}
\ph_{\sigma(5) \sigma(6) \sigma(7)}
\end{equation}
where the sum is taken over all permutations $\sigma$ of the group
$S_7$ of seven letters.

\smallskip

\begin{thm} \label{Bijevolutionthm}
The evolution of $B_{ij}$ under the flow~\eqref{generalfloweq} is
given by
\begin{equation} \label{Bijevolutioneq}
\ddt B_{ij} = \tr_g(h) B_{ij} + h^l_i B_{lj} + h^l_j B_{il}
\end{equation}
Note that the evolution of $B_{ij}$ deponds only on the
symmetric $2$-tensor $h_{ij}$ and not on the vector field $X^k$.
\end{thm}
\begin{proof}
We will need to appeal to various identities for \Gs s that are
proved in Section~\ref{g2identitiessec}. We start by computing
\begin{eqnarray*}
\ddt B_{ij} & = & \frac{1}{24} \, \sum_{\sigma \in S_7} \sgn(\sigma)
\, \left( \ddt \ph_{i \sigma(1) \sigma(2)} \right) \ph_{j \sigma(3)
\sigma(4)} \ph_{\sigma(5) \sigma(6) \sigma(7)} \\ & & {} + 
\frac{1}{24} \, \sum_{\sigma \in S_7} \sgn(\sigma)
\, \ph_{i \sigma(1) \sigma(2)} \left( \ddt \ph_{j \sigma(3)
\sigma(4)} \right) \ph_{\sigma(5) \sigma(6) \sigma(7)} \\ & & {} +
\frac{1}{24} \, \sum_{\sigma \in S_7} \sgn(\sigma)
\, \ph_{i \sigma(1) \sigma(2)} \ph_{j \sigma(3) \sigma(4)}
\left( \ddt \ph_{\sigma(5) \sigma(6) \sigma(7)} \right)
\end{eqnarray*}
We substitute the evolution equation~\eqref{generalfloweq} into
this expression:
\begin{eqnarray}
\nonumber 24 \ddt B_{ij} & = & \sum_{\sigma \in S_7}
\sgn(\sigma) \, \left( h^l_i \ph_{l\sigma(1)\sigma(2)} +
h^l_{\sigma(1)} \ph_{il\sigma(2)} + h^l_{\sigma(2)}
\ph_{i\sigma(1)l} \right) \ph_{j \sigma(3) \sigma(4)}
\ph_{\sigma(5) \sigma(6) \sigma(7)} \\ \nonumber & & {} +
\sum_{\sigma \in S_7} \sgn(\sigma) \, \ph_{i
\sigma(1) \sigma(2)} \left( h^l_j \ph_{l\sigma(3)\sigma(4)} +
h^l_{\sigma(3)} \ph_{jl\sigma(4)} + h^l_{\sigma(4)}
\ph_{j\sigma(3)l} \right) \ph_{\sigma(5) \sigma(6)
\sigma(7)} \\ \nonumber & & {} + \sum_{\sigma \in S_7}
\sgn(\sigma) \, \ph_{i \sigma(1) \sigma(2)}
\ph_{j\sigma(3)\sigma(4)} \left( h^l_{\sigma(5)}
\ph_{l\sigma(6)\sigma(7)} + h^l_{\sigma(6)}
\ph_{\sigma(5)l\sigma(7)} + h^l_{\sigma(7)}
\ph_{\sigma(5)\sigma(6)l} \right) \\ \nonumber & & {}+ \sum_{\sigma
\in S_7} \sgn(\sigma) \, X^l \ps_{li\sigma(1)\sigma(2)}
\ph_{j \sigma(3) \sigma(4)} \ph_{\sigma(5) \sigma(6) \sigma(7)} \\
\nonumber & & {}+ \sum_{\sigma \in S_7} \sgn(\sigma) \ph_{i
\sigma(1) \sigma(2)} X^l \ps_{lj\sigma(3)\sigma(4)} \ph_{\sigma(5)
\sigma(6) \sigma(7)} \\ \label{Bijtempeq} & & {}+ \sum_{\sigma \in
S_7} \sgn(\sigma) \, \ph_{i \sigma(1) \sigma(2)}
\ph_{j\sigma(3)\sigma(4)} X^l \ps_{l\sigma(5)\sigma(6)\sigma(7)}
\end{eqnarray}
where we have separated the expression into three terms involving
$h_{ij}$ and three terms involving $X^l$.

Let us consider the $X^l$ terms first. Equation~\eqref{vfzeroeq2}
tells us immediately that the last term is zero. Let $\tau$ be the
permutation in $S_7$ which sends $(1,2,3,4,5,6,7) \mapsto
(3,4,1,2,5,6,7)$. This permutation is even, so $\sgn(\sigma
\circ \tau) = \sgn(\sigma)$. Therefore if we write $\sigma' =
\sigma \circ \tau$, then the first two terms involving $X^l$ can be
written as
\begin{eqnarray*}
& & \sum_{\sigma \in S_7} \sgn(\sigma) \, X^l
\ps_{li\sigma(1)\sigma(2)} \ph_{j \sigma(3) \sigma(4)} \ph_{\sigma(5)
\sigma(6) \sigma(7)} \\ & & {}+ \sum_{\sigma' \in S_7}
\sgn(\sigma') \ph_{i \sigma'(3) \sigma'(4)} X^l
\ps_{lj\sigma'(1)\sigma'(2)} \ph_{\sigma'(5) \sigma'(6) \sigma'(7)}
\end{eqnarray*}
which also vanishes because of~\eqref{vfzeroeq}. Hence all the $X^l$
terms are zero and thus the vector field $X^l$ does not affect the
evolution of $B_{ij}$. Therefore it does not affect the
evolution of the metric $g_{ij}$ or the volume form $\vol$ either.
It is well known that infinitesmal variations in the $\wths$
direction do not affect the metric. See~\cite{BS, K1} for other
demonstrations of this fact.

We now return to the terms in~\eqref{Bijtempeq} involving $h_{ij}$.
We start with the first term on the first line, together with the
first term on the second line. These are:
\begin{eqnarray*}
& & \sum_{\sigma \in S_7} \sgn(\sigma) \, h^l_i
\ph_{l\sigma(1)\sigma(2)} \ph_{j \sigma(3) \sigma(4)} \ph_{\sigma(5)
\sigma(6) \sigma(7)} \\ & & {} + \sum_{\sigma \in S_7}
\sgn(\sigma) \, \ph_{i \sigma(1) \sigma(2)} h^l_j
\ph_{l\sigma(3)\sigma(4)} \ph_{\sigma(5) \sigma(6)
\sigma(7)}
\end{eqnarray*}
which, using~\eqref{Bijdeteq}, are seen to be equal to
\begin{equation*}
24 \left( h^l_i B_{lj} + h^l_j B_{il} \right)
\end{equation*}

Next, we consider in~\eqref{Bijtempeq} the remaining two terms from
the first line, together with the remaining two terms from the
second line. These are:
\begin{eqnarray*}
& & \sum_{\sigma \in S_7} \sgn(\sigma) \, \left(
h^l_{\sigma(1)} \ph_{il\sigma(2)} + h^l_{\sigma(2)} \ph_{i\sigma(1)l}
\right) \ph_{j \sigma(3) \sigma(4)} \ph_{\sigma(5) \sigma(6)
\sigma(7)} \\ & & {} + \sum_{\sigma \in S_7} \sgn(\sigma) \,
\ph_{i \sigma(1) \sigma(2)} \left( h^l_{\sigma(3)} \ph_{jl\sigma(4)} +
h^l_{\sigma(4)} \ph_{j\sigma(3)l} \right) \ph_{\sigma(5) \sigma(6)
\sigma(7)}
\end{eqnarray*}
If we consider the permutation $\tau$ which interchanges $1$ and $2$,
and denote $\sigma' = \sigma \circ \tau$, then the first line can be
rewritten as
\begin{eqnarray*}
& & \sum_{\sigma \in S_7} \sgn(\sigma) \, h^l_{\sigma(1)}
\ph_{il\sigma(2)} \ph_{j \sigma(3) \sigma(4)} \ph_{\sigma(5) \sigma(6)
\sigma(7)} \\ & & {} - \sum_{\sigma' \in S_7} \sgn(\sigma') \,
h^l_{\sigma'(1)} \ph_{i\sigma'(2)l} \ph_{j \sigma'(3) \sigma'(4)}
\ph_{\sigma'(5) \sigma'(6) \sigma'(7)} \\ & = & 2 \sum_{\sigma \in S_7}
\sgn(\sigma) \, h^l_{\sigma(1)} \ph_{il\sigma(2)} \ph_{j
\sigma(3) \sigma(4)} \ph_{\sigma(5) \sigma(6) \sigma(7)}
\end{eqnarray*}
by the skew-symmetry of $\ph_{ijk}$. The two terms on the second line
can be similarly combined as
\begin{eqnarray*}
& & 2 \sum_{\sigma \in S_7} \sgn(\sigma) \, h^l_{\sigma(3)}
\ph_{jl\sigma(4)} \ph_{i \sigma(1) \sigma(2)} \ph_{\sigma(5) \sigma(6)
\sigma(7)} \\ & = & 2 \sum_{\sigma \in S_7} \sgn(\sigma) \,
h^l_{\sigma(1)} \ph_{jl\sigma(2)} \ph_{i \sigma(3) \sigma(4)}
\ph_{\sigma(5) \sigma(6) \sigma(7)}
\end{eqnarray*}
by interchanging $1$ with $2$ and $3$ with $4$. We can now
apply~\eqref{wtwscubedcoreq} to these two expressions (with $\alpha_k
= h^l_k$). These terms then become
\begin{eqnarray*}
& & 2 \left( -4(B_{lj} h^l_i - B_{ji} h^l_l) + 24 (\ps_{iljm} g^{km}
h^l_k \sqrt{\det(g)} ) \right) \\ & & {} + 2 \left( -4(B_{li} h^l_j -
B_{ij} h^l_l) + 24 (\ps_{jlim} g^{km} h^l_k \sqrt{\det(g)} ) \right)
\\ & = & 16 \tr_g(h) B_{ij} - 8 \, h^l_i B_{lj} - 8 \, h^l_j B_{il} 
\end{eqnarray*}
using the symmetry of $B_{ij}$ and the skew-symmetry of $\ps_{ijkl}$.

All that remains now in~\eqref{Bijtempeq} are the final three terms
involving $h_{ij}$. These are:
\begin{equation*}
\sum_{\sigma \in S_7} \sgn(\sigma) \, \ph_{i \sigma(1)
\sigma(2)} \ph_{j\sigma(3)\sigma(4)} \left( h^l_{\sigma(5)}
\ph_{l\sigma(6)\sigma(7)} + h^l_{\sigma(6)} \ph_{\sigma(5)l\sigma(7)}
+ h^l_{\sigma(7)} \ph_{\sigma(5)\sigma(6)l} \right)
\end{equation*}
If we consider the cyclic permutation $5 \mapsto 6 \mapsto 7 \mapsto
5$, which is even, we see that the three terms above are actually the
same, so we can combine them into
\begin{equation*}
3 \sum_{\sigma \in S_7} \sgn(\sigma) \, \ph_{i \sigma(1)
\sigma(2)} \ph_{j\sigma(3)\sigma(4)} \ph_{l\sigma(5)\sigma(6)}
h^l_{\sigma(7)}
\end{equation*}
Now we can use~\eqref{wtwscubedeq}, (again with $\alpha_k = h^l_k$),
to write this expression as
\begin{equation*}
3 \, \left( \frac{8}{3} \right) \left( B_{ij} h^l_l + B_{il} h^l_j +
B_{jl} h^l_i \right) = 8 \tr_g(h) B_{ij} + 8 \, h^l_i B_{lj} + 8 \, h^l_j
B_{il}
\end{equation*}

Finally, we combine the results of these calculations of the
terms in~\eqref{Bijtempeq} to obtain
\begin{eqnarray*}
24 \ddt B_{ij} & = & 24 \left( h^l_i B_{lj} + h^l_j B_{il} \right) +
\left( 16 \tr_g(h) B_{ij} - 8 \, h^l_i B_{lj} - 8 \, h^l_j B_{il} \right)
\\ & & {} + \left( 8 \tr_g(h) B_{ij} + 8 \, h^l_i B_{lj} + 8 \, h^l_j
B_{il} \right) \\ & = & 24 \left( \tr_g(h) B_{ij} + h^l_i B_{lj} +
h^l_j B_{il} \right)
\end{eqnarray*}
and the proof is complete.
\end{proof}

\begin{cor} \label{metricevolutioncor}
The evolution of the metric $g_{ij}$ under the
flow~\eqref{generalfloweq} is given by
\begin{equation} \label{metricevolutioneq}
\boxed{\, \, \ddt g_{ij} = 2 h_{ij} \, \, }
\end{equation}
\end{cor}
\begin{proof}
We begin by computing the evolution of $\det(B)$ using
Theorem~\ref{Bijevolutionthm} and Lemma~\ref{volumeflowlemma} (which
applies to any symmetric matrix $g_{ij}$ with inverse $g^{ij}$.)
\begin{eqnarray*}
\ddt \det(B) & = & \left( \ddt B_{ij} \right) B^{ij} \det(B) \\ & = &
\left( \tr_g(h) B_{ij} + h^l_i B_{lj} + h^l_j B_{il} \right) B^{ij} \det(B)
\\ & = & \left( \tr_g(h) \delta^i_i + h^l_i \delta^i_l +
h^l_j \delta^j_l \right) \det(B) \\ & = & 9 \tr_g(h) \det(B)
\end{eqnarray*}
We now use this result to differentiate~\eqref{gijeq}:
\begin{eqnarray*}
\ddt g_{ij} & = & \ddt \left( \frac{1}{6^{\frac{2}{9}}} \,
\frac{B_{ij}}{\det(B)^{\frac{1}{9}}} \right) \\ & = &
\frac{1}{6^{\frac{2}{9}}} \, \left( \frac{\ddt B_{ij}}{
{\det(B)^{\frac{1}{9}}} } - \frac{1}{9} \frac{B_{ij} \ddt \det(B)} {
{\det(B)^{\frac{10}{9}}} } \right) \\ & = & \frac{1}{6^{\frac{2}{9}}}
\, \left( \frac{ \tr_g(h) B_{ij} + h^l_i B_{lj} + h^l_j B_{il} }{
{\det(B)^{\frac{1}{9}}} } - \frac{1}{9} \frac{B_{ij} 9 \tr_g(h) \det(B)} {
{\det(B)^{\frac{10}{9}}} } \right) \\ & = & \tr_g(h)\, g_{ij} + h^l_i g_{lj} +
h^l_j g_{il} - \tr_g(h)\, g_{ij} \, \, = \, \, 2 h_{ij}
\end{eqnarray*}
as claimed.
\end{proof}

\begin{cor} \label{invmetricandvolevolutioncor}
The evolution of the inverse $g^{ij}$ of the metric and the evolution
of the volume form $\vol$ under the flow~\eqref{generalfloweq} are
given by
\begin{eqnarray} \label{invmetricevolutioneq}
& & \boxed{\, \, \ddt g^{ij} = -2 h^{ij} \, \, } \\
\label{volevolutioneq} & & \boxed{\, \, \ddt \vol = \tr_g(h)\, \vol \, \, }
\end{eqnarray}
\end{cor}
\begin{proof}
These follows directly from Lemmas~\ref{invmetriclemma}
and~\ref{volumeflowlemma}.
\end{proof}

\smallskip

Finally we consider the evolution of the Christoffel symbols $\Gamma^k_{ij}$.

\smallskip

\begin{prop} \label{christoffelevolutionprop}
The evolution of the Christoffel symbols $\Gamma^k_{ij}$ 
under the flow~\eqref{generalfloweq} is given by
\begin{equation} \label{christoffelevolutioneq}
\ddt \Gamma^k_{ij} = g^{kl} \left( \nab{i} h_{jl} +
\nab{j} h_{il} - \nab{l} h_{ij} \right)
\end{equation}
\end{prop}
\begin{proof}
This result is standard and can be found, for example, in~\cite{CK}, Lemma 3.2.
(Note that by our Corollary~\ref{metricevolutioncor} above, their $h_{ij}$
in~\cite{CK} is replaced by our $2 \, h_{ij}$.)
\end{proof}

\subsection{Evolution of the dual $4$-form $\ps_{ijkl}$}
\label{psievolutionsec}

We proceed now to the computation of the evolution equation for
the dual $4$-form $\ps = \st \ph$. To do this directly would be
quite complicated, because we would need to compute the evolution
of the Hodge star operator $\st$. We can avoid this by using the
final equation from Lemma~\ref{g2identities1lemma}, which can be
rearranged into the form
\begin{equation} \label{psitempeq}
\ps_{ijkl} = g_{ik} g_{jl} - g_{il} g_{jk} - \ph_{ij\alpha}
\ph_{kl\beta} g^{\alpha\beta}
\end{equation}
This essentially says that, up to some correction terms involving
only the metric, the form $\ps$ is the contraction of $\ph$ with
itself.

\smallskip

\begin{thm} \label{psievolutionthm}
The evolution of the $4$-form $\ps_{ijkl}$ under the
flow~\eqref{generalfloweq} is given by
\begin{equation} \label{psievolutioneq}
\fbox{
  \begin{minipage}{0.6\textwidth}
    \leftline{$\, \, \ddt \ps_{ijkl} = h^m_i \ps_{mjkl} + h^m_j
\ps_{imkl} + h^m_k \ps_{ijml} + h^m_l \ps_{ijkm}$}
\vskip 0.05in
    \leftline{$\qquad \qquad \, \, {}- X_i \ph_{jkl} + X_j
\ph_{ikl} - X_k \ph_{ijl} + X_l \ph_{ijk}$}
  \end{minipage}
}
\end{equation}
where $X_k = g_{kl} X^l$.
\end{thm}
\begin{proof}
We differentiate~\eqref{psitempeq} with respect to time:
\begin{eqnarray} \label{psitemp2eq}
\ddt \ps_{ijkl} & = & \left( \ddt g_{ik} \right) g_{jl} + g_{ik}
\left( \ddt g_{jl} \right) - \left( \ddt g_{il} \right) g_{jk} -
g_{il} \left( \ddt g_{jk} \right) \\ \nonumber & & {} -
\ph_{ij\alpha} \ph_{kl\beta} \left( \ddt g^{\alpha\beta} \right)-
\left( \ddt \ph_{ij\alpha} \right) \ph_{kl\beta} g^{\alpha \beta} -
\ph_{ij\alpha} \left( \ddt \ph_{kl\beta} \right) g^{\alpha\beta}
\end{eqnarray}
Using Corollary~\ref{metricevolutioncor} and
equation~\eqref{invmetricevolutioneq}, the first five terms
in~\eqref{psitemp2eq} are
\begin{equation*}
2 \, h_{ik} g_{jl} + 2 \, h_{jl} g_{ik} - 2 \, h_{il} g_{jk} - 2
\, h_{jk} g_{il} + 2 \, \ph_{ij\alpha} \ph_{kl\beta} h^{\alpha
\beta}
\end{equation*}
Meanwhile the next to last term in~\eqref{psitemp2eq},
using~\eqref{generalfloweq} and Lemma~\ref{g2identities1lemma} is
\begin{eqnarray*}
& & {}- \left( h^m_i \ph_{mj\alpha} + h^m_j \ph_{im\alpha} +
h^m_{\alpha} \ph_{ijm} + X^m \ps_{mij\alpha} \right) \ph_{kl\beta}
g^{\alpha \beta} \\ & = & - h^m_i \left( g_{mk} g_{jl} - g_{ml}
g_{jk} - \ps_{mjkl} \right) - h^m_j \left( g_{ik} g_{ml} - g_{il}
g_{km} - \ps_{imkl} \right) \\ & & {}- h^{m\beta} \ph_{ijm}
\ph_{kl\beta} - X^m \ps_{mij\alpha} \ph_{kl\beta}
g^{\alpha \beta} \\ & = & - h_{ik} g_{jl} + h_{il} g_{jk} - 
h_{jl} g_{ik} + h_{jk} g_{il} + h^m_i \ps_{mjkl} + h^m_j
\ps_{imkl} \\ & & {}  {}- h^{m\beta} \ph_{ijm}
\ph_{kl\beta} - X^m \ps_{mij\alpha} \ph_{kl\beta}
g^{\alpha \beta}
\end{eqnarray*}
Finally, the last term in~\eqref{psitemp2eq} is the same as the
next to last term if we interchange $i$ with $k$ and $j$ with $l$,
so it equals
\begin{eqnarray*}
& & - h_{ik} g_{jl} + h_{jk} g_{il} - h_{jl} g_{ik} + h_{il} g_{jk} +
h^m_k \ps_{mlij} + h^m_l \ps_{kmij} \\ & & {} {}- h^{m\beta} \ph_{klm}
\ph_{ij\beta} - X^m \ps_{mkl\alpha} \ph_{ij\beta} g^{\alpha \beta}
\end{eqnarray*}
When we add up these three sets of terms, all the `$g \cdot h$'
terms and the `$\ph \cdot \ph$' terms cancel, and we are left with
\begin{eqnarray} \nonumber
\ddt \ps_{ijkl} & = & h^m_i \ps_{mjkl} + h^m_j \ps_{imkl} + h^m_k
\ps_{ijml} + h^m_l \ps_{ijkm} \\ \label{psievolutiontempeq} & & {}- X^m
\ps_{mij\alpha} \ph_{kl\beta} g^{\alpha\beta} - X^m \ps_{mkl\alpha}
\ph_{ij\beta} g^{\alpha \beta}
\end{eqnarray}
We deal with the $X^m$ terms above using the final equation of
Lemma~\ref{g2identities2lemma}. They become
\begin{eqnarray*}
& & -X^m \left( g_{km} \ph_{lij} + g_{ki} \ph_{mlj} + g_{kj}
\ph_{mil} - g_{lm} \ph_{kij} - g_{li} \ph_{mkj} - g_{lj} \ph_{mik}
\right) \\ & & {}- X^m \left( g_{im} \ph_{jkl} + g_{ik} \ph_{mjl} +
g_{il} \ph_{mkj} - g_{jm} \ph_{ikl} - g_{jk} \ph_{mil} - g_{jl}
\ph_{mki} \right) \\ & = & -X_k \ph_{lij} + X_l \ph_{kij} - X_i
\ph_{jkl} + X_j \ph_{ikl}
\end{eqnarray*}
because the other terms cancel in pairs. Substituting this expression
into~\eqref{psievolutiontempeq} above completes the proof.
\end{proof}

\smallskip

\begin{rmk} \label{expectedrmk}
Let the infinitesmal deformation of a \G\ $3$-form $\ph(t)$
be given by
\begin{equation*}
\ddt \ph(t) = \eta_1 + \eta_7 + \eta_{27}
\end{equation*}
where $\eta_k$ belongs to the subspace $\Omega^3_k$ associated to
$\ph(t)$. Then one can show that the infinitesmal deformation of the
associated $4$-form $\ps(t)$ is
\begin{equation*}
\ddt \ps(t) = \frac{4}{3} \left( \st_{\ph(t)} \eta_1 \right) + \left(
\st_{\ph(t)} \eta_7 \right) - \left( \st_{\ph(t)} \eta_{27} \right)
\end{equation*}
This is mentioned in~\cite{Br3, J1, J4} and an explicit proof is given
in~\cite{Hi1}. The purpose of this remark is to clarify that
Theorem~\ref{psievolutionthm} agrees with this result. Let $\ddt \ph =
\eta$ be an arbitrary $3$-form. It can be written uniquely as
\begin{eqnarray*}
\eta & = & \frac{1}{2} h^l_l \ph_{ljk} \, \dx{i} \wedge \dx{j} \wedge
\dx{k} + \frac{1}{6} X^l \ps_{lijk} \, \dx{i} \wedge \dx{j} \wedge \dx{k}
\\ & = & h^l_i \, \dx{i} \wedge \left( \ddx{l} \hk \ph \right) + (X \hk \ps)
\end{eqnarray*}
for some vector field $X^l$ and some symmetric $2$-tensor $h_{ij}$. We
can write $h_{ij} = \frac{1}{7} \tr_g(h)\, g_{ij} + h^0_{ij}$ where
$h^0_{ij}$ is the trace free part, and $\tr_g(h) = g^{ij} h_{ij}$.
Taking the Hodge star, and
using~\eqref{iosrelationseq} and Proposition~\ref{symmstarprop} gives
\begin{equation*}
\st \eta = f^m_k \dx{k} \wedge \left( \ddx{m} \hk \ps \right) - \Xs \wedge \ph
\end{equation*}
where $f_{ij} = \frac{1}{4} \tr_g(h)\, g_{ij} - h_{ij} = \frac{1}{4} \tr_g(h)\,
g_{ij} - \frac{1}{7} \tr_g(h)\, g_{ij} - h^0_{ij} = \frac{3}{4}
\left( \frac{1}{7} \tr_g(h)\, g_{ij} \right) - h^0_{ij}$.
In coordinates, we have
\begin{eqnarray*}
\st \eta & =  & f^m_k \, \dx{k} \wedge \left( \ddx{m} \hk \ps \right) - \left( \Xs
\wedge \ph \right) \\ & & = \frac{1}{6} f^m_k \ps_{mjkl} \, \dx{i}
\wedge \dx{j} \wedge \dx{k} \wedge \dx{l} - \frac{1}{6} X_i \ph_{jkl}
\, \dx{i} \wedge \dx{j} \wedge \dx{k} \wedge \dx{l} \\ & & =
\frac{1}{24} \left( f^m_i \ps_{mjkl} + f^m_j \ps_{imkl} + f^m_k
\ps_{ijml} + f^m_l \ps_{ijkm} \right) \, \dx{i} \wedge \dx{j} \wedge
\dx{k} \wedge \dx{l} \\ & & \quad {} + \frac{1}{24} \left( -X_i
\ph_{jkl} + X_j \ph_{ikl} - X_k \ph_{ijl} + X_l \ph_{ijk} \right) \,
\dx{i} \wedge \dx{j} \wedge \dx{k} \wedge \dx{l}
\end{eqnarray*}
Comparing with~\eqref{psievolutioneq}, we see that
\begin{equation*}
\left(\ddt \ps\right)_1 = \frac{4}{3} \st \left(\ddt \ph\right)_1
\quad ; \quad \left(\ddt \ps\right)_7 = \st \left(\ddt \ph\right)_7
\quad ; \quad \left(\ddt \ps\right)_{27} = - \st \left(\ddt
\ph\right)_{27}
\end{equation*}
as expected. The approach we have adopted is advantageous, because in
this setup
\begin{eqnarray*}
& \text{If } \qquad & \ddt \ph \, \, = \, \, h^k_m \, \dx{k} \wedge
\left( \ddx{m} \hk \ph \right) + (X \hk \ps) \\ & \text{then } \qquad
& \ddt \ps \, \, = \, \, h^k_m \, \dx{k} \wedge \left( \ddx{m} \hk \ps
\right) + \st (X \hk \ps)
\end{eqnarray*}
so the two equations look much more symmetric this way.
\end{rmk}

\subsection{Evolution of the torsion forms}
\label{torsionevolutionsec}

In this section we derive the evolution equations for the four torsion
tensors of a \Gs\ under a general flow described by a symmetric
tensor $h_{ij}$ and a vector field $X^k$. We begin with the evolution of
$\nab{l} \ph_{ijk}$.

\smallskip

\begin{lemma} \label{nabphevolutionlemma}
The evolution of $\nab{l} \ph_{ijk}$ under the
flow~\eqref{generalfloweq} is given by
\begin{eqnarray} \nonumber
\ddt \left( \nab{l} \ph_{ijk} \right) & = & h^m_i (\nab{l} \ph_{mjk}) 
+ h^m_j (\nab{l} \ph_{imk}) + h^m_k (\nab{l} \ph_{ijm}) +
X^m (\nab{l} \ps_{mijk}) \\ \nonumber & & {} +
(\nab{s} h_{il}) g^{ms} \ph_{mjk} + (\nab{s} h_{jl}) g^{ms} \ph_{imk} +
(\nab{s} h_{kl}) g^{ms} \ph_{ijm} \\ \nonumber & & {} - 
(\nab{i} h_{ls}) g^{ms} \ph_{mjk} - (\nab{j} h_{ls}) g^{ms} \ph_{imk} -
(\nab{k} h_{ls}) g^{ms} \ph_{ijm} \\ \label{nabphevolutioneq} & & {} +
(\nab{l} X^m) \ps_{mijk}
\end{eqnarray}
in terms of $h_{ij}$ and $X^k$.
\end{lemma}
\begin{proof}
Recall that
\begin{equation*}
\nab{l} \ph_{ijk} = \ddx{l} \ph_{ijk} - \Gamma^m_{li} \ph_{mjk}
- \Gamma^m_{lj} \ph_{imk} - \Gamma^m_{lk} \ph_{ijm}
\end{equation*}
We differentiate this equation with respect to $t$ to obtain
\begin{eqnarray*}
\ddt \left( \nab{l} \ph_{ijk} \right) & = & \ddx{l} \left( \ddt \ph_{ijk}
\right) - \left(\ddt \Gamma^m_{li}\right) \ph_{mjk} - \left( \ddt \Gamma^m_{lj}
\right) \ph_{imk} - \left( \ddt \Gamma^m_{lk} \right) \ph_{ijm} \\ & & {}
- \Gamma^m_{li} \left( \ddt \ph_{mjk} \right) - \Gamma^m_{lj} \left( \ddt
\ph_{imk} \right) - \Gamma^m_{lk} \left( \ddt \ph_{ijm} \right) \\
& = & \nab{l} \left( \ddt \ph_{ijk} \right)
 - \left(\ddt \Gamma^m_{li}\right) \ph_{mjk} - \left( \ddt \Gamma^m_{lj}
\right) \ph_{imk} - \left( \ddt \Gamma^m_{lk} \right) \ph_{ijm} 
\end{eqnarray*}
Now we substitute~\eqref{generalfloweq} and~\eqref{christoffelevolutioneq} to
get
\begin{eqnarray*}
\ddt \left( \nab{l} \ph_{ijk} \right) & = & \nab{l} \left( h^m_i \ph_{mjk}
+ h^m_j \ph_{imk} + h^m_k \ph_{ijm} + X^m \ps_{mijk} \right) \\ & & 
{} - g^{ms} \left( \nab{l} h_{is} + \nab{i} h_{ls} - \nab{s} h_{il} \right)
\ph_{mjk} \\ & & - g^{ms} \left( \nab{l} h_{js} + \nab{j} h_{ls} -
\nab{s} h_{jl} \right) \ph_{imk} \\ & & - g^{ms} \left( \nab{l} h_{ks} +
\nab{k} h_{ls} - \nab{s} h_{kl} \right) \ph_{ijm}  
\end{eqnarray*}
We use the product rule on the first line, and see that all the terms involving
$\nab{l} h$ cancel in pairs. The result now follows.
\end{proof}

\smallskip

We are now ready to compute the evolution equation of the full torsion
tensor $T_{pq}$ for a general flow of \Gs s.

\smallskip

\begin{thm} \label{fulltorsionevolutionthm}
The evolution of the full torsion tensor $T_{pq}$ under the
flow~\eqref{generalfloweq} is given by
\begin{equation} \label{fulltorsionevolutioneq}
\boxed{\, \, \ddt T_{pq} = T_{pl} g^{lm} h_{mq} + T_{pl} g^{lm} X_{mq}
+ (\nab{k} h_{ip}) g^{ka} g^{ib} \ph_{abq} + \nab{p} X_q \, \, }
\end{equation}
where $X_j = g_{jk} X^k$ and $X_{ij} = X^k \ph_{kij}$ is the element of
$\wtws$ corresponding to $X$.
\end{thm}
\begin{proof}
We start with Lemma~\ref{fulltorsionlemma}, and differentiate:
\begin{eqnarray} \nonumber
\ddt T_{pq} & = & \frac{1}{24} \ddt \left( \nab{p} \ph_{ijk} \right) \ps_{qbcd}
g^{ib} g^{jc} g^{kd} + \frac{1}{24} (\nab{p} \ph_{ijk}) \left( \ddt
\ps_{qbcd} \right) g^{ib} g^{jc} g^{kd} \\ \label{fulltorsionevolutiontempeq}
& & {} + \frac{3}{24} (\nab{p} \ph_{ijk}) 
\ps_{qbcd} \left( \ddt g^{ib} \right) g^{jc} g^{kd}
\end{eqnarray}
where we have relabelled indices and used the skew-symmetry of $\ph$ and $\ps$
to combine the three terms involving derivatives of $g$. By
Lemma~\ref{nabphevolutionlemma}, the first term
in~\eqref{fulltorsionevolutiontempeq} is
\begin{eqnarray*}
& & \frac{1}{24} \left( h^m_i (\nab{p} \ph_{mjk}) + h^m_j (\nab{p} \ph_{imk})
+ h^m_k (\nab{p} \ph_{ijm}) + X^m (\nab{p} \ps_{mijk}) \right) \ps_{qbcd}
g^{ib} g^{jc} g^{kd} \\ & & {} + \frac{1}{24} \left(
(\nab{s} h_{ip}) g^{ms} \ph_{mjk} + (\nab{s} h_{jp}) g^{ms} \ph_{imk} +
(\nab{s} h_{kp}) g^{ms} \ph_{ijm} \right) \ps_{qbcd}
g^{ib} g^{jc} g^{kd} \\ & & {} - \frac{1}{24} \left(
(\nab{i} h_{ps}) g^{ms} \ph_{mjk} + (\nab{j} h_{ps}) g^{ms} \ph_{imk} +
(\nab{k} h_{ps}) g^{ms} \ph_{ijm} \right) \ps_{qbcd}
g^{ib} g^{jc} g^{kd} \\ & & {} + \frac{1}{24} (\nab{p} X^m) \ps_{mijk} 
\ps_{qbcd} g^{ib} g^{jc} g^{kd}
\end{eqnarray*}
Again exploiting the skew-symmetry of $\ph$ and $\ps$ and relabelling indices,
the above expression simplifies to
\begin{eqnarray*}
& & \frac{3}{24} h^m_i (\nab{p} \ph_{mjk}) \ps_{qbcd} g^{ib} g^{jc} g^{kd} 
 + \frac{1}{24} X^m (\nab{p} \ps_{mijk}) \ps_{qbcd}
g^{ib} g^{jc} g^{kd} \\ & & {}+ \frac{3}{24} (\nab{s} h_{ip}) g^{ms} \ph_{mjk} 
\ps_{qbcd} g^{ib} g^{jc} g^{kd} - \frac{3}{24}
(\nab{i} h_{ps}) g^{ms} \ph_{mjk} \ps_{qbcd} g^{ib} g^{jc} g^{kd} \\ & & {}+
\frac{1}{24} (\nab{p} X^m) \ps_{mijk} \ps_{qbcd} g^{ib} g^{jc} g^{kd}
\end{eqnarray*}
Now we use Proposition~\ref{nabphirelnabpsiprop} on the second term,
Lemma~\ref{g2identities2lemma} on the third and fourth terms, and
Lemma~\ref{g2identities3lemma} on the fifth term above to obtain
\begin{eqnarray*}
& & \frac{3}{24} (\nab{p} \ph_{mjk}) \ps_{qbcd} h^{mb} g^{jc} g^{kd} + 
\frac{3}{24} X^m (\nab{p} \ph_{mjk}) \ph_{qcd} g^{jc} g^{kd} \\ & & {}
+ \frac{3}{24} (\nab{s} h_{ip}) g^{ms} (-4 \ph_{mqb}) g^{ib}
 - \frac{3}{24} (\nab{i} h_{ps}) g^{ms} (-4 \ph_{mqb}) g^{ib} \\ & & {}+
\frac{1}{24} (\nab{p} X^m) (24 \, g_{mq})
\end{eqnarray*}
Finally using the symmetry of $h$ and skew-symmetry of $\ph$, the third and
fourth terms above combine, and hence we have written the first term
of~\eqref{fulltorsionevolutiontempeq} as
\begin{equation} \label{fulltorsionevolutiontempeq2}
\frac{3}{24} (\nab{p} \ph_{mjk}) \ps_{qbcd} h^{mb} g^{jc} g^{kd} + 
\frac{3}{24} X^m (\nab{p} \ph_{mjk}) \ph_{qcd} g^{jc} g^{kd} 
+ (\nab{i} h_{ps}) g^{ms} \ph_{mqb} g^{ib} + \nab{p} X_q
\end{equation}
We now consider the second term of~\eqref{fulltorsionevolutiontempeq}.
Using Theorem~\ref{psievolutionthm}. It is
\begin{eqnarray*}
& & \frac{1}{24} (\nab{p} \ph_{ijk}) \left( h^m_q \ps_{mbcd} + h^m_b
\ps_{qmcd} + h^m_c \ps_{qbmd} + h^m_d \ps_{qbcm} \right) g^{ib} g^{jc} g^{kd}
\\ & & + \frac{1}{24} (\nab{p} \ph_{ijk}) \left( - X_q \ph_{bcd} + X_b
\ph_{qcd} - X_c \ph_{qbd} + X_d \ph_{qbc} \right) g^{ib} g^{jc} g^{kd}
\end{eqnarray*}
As usual, by exploiting symmetries, this simplifies to
\begin{eqnarray*}
& & \frac{1}{24} (\nab{p} \ph_{ijk}) h^m_q \ps_{mbcd} g^{ib} g^{jc} g^{kd}
+ \frac{3}{24} (\nab{p} \ph_{ijk}) h^m_b \ps_{qmcd} g^{ib} g^{jc} g^{kd}
\\ & & - \frac{1}{24} X_q (\nab{p} \ph_{ijk}) \ph_{bcd} g^{ib} g^{jc} g^{kd}
+ \frac{3}{24} X_b (\nab{p} \ph_{ijk}) \ph_{qcd} g^{ib} g^{jc} g^{kd}
\end{eqnarray*}
The third term above vanishes by Proposition~\ref{g2derivativeidentitiesprop}.
Therefore we see that the second term of~\eqref{fulltorsionevolutiontempeq}
can be written as
\begin{equation} \label{fulltorsionevolutiontempeq3}
T_{pm} h^m_q + \frac{3}{24} (\nab{p} \ph_{ijk}) \ps_{qmcd} h^{im} g^{jc} g^{kd}
+ \frac{3}{24} X^i (\nab{p} \ph_{ijk}) \ph_{qcd} g^{jc} g^{kd}
\end{equation}
Finally, using equation~\eqref{invmetricevolutioneq}, the third term
of~\eqref{fulltorsionevolutiontempeq} is
\begin{equation} \label{fulltorsionevolutiontempeq4}
-\frac{6}{24} (\nab{p} \ph_{ijk}) \ps_{qbcd} h^{ib} g^{jc} g^{kd}
\end{equation}
Adding the expressions~\eqref{fulltorsionevolutiontempeq2},
\eqref{fulltorsionevolutiontempeq3},
and~\eqref{fulltorsionevolutiontempeq4} gives
\begin{equation} \label{fulltorsionevolutiontempeq5}
\ddt T_{pq} = T_{pm} h^m_q + \frac{1}{4} X^i (\nab{p} \ph_{ijk}) \ph_{qcd}
g^{jc} g^{kd} + (\nab{i} h_{ps}) g^{ms} \ph_{mqb} g^{ib} + \nab{p} X_q
\end{equation}
All that remains now is to rewrite the second term in the above expression.
Using $\nab{p} \ph_{ijk} = T_{pn} g^{nm} \ps_{mijk}$ and
Lemma~\ref{g2identities2lemma}, we see
\begin{eqnarray*}
\frac{1}{4} X^i (\nab{p} \ph_{ijk}) \ph_{qcd} g^{jc} g^{kd} & = & \frac{1}{4}
X^i T_{pn} g^{nm} \ps_{mijk} \ph_{qcd} g^{jc} g^{kd} \\ & = & \frac{1}{4}
X^i T_{pn} g^{nm} (-4 \ph_{qmi}) = T_{pn} g^{nm} X_{mq}
\end{eqnarray*}
Substituting this into~\eqref{fulltorsionevolutiontempeq5} and relabelling some
indices gives
\begin{equation*}
\ddt T_{pq} = T_{pl} g^{lm} h_{mq} + T_{pl} g^{lm} X_{mq}
+ (\nab{k} h_{ip}) g^{ka} g^{ib} \ph_{abq} + \nab{p} X_q
\end{equation*}
which is what we wanted to prove.
\end{proof}

\smallskip

We can now derive the evolution equations of the four independent
torsion forms using Proposition~\ref{torsionformsprop}. We use
$\langle \cdot, \cdot \rangle$ to denote the matrix norm: $\langle A,
B \rangle = A_{ij} B_{kl} g^{ik} g^{jl}$. Note that in the case of
skew-symmetric matrices, this differs from their norm as $2$-forms by
a factor of $2$. Furthermore, $[ A, B]_{pq} = A_{pl} g^{lm} B_{mq} -
B_{pl} g^{lm} A_{mq}$ is the matrix commutator while $\{ A, B\}_{pq} =
A_{pl} g^{lm} B_{mq} + B_{pl} g^{lm} A_{mq}$ is the anti-commutator.

\smallskip

\begin{prop} \label{symmetrictorsionevolutionprop}
The evolution equations of the scalar torsion $\tau_0$ and the symmetric
traceless torsion $\tau_3$ under the flow~\eqref{generalfloweq} are
given by
\begin{equation} \label{tauzeroevolutioneq}
\ddt \tau_0 = -\frac{1}{7} \, \tr_g(h) \tau_0 + \frac{4}{7}\, \langle
h, \tau_3 \rangle - \frac{4}{7} \, \langle X, \tau_1 \rangle + \frac{4}{7}
\, g^{pq} (\nab{p} X_q)
\end{equation}
and
\begin{eqnarray} \nonumber 
\ddt (\tau_3)_{ij} & = & \left( - \frac{1}{28} \tr_g(h) \tau_0 +
\frac{1}{7} \langle h, \tau_3 \rangle - \frac{1}{7} \langle X, \tau_1 
\rangle + \frac{1}{7} \, g^{pq} (\nab{p} X_q) \right) g_{ij} \\
\nonumber & & {}+ \frac{1}{4} \, \tau_0 h_{ij} +
\frac{1}{2} \{h, \tau_3\}_{ij} + \frac{1}{2} \, [h, \tau_1]_{ij} -
\frac{1}{4} \, [h, \tau_2]_{ij} \\ \nonumber & & {}- \frac{1}{2} \, 
[X, \tau_3]_{ij} - \frac{1}{2} \, \{ X, \tau_1 \}_{ij} + \frac{1}{4} \,
\{ X, \tau_2 \}_{ij} \\ \nonumber & & {} -\frac{1}{2}
(\nab{k} h_{li}) g^{ka} g^{lb} \ph_{abj} -\frac{1}{2}
(\nab{k} h_{lj}) g^{ka} g^{lb} \ph_{abi} \\ \label{tauthreeevolutioneq}
& & {} - \frac{1}{2} (\nab{i} X_j + \nab{j} X_i) 
\end{eqnarray}
\end{prop}
\begin{proof}
We have $\tau_0 = \frac{4}{7} \, g^{pq} T_{pq}$. Differentiating and using
Theorem~\ref{fulltorsionevolutionthm} and~\eqref{invmetricevolutioneq},
\begin{eqnarray*}
\ddt \tau_0 & = & \frac{4}{7} \left( \ddt g^{pq} \right) T_{pq}
+ \frac{4}{7} \, g^{pq} \left( \ddt T_{pq} \right) \\ & = & {}-
\frac{8}{7} \, h^{pq} T_{pq} + \frac{4}{7} \, g^{pq} \left( T_{pl}
g^{lm} h_{mq} + T_{pl} g^{lm} X_{mq} + (\nab{k} h_{ip}) g^{ka} g^{ib}
\ph_{abq} + \nab{p} X_q \right) \\ & = & -\frac{4}{7} \, \langle T, h 
\rangle - \frac{4}{7} \langle T, X \rangle + \frac{4}{7} (\nab{k} h^{bq})
\ph_{abq} + \frac{4}{7} \, g^{pq} (\nab{p} X_q)
\end{eqnarray*}
where we have used the symmetry of $h$ and the skew-symmetry of $X$. The
third term vanishes by the skew-symmetry of $\ph$. The result now follows
by recalling that $T_{lm} = \frac{1}{4} \, \tau_0 g_{lm} - (\tau_3)_{lm}
+ (\tau_1)_{lm} - \frac{1}{2} \, (\tau_2)_{lm}$ and that this
decomposition is orthogonal with respect to the matrix inner product.
Equation~\eqref{tauthreeevolutioneq} is proved similarly using
$(\tau_3)_{ij} = \frac{1}{4} \, \tau_0 g_{ij} - \frac{1}{2} (T_{ij} + 
T_{ji})$ and substituting~\eqref{tauzeroevolutioneq} into the
computation. We omit the details.
\end{proof}

\smallskip

Before we can compute the evolution equations for $\tau_1$ and
$\tau_2$, we need the following preliminary result. Define a 
linear map $P$ which takes $2$-tensors to skew-symmetric $2$-tensors
by $(P(A))_{pq} = A_{ij} g^{ia} g^{jb} \ps_{abpq}$. Clearly the symmetric
$2$-tensors are in the kernel of $P$, and by Proposition~\ref{wtwprop},
we have $P(C_7 + C_{14}) = -4 \, C_7 + 2 \, C_{14}$ where 
$C_7 \in \wtws$ and $C_{14} \in \wtwf$.

\smallskip

\begin{lemma} \label{skewtorsionpreliminarylemma}
If $C = C_7 + C_{14}$ is a skew-symmetric tensor, then the evolution of
the skew-symmetric tensor $P(C)$ under the flow~\eqref{generalfloweq} is
given by
\begin{eqnarray} \nonumber
\ddt (P(C))_{ij} & = & ( P ( \ddt C ) )_{ij} + 6 \,
\pi_7(\{h, C_{14}\})_{ij} - 6 \, \pi_{14}(\{h, C_7\})_{ij} \\ 
\label{skewtorsionpreliminaryeq} & & \qquad {}- 2 \, \pi_7 ([X, C_{14}])_{ij}
+ 2 \, \pi_{14}([X, C_7])_{ij}
\end{eqnarray}
where $\pi_7$ and $\pi_{14}$ denote the projections onto $\wtws$ and
$\wtwf$, respectively.
\end{lemma}
\begin{proof}
Using~\eqref{invmetricevolutioneq} and~\eqref{psievolutioneq}, we see
that $\ddt \left( C_{ab} g^{ap} g^{bq} \ps_{pqij} \right)$ equals
\begin{eqnarray} \nonumber
& & \left( \ddt C_{ab} \right) g^{ap} g^{bq} \ps_{pqij} + 2 \, C_{ab}
\left( \ddt g^{ap} \right) g^{bq} \ps_{pqij} + C_{ab} g^{ap} g^{bq}
\left( \ddt \ps_{pqij} \right) \\ \nonumber & = & ( P ( \ddt C ) )_{ij}
- 4 \, C_{ab} h^{ap} g^{bq} \ps_{pqij}  + C_{ab} g^{ap} g^{bq}
( h^l_p \ps_{lqij} + h^l_q \ps_{plij}) \\ \nonumber & & \quad {}+
C_{ab} g^{ap} g^{bq} (h^l_i \ps_{pqlj} + h^l_j \ps_{pqil} - X_p \ph_{qij}
+ X_q \ph_{pij} - X_i \ph_{pqj} + X_j \ph_{pqi} ) \\ \nonumber
& = & ( P ( \ddt C ) )_{ij} - 2 \, C_{ab} h^{ap} g^{bq} \ps_{pqij}
 + h^l_i (P(C))_{lj} + (P(C))_{il} h^l_j \\ & &
\label{skewtorsionpreliminarytempeq} \quad {} +2 
(C_{ab} X^b g^{ap}) \ph_{pij} - 6 \, (C_7)_j X_i + 6 \, (C_7)_i X_j  
\end{eqnarray}
where we have used the skew-symmetry of $C$ and of $\ph$ and relabeled
indices to combine terms. The second term above can be written as
\begin{eqnarray*}
& & -2 \, h_{al} g^{lm} C_{mb} g^{ap} g^{bq} \ps_{pqij} = {}- 
(h_{al} g^{lm} C_{mb} + C_{al} g^{lm} h_{mb} ) g^{ap} g^{bq} \ps_{pqij} 
\\ & = & {}- \{ h, C \}_{ab} g^{ap} g^{bq} \ps_{pqij} = {}-
P ( \{ h, C \} )_{ij} = 4 \, (\pi_7 \{h, C \})_{ij} - 2 \, (\pi_{14} 
\{h, C \})_{ij} \\ & = & 4 \, (\pi_7 \{h, C_7 \})_{ij} + 
4 \, (\pi_7 \{h, C_{14} \})_{ij} - 2 \, (\pi_{14} 
\{h, C_7 \})_{ij} - 2 \, (\pi_{14} \{h, C_{14} \})_{ij}
\end{eqnarray*}
Meanwhile the third and fourth terms
of~\eqref{skewtorsionpreliminarytempeq} become
\begin{eqnarray*}
& & \{ h, P(C) \}_{ij} = \{ h, -4\, C_7 + 2\, C_{14} \}_{ij} \\
& = &  {}-4 \, (\pi_7 \{h, C_7 \})_{ij} + 
2 \, (\pi_7 \{h, C_{14} \})_{ij} - 4 \, (\pi_{14} 
\{h, C_7 \})_{ij} + 2 \, (\pi_{14} \{h, C_{14} \})_{ij}
\end{eqnarray*}
Combining these expressions, after some cancellation we see that
\begin{eqnarray} \nonumber
\ddt \left( C_{ab} g^{ap} g^{bq} \ps_{pqij} \right) & = & 
( P ( \ddt C ) )_{ij} + 6 \, \pi_7(\{h, C_{14}\})_{ij} - 6 \,
\pi_{14}(\{h, C_7\})_{ij} \\ \label{skewtorsionpreliminarytempeq2}
& & \quad +2  (C_{ab} X^b g^{ap}) \ph_{pij} - 6 \, (C_7)_j X_i +
6 \, (C_7)_i X_j
\end{eqnarray}
Consider now the third to last term above. In the notation of
Proposition~\ref{wtwactionvecsprop}, this is $2 \, C(X)_{ij}$, which
by~\eqref{wtwsactioneq} and~\eqref{wtwfactioneq} is
\begin{eqnarray*}
& & 2 ( C_7(X)_{ij} + C_{14}(X)_{ij} ) = 2 ( -\frac{1}{2} \, [C_7, X]_{ij}
 - \frac{3}{2} \, (C_7)_i X_j + \frac{3}{2} \, (C_7)_j X_i + 
[C_{14}, X]_{ij} ) \\ & = & - [C_7, X]_{ij} + 2 \, [C_{14}, X]_{ij}
- 3\, (C_7)_i X_j + 3 \, (C_7)_j X_i 
\end{eqnarray*}
Hence the final three terms of~\eqref{skewtorsionpreliminarytempeq2} are
\begin{eqnarray*}
& & {}- [C_7, X]_{ij} + 2 \, [C_{14}, X]_{ij} + 3\, (C_7)_i X_j - 3 \,
(C_7)_j X_i \\ & = & {}- [C_7, X]_{ij} + 2 \, [C_{14}, X]_{ij}
+ 3 \left( {}-\frac{1}{3} [C_7, X]_{ij} + \frac{2}{3} \,
(C_7 \times X)_{ij} \right) \\ & = & {}- 2 \, [C_7, X]_{ij} + 2 \,
[C_{14}, X]_{ij} + 2\, (C_7 \times X)_{ij}
\end{eqnarray*}
using~\eqref{wtwsactioneq2}. By equation~\eqref{pi7commutatorcrosseq}, the
first and third terms above combine to give $-2 \, \pi_{14} ([C_7, X])_{ij}$.
Also, equation~\eqref{wtwfactioneq} says $[C_{14}, X] = \pi_7 ([C_{14}, X])$.
Equation~\eqref{skewtorsionpreliminaryeq} now follows since the
commutator $[ \cdot, \cdot ]$ is skew-symmetric in its arguments.
\end{proof}

\smallskip

\begin{prop} \label{skewtorsionevolutionprop}
The evolution equations of the vector torsion $\tau_1$ (as a $\wtws$-form) and
the Lie algebra torsion $\tau_2$ under the flow~\eqref{generalfloweq} are
given by
\begin{eqnarray} \nonumber 
\ddt (\tau_1)_{ij} & = & \pi_7 ( \frac{1}{2} \, [h, \tau_3]_{ij} + \frac{1}{2} \,
\{h, \tau_1\}_{ij} + \frac{1}{4} \, \{h, \tau_2\}_{ij} ) +
\pi_{14} ( \{h, \tau_1\}_{ij} ) \\ \nonumber & & {}+ \frac{1}{4} \, \tau_0 X_{ij}
+ \pi_7 ( - \frac{1}{2} \, \{X, \tau_3\}_{ij} - \frac{1}{2} \, [X, \tau_1]_{ij}
+ \frac{1}{12} \, [X, \tau_2]_{ij} ) + \pi_{14} ( - \frac{1}{3} \,
[X, \tau_1]_{ij} ) \\ \nonumber & & {} - \frac{1}{6} \, g^{pq}
(\nab{p} h_{qk}) g^{kl} \ph_{lij} + \frac{1}{6} \, (\nab{k} \tr_g(h))g^{kl}
\ph_{lij} \\ \label{tauoneevolutioneq} & & {} + \frac{1}{6} \, (\nab{p} X_q)
g^{pa} g^{qb} \ph_{abk} g^{kl} \ph_{lij}
\end{eqnarray}
and
\begin{eqnarray} \nonumber 
\ddt (\tau_2)_{ij} & = & \pi_7 ( \{h, \tau_2\}_{ij} ) + \pi_{14}
( {}- [h, \tau_3]_{ij} + \{h, \tau_1\}_{ij} + \frac{1}{2} \, \{h, \tau_2\}_{ij} )
\\ \nonumber & & {}+ \pi_7 ( - \frac{1}{3} \, [X, \tau_2]_{ij} ) + \pi_{14}
( \{X, \tau_3\}_{ij} + \frac{1}{3} \, [X, \tau_1]_{ij} )
\\ \nonumber & & {} - (\nab{k} h_{li}) g^{ka} g^{lb} \ph_{abj} +
(\nab{k} h_{lj}) g^{ka} g^{lb} \ph_{abi} \\ \nonumber & & {}- \frac{1}{3} \, g^{pq}
(\nab{p} h_{qk}) g^{kl} \ph_{lij} + \frac{1}{3} \, (\nab{k} \tr_g(h))g^{kl}
\ph_{lij} \\ \label{tautwoevolutioneq} & & {} - \nab{i} X_j + \nab{j} X_i
+ \frac{1}{3} \, (\nab{p} X_q) g^{pa} g^{qb} \ph_{abk} g^{kl} \ph_{lij}
\end{eqnarray}
\end{prop}
\begin{proof}
This is a tedious but straightforward computation. By
Proposition~\ref{torsionformsprop},
\begin{eqnarray*}
\ddt (\tau_1)_{ij} & = & \frac{1}{3} \left( \ddt C_{ij} \right) - \frac{1}{6} 
\left( \ddt (P(C))_{ij} \right) \\ \ddt (\tau_2)_{ij} & = & -\frac{4}{3} \left(
\ddt C_{ij} \right) - \frac{1}{3} \left( \ddt (P(C))_{ij} \right)
\end{eqnarray*}
where $C_{ij} = \frac{1}{2} (T_{ij} - T_{ji}) = (\tau_1)_{ij} - \frac{1}{2} \,
(\tau_2)_{ij}$. Now we use Theorem~\ref{fulltorsionevolutionthm}, together with
Lemma~\ref{skewtorsionpreliminarylemma} and repeated applications of
Proposition~\ref{wtwprop}. We also need Lemmas~\ref{g2identities1lemma}
and~\ref{g2identities2lemma} to put the terms involving $\nab{}h$ and $\nab{}X$
in the form given in~\eqref{tauoneevolutioneq} and~\eqref{tautwoevolutioneq}.
We omit the details.
\end{proof}

\smallskip

\begin{rmk} \label{torsionevolutionfinalrmk}
The evolution equations given in Propositions~\ref{symmetrictorsionevolutionprop} and~\ref{skewtorsionevolutionprop} are clearly quite complicated.
Ideally one could find an $h_{ij}$ and $X^k$ depending on the four
torsion tensors (and their covariant derivatives) for which these evolution
equations become simpler. The sequel~\cite{KY} to this paper will
discuss several specific flows of \Gs s.
\end{rmk}

\section{Bianchi-type identities in \G-geometry} \label{bianchisec}

In this section, we apply the evolution equations from
Section~\ref{generalflowssec} to derive Bianchi-type identities for
manifolds with \Gs. As a consequence, we obtain explicit formulas for
the Ricci tensor and part of the Riemann curvature tensor in terms of
the full torsion tensor. This leads to new simple proofs of some known
results in \G-geometry.

\subsection{Identities via diffeomorphism invariance} \label{diffeoinvsec}

Let $\alpha$ be a smooth tensor on a manifold $M$, defined 
entirely in terms of some other tensor $\beta$. We write $\alpha =
\alpha[\beta]$. Suppose $f_t$ is a 
one-parameter family of diffeomorphisms of $M$ such that $f_0$ is the
identity, and $\ddt f_t = Y$ for some smooth vector field $Y$ on $M$.
Then, by diffeomorphism invariance, we have
\begin{equation*}
f_t^*( \alpha[\beta] ) = \alpha[ f_t^*(\beta) ]
\end{equation*}
where $f_t^*$ denotes the pullback by $f_t$. Differentiating with respect
to $t$ and setting $t=0$ gives
\begin{equation} \label{diffeoinveq}
\mathcal{L}_Y \alpha = (D\alpha) [\mathcal{L}_Y \beta]
\end{equation}
where $D\alpha$ is the linearization of the tensor $\alpha$ as a function
of $\beta$. Explicitly, if $\beta(t)$ satisfies $\ddt \beta(t) = \beta_1$,
then $(D\alpha)(\beta_1) = \ddt (\alpha[\beta(t)])$.
Since~\eqref{diffeoinveq} holds for any vector field $Y$, this yields
identities involving $\alpha$. This idea was exploited by Kazdan
in~\cite{Kaz} to show that the first and second Bianchi identities of
Riemannian geometry were consequences of the diffeomorphism invariance
of the Riemann curvature tensor $R_{ijkl}$ as a function of the metric
$g_{ij}$. A discussion can be found in~\cite{CK}.

We will apply this idea to the setting of \G-geometry, to derive
`Bianchi-type' identities. First, we note that from~\eqref{liederivativeeq},
we have
\begin{equation} \label{liepheq}
(\mathcal{L}_Y \ph)_{ijk} = (\nab{Y}\ph)_{ijk} + (\nab{i} Y^l)\ph_{ljk}
+ (\nab{j} Y^l)\ph_{ilk} + (\nab{k} Y^l)\ph_{ijl}
\end{equation}
This should be compared to the well-known analogous expression for the 
Riemannian metric:
\begin{equation} \label{liemetriceq}
(\mathcal{L}_Y g)_{ij} = (\nab{i} Y^l)g_{lj} + (\nab{j} Y^l)g_{il} =
\nab{i} Y_j + \nab{j} Y_i
\end{equation}
which also follows from~\eqref{liederivativeeq} since $\nab{} g = 0$.

\smallskip

\begin{prop} \label{firstbianchiprop}
The diffeomorphism invariance of the metric $g$ as a function
of the $3$-form $\ph$ is equivalent to the vanishing of the
$\wtho \oplus \wtht$ component of $\nab{Y}\ph$ for
any vector field $Y$. This is the fact which was proved earlier
in Lemma~\ref{torsionsymmetrieslemma}. See also the discussion following
Remark~\ref{torsionsymmetriesrmk1}.
\end{prop}
\begin{proof}
Equation~\eqref{diffeoinveq} for $g[\ph]$ says
\begin{equation} \label{firstbianchitempeq}
(\mathcal{L}_Y g)_{ij} = \nab{i} Y_j + \nab{j} Y_i = (Dg)[\mathcal{L}_Y \ph]
\end{equation}
By equation~\eqref{metricevolutioneq}, the right hand side above equals
$2 \, h_{ij}$ where the symmetric $2$-tensor $h_{ij}$ and the vector field
$X^k$ are defined by the decomposition of the $3$-form $\mathcal{L}_Y \ph$ as
\begin{equation} \label{liepheq2}
(\mathcal{L}_Y \ph)_{ijk} = \frac{1}{2} h_{im} g^{ml} \ph_{ljk} + X^m
\ps_{mijk}
\end{equation}
Consider the last three terms of~\eqref{liepheq}. They can be written as
\begin{equation} \label{liephtempeq}
A_{im} g^{ml}\ph_{ljk} + A_{jm} g^{ml}\ph_{ilk} + A_{km} g^{ml}\ph_{ijl}
\end{equation}
where $A_{ij} = \nab{i} Y_j$. From the discussion following
Remark~\ref{wtwvecsrmk}, the $\wtho \oplus \wtht$ component
of this expression is $\frac{1}{2} S_{im} g^{ml} \ph_{ljk}$ where
$S_{ij} = \frac{1}{2} (A_{ij} + A_{ji}) = \frac{1}{2} (\nab{i} Y_j + 
\nab{j} Y_i)$. Since $2 \, S_{ij}$ is already equal to the left
hand side of~\eqref{firstbianchitempeq}, we see that
equation~\eqref{firstbianchitempeq} holds for all $Y$ if and
only if the $\wtho \oplus \wtht$ component of $\nab{Y}\ph$ vanishes for 
all $Y$.
\end{proof}

We pause here to note that this result enables us to explicitly describe
the decomposition of $\mathcal{L}_Y \ph$ into a symmetric tensor $h_{ij}$
and a vector field $X^k$ as given in~\eqref{liepheq2}.
We already remarked in the proof of Proposition~\ref{firstbianchiprop} that 
\begin{equation} \label{liepheq3}
h_{ij} = \frac{1}{2} \, (\nab{i} Y_j + \nab{j} Y_i)
\end{equation}
Also, from~\eqref{skewthreeformeq} and~\eqref{liephtempeq}, we see that
the $\wths$ component of the last three terms of~\eqref{liepheq} is
\begin{eqnarray*}
Z^n \ps_{nijk} \qquad \text{ where} \qquad Z^n & = & -\frac{1}{2} \, \left( \frac{1}{2}
(\nab{a} Y_b - \nab{b} Y_a ) \right) g^{ai} g^{bj} \ph_{ijk} g^{kn} \\ & = &
-\frac{1}{2} \, (\nab{a} Y_b) g^{ai} g^{bj} \ph_{ijk} g^{kn}
\end{eqnarray*}
using the skew-symmetry of $\ph_{ijk}$. Now combining this with
$\nab{l} \ph_{ijk} = T_{lm} g^{mn} \ps_{nijk}$, we see that $X^k$ is
\begin{equation} \label{liepheq4}
X^k = Y^l T_{lm} g^{mk} -\frac{1}{2} \, (\nab{a} Y_b) g^{ai} g^{bj} \ph_{ijm} g^{mk}
\end{equation}
in equation~\eqref{liepheq2}.

Using~\eqref{nablapsieq2} and Theorem~\ref{psievolutionthm}, one can consider
the analogous calculation of diffeomorphism invariance for the $4$-form $\ps$ as a
function of $\ph$. It is straightforward to check that in this case no new information is
obtained. Therefore we turn our attention to the full torsion tensor $T_{lm}$.

\smallskip

\begin{thm} \label{secondbianchithm}
The diffeomorphism invariance of the full torsion tensor $T$ as a function
of the $3$-form $\ph$ is equivalent to the following identity:
\begin{equation} \label{secondbianchieq}
\boxed{\, \, \nab{i} T_{jl} - \nab{j} T_{il} = T_{ia} T_{jb} g^{am} g^{bn} \ph_{mnl}
+ \frac{1}{2} R_{ijab} g^{am} g^{bn} \ph_{mnl} \, \, }
\end{equation}
\end{thm}
\begin{proof}
Equation~\eqref{diffeoinveq} for $T[\ph]$ says
\begin{equation} \label{secondbianchitempeq}
(\mathcal{L}_Y T)_{ij} = Y^l (\nab{l} T_{ij}) + (\nab{i} Y^l) T_{lj} + (\nab{j} Y^l) T_{il}
= (DT)[\mathcal{L}_Y \ph]
\end{equation}
By Theorem~\ref{fulltorsionevolutionthm}, the right hand side above equals
\begin{equation} \label{secondbianchitempeq2}
T_{il} g^{lm} h_{mj} + T_{il} g^{lm} X_{mj} + (\nab{k} h_{li}) g^{ka} g^{lb} \ph_{abj}
+ \nab{i} X_j
\end{equation}
where the symmetric $2$-tensor $h_{ij}$ and the vector field
$X^k$ are given by~\eqref{liepheq3} and~\eqref{liepheq4}, respectively.
Remark~\ref{wtwvecsrmk} shows that
\begin{eqnarray} \nonumber
X_{mj} & = & X^k \ph_{kmj} = Y^l T_{ln} g^{nk} \ph_{kmj} 
-\frac{1}{2} \, (\nab{a} Y_b) g^{ap} g^{bq} \ph_{pqn} g^{nk} \ph_{kmj} \\ \nonumber
& = &  Y^l T_{ln} g^{nk} \ph_{kmj} 
-\frac{1}{2} \, (\nab{a} Y_b) g^{ap} g^{bq} ( g_{pm} g_{qj} - g_{pj} g_{qm} - \ps_{pqmj})
\\ \label{secondbianchitempeq3} & = &  Y^l T_{ln} g^{nk} \ph_{kmj} 
-\frac{1}{2} \, \nab{m} Y_j + \frac{1}{2} \, \nab{j} Y_m
+ \frac{1}{2} \, (\nab{a} Y_b) g^{ap} g^{bq} \ps_{pqmj}
\end{eqnarray}
Substitute~\eqref{liepheq3}, ~\eqref{liepheq4}, and~\eqref{secondbianchitempeq3}
into~\eqref{secondbianchitempeq2} to obtain
\begin{eqnarray*}
& & T_{il} g^{lm} \left( \frac{1}{2} \, (\nab{m} Y_j + \nab{j} Y_m) \right) \\ & & {} +
T_{il} g^{lm} \left(  Y^b T_{bn} g^{nk} \ph_{kmj} 
-\frac{1}{2} \, \nab{m} Y_j + \frac{1}{2} \, \nab{j} Y_m
+ \frac{1}{2} \, (\nab{a} Y_b) g^{ap} g^{bq} \ps_{pqmj} \right) \\ & & 
{}+ \nab{k} \left( \frac{1}{2} \, (\nab{l} Y_i + \nab{i} Y_l) \right) g^{ka} g^{lb}
\ph_{abj} + \nab{i} \left( Y^b T_{bj} -\frac{1}{2}(\nab{a} Y_b) g^{ap} b^{bq} \ph_{pqj}
\right)
\end{eqnarray*}
Now we expand the above expression using $\nab{l} \ph_{abc} = 
T_{lm} g^{mn} \ps_{nabc}$ and the product rule, and collect terms. After some cancellation, we are left with
\begin{eqnarray*}
& & (\nab{j} Y^l) T_{il} + (\nab{i} Y^b) T_{bj} + Y^b (\nab{i} T_{bj}) 
+ Y^b T_{bn} T_{il} g^{nk} g^{lm} \ph_{kmj} \\ & & {}+ \frac{1}{2} (\nab{k} \nab{l} Y_i +
\nab{k} \nab{i} Y_l) g^{ka} g^{lb} \ph_{abj} - \frac{1}{2} (\nab{i} \nab{a} Y_b)
g^{ap} g^{bq} \ph_{pqj}
\end{eqnarray*}
for the right hand side of~\eqref{secondbianchitempeq}. Substituting this expression
into~\eqref{secondbianchitempeq}, cancelling terms, rearranging, and relabeling indices
gives
\begin{eqnarray} \label{secondbianchitempeq4}
Y^l (\nab{l} T_{ij}) - Y^l (\nab{i} T_{lj}) & = & Y^l T_{la} T_{ib} g^{am} g^{bn} \ph_{mnj}
\\ & & \nonumber {}+ \frac{1}{2} (\nab{k} \nab{l} Y_i + \nab{k} \nab{i} Y_l -
\nab{i} \nab{k} Y_l) g^{ka} g^{lb} \ph_{abj}
\end{eqnarray}
Consider the last term in~\eqref{secondbianchitempeq4}. Using the Ricci
identity~\eqref{ricciidentityeq} and the skew-symmetry of $\ph$, we have
\begin{eqnarray*}
\frac{1}{2} (\nab{k} \nab{l} Y_i) g^{ka} g^{lb} \ph_{abj} & = & \frac{1}{4} (\nab{k}
\nab{l} Y_i- \nab{l} \nab{k} Y_i) g^{ka} g^{lb} \ph_{abj} \\ & = & -\frac{1}{4} R_{klim}
Y^m g^{ka} g^{lb} \ph_{abj}
\end{eqnarray*}
and also
\begin{eqnarray*}
\frac{1}{2} (\nab{k} \nab{i} Y_l - \nab{i} \nab{k} Y_l) g^{ka} g^{lb} \ph_{abj}
& = & -\frac{1}{2} R_{kilm} Y^m g^{ka} g^{lb} \ph_{abj} \\ & = & -\frac{1}{4} (R_{kilm}
- R_{likm}) Y^m g^{ka} g^{lb} \ph_{abj}
\end{eqnarray*}
Putting these two expressions together, the last term of~\eqref{secondbianchitempeq4}
becomes
\begin{eqnarray*}
-\frac{1}{4} (R_{klim} + R_{kilm} - R_{likm}) Y^m g^{ka} g^{lb} \ph_{abj} & = & 
-\frac{1}{4} (-R_{mikl} + R_{mlik} + R_{mkli}) Y^m g^{ka} g^{lb} \ph_{abj} \\ 
& = & -\frac{1}{4} (-R_{mikl} -R_{mikl}) Y^m g^{ka} g^{lb} \ph_{abj} \\ & = &
 \frac{1}{2} R_{mikl} Y^m g^{ka} g^{lb} \ph_{abj}
\end{eqnarray*}
where we have used the symmetries of $R_{ijkl}$ in the first line, and the Riemannian
first Bianchi identity~\eqref{riemannbianchieq} in the second line.
Hence~\eqref{secondbianchitempeq4} becomes
\begin{equation*}
Y^l (\nab{l} T_{ij}) - Y^l (\nab{i} T_{lj}) = Y^l T_{la} T_{ib} g^{am} g^{bn} \ph_{mnj}
+ \frac{1}{2} Y^m R_{mikl} g^{ka} g^{lb} \ph_{abj}
\end{equation*}
Since this must hold for any $Y$, after relabelling indices this is
exactly~\eqref{secondbianchieq}.
\end{proof}

\smallskip

\begin{rmk} \label{secondbianchirmk}
The identity~\eqref{secondbianchieq} can also be established directly by
using~\eqref{fulltorsioneq}, ~\eqref{nablapsieq2}, Lemma~\ref{g2identities2lemma}, and
the Ricci identities. However the proof above shows
that~\eqref{secondbianchieq} is equivalent to the diffeomorphism invariance of the full
torsion tensor $T_{lm}$, and as such it deserves to be called a Bianchi-type identity
for \G-geometry.
\end{rmk}

\smallskip

\begin{cor} \label{contractedsecondbianchicor}
The following identity holds
\begin{equation} \label{contractedsecondbianchieq}
\boxed{\, \, \frac{7}{4} \nab{i} \tau_0 - g^{jl} \nab{j} T_{il} =
-6 \, T_{ia} (\tau_1)^a - \frac{1}{2} R_{abji} g^{am} g^{bn} g^{jl} \ph_{mnl} \, \, }
\end{equation}
We call it the contracted Bianchi-type identity for \G-geometry.
\end{cor}
\begin{proof}
This follows easily by contracting~\eqref{secondbianchieq} on $j$ and $l$, and using
the fact that $T_{lm} = \frac{1}{4} \tau_0 \, g_{lm} - (\tau_3)_{lm} + (\tau_1)_{lm} -
\frac{1}{2} (\tau_2)_{lm}$ and equation~\eqref{wtwvecseq} for $\tau_1$. 
\end{proof}

\smallskip

\begin{rmk} \label{contractedsecondbianchirmk}
The usefulness of~\eqref{contractedsecondbianchieq} is that it shows that the 
`divergence-like' expression $g^{jl} \nab{j} T_{il}$ is equal to a gradient plus
zero order terms in the torsion and curvature. This can potentially be used to simplify
expressions, in the same way that the contracted second Bianchi identity of
Riemannian geometry is used in Ricci flow.
\end{rmk}

\subsection{Curvature formulas in terms of the torsion tensor}
\label{curvaturesec}

We now examine some consequences of Theorem~\ref{secondbianchithm}. For $i$ and $j$
fixed, the Riemann curvature tensor $R_{ijkl}$ is skew-symmetric in $k$ and $l$. 
Hence we can use~\eqref{wtwsprojeq} and~\eqref{wtwfprojeq} to decompose it as
\begin{equation*}
R_{ijkl} = (\pi_7 (\mathrm{Riem}))_{ijkl} + (\pi_{14} (\mathrm{Riem}))_{ijkl}
\end{equation*}
where
\begin{eqnarray} \label{riem7eq}
(\pi_7 (\mathrm{Riem}))_{ijkl} & = & \frac{1}{3} R_{ijkl} - \frac{1}{6} R_{ijab} g^{ap} 
g^{bq} \ps_{pqkl} \\ \label{riem14eq} (\pi_{14} (\mathrm{Riem}))_{ijkl} & = & \frac{2}{3}
R_{ijkl} + \frac{1}{6} R_{ijab} g^{ap} g^{bq} \ps_{pqkl}
\end{eqnarray}
Furthermore, by Remark~\ref{wtwvecsrmk}, we also know that
\begin{equation} \label{pi7riemanntempeq}
(\pi_7 (\mathrm{Riem}))_{ijkl} = (\pi_7 (\mathrm{Riem}))^m_{ij} \, \ph_{mkl}
\end{equation}
where
\begin{equation*}
(\pi_7 (\mathrm{Riem}))^m_{ij} = \frac{1}{6} R_{ijkl} g^{ka} g^{lb} \ph_{abc} g^{cm}
\end{equation*}
Thus the identity~\eqref{secondbianchieq} can be expressed as
\begin{equation} \label{pi7riemmaneq}
3 \, (\pi_7 (\mathrm{Riem}))^m_{ij} g_{ml} = \nab{i} T_{jl} - \nab{j} T_{il} -
T_{ia} T_{jb} g^{am} g^{bn} \ph_{mnl}
\end{equation}

\smallskip

\begin{rmk} \label{bianchibryantrmk}
The fact that the $\wtws$ part of the Riemann curvature tensor can be expressed
entirely in terms of the full torsion tensor $T_{lm}$ was first demonstrated
using frame bundle calculations in section [4.5] of~\cite{Br3}.
\end{rmk}

\smallskip

\begin{cor} \label{ambrosesingercor}
If the \Gs\ $\ph$ is torsion-free, then the Riemann curvature tensor $R_{ijkl} \in
\mathrm{Sym}^2 (\Omega^2)$ actually takes values $\mathrm{Sym}^2(\wtwf)$, where
$\wtwf \cong \lieg$, the Lie algebra of $\G$.
\end{cor}
\begin{proof}
Setting $T = 0$ in~\eqref{pi7riemmaneq} shows the for fixed $i$, $j$, we have
$R_{ijkl} \in \wtwf$ as a skew-symmetric tensor in $k$, $l$. The result now follows
from the symmetry $R_{ijkl} = R_{klij}$.
\end{proof}

\smallskip

\begin{rmk} \label{ambrosesingerrmk}
This result is well-known. When $T=0$, the three form $\ph$ is parallel and
the Riemannian holonomy of the associated metric is contained in the group $\G$. It is a general fact that the Riemann curvature tensor of a metric with holonomy
contained in a group $\mathrm{G}$ is an element of $\mathrm{Sym}^2(\mathfrak{g})$, where
$\mathfrak{g}$ is the Lie algebra of $\mathrm{G}$. Here we have a direct proof of this
fact in the \G\ case.
\end{rmk}

\smallskip

\begin{lemma} \label{riccilemma}
Let $Q_{ijkl} = R_{ijab} g^{ap} g^{bq} \ps_{pqkl}$. Then we have $Q_{ijkl} g^{il} = 0$.
\end{lemma}
\begin{proof}
We use the Riemannian Bianchi identity~\eqref{riemannbianchieq} and compute:
\begin{eqnarray*}
Q_{ijkl} g^{il} & = & R_{ijab} g^{ap} g^{bq} \ps_{pqkl} g^{il} \\ & = & -(R_{jaib} + 
R_{aijb}) g^{ap} g^{bq} \ps_{pqkl} g^{il} \\ & = & R_{ajib} g^{ap} g^{bq} \ps_{pqkl} g^{il}
+ R_{bjai} g^{ap} g^{bq} \ps_{pqkl} g^{il} \\ & = & {}-R_{ajib} g^{il} g^{bq} \ps_{lqkp}
g^{ap} - R_{bjai} g^{ap} g^{il} \ps_{plkq} g^{bq} \\ & = & {}-Q_{ajkp} g^{ap} - Q_{bjkq}
g^{bq}
\end{eqnarray*}
and thus $3 \, Q_{ijkl} g^{il} = 0$.
\end{proof}

\smallskip

\begin{rmk} \label{riccilemmaCIrmk}
This is essentially Proposition [3.2] in Cleyton-Ivanov~\cite{CI}.
\end{rmk}

\smallskip

\begin{cor} \label{riccicor}
The Ricci tensor $R_{jk}$ can be expressed as
\begin{equation*}
R_{jk} = R_{ijkl} g^{il} = 3 \, (\pi_7 (\mathrm{Riem}))_{ijkl} g^{il} = \frac{3}{2} \,
(\pi_{14} (\mathrm{Riem}))_{ijkl} g^{il}
\end{equation*}
\end{cor}
\begin{proof}
This follows immediately from Lemma~\ref{riccilemma} and equations~\eqref{riem7eq}
and~\eqref{riem14eq}.
\end{proof}

\smallskip

\begin{cor} \label{ricciflatcor}
The metric of a torsion-free \Gs\ is necessarily Ricci-flat.
\end{cor}
\begin{proof}
From Corollary~\ref{ambrosesingercor}, we have $\pi_7 (\mathrm{Riem}) = 0$.
The result then follows from Corollary~\ref{riccicor}.
\end{proof}

\smallskip

\begin{rmk} \label{ricciflatrmk} When the \Gs\ is torsion-free, the
  associated metric $g$ has holonomy contained in \G, and it is a
  classical fact due to Bonan~\cite{Bo} that such a metric is Ricci-flat. The
  original proof depends on the fact that the conditions on the
  Riemann curvature tensor implied by the Lie algebra structure of
  $\lieg$ are too strong to have non-zero Ricci tensor. Here we have a
  direct proof. We were also informed of a simple
  representation-theoretic proof of this fact by Bryant~\cite{BrP}.
\end{rmk}

\smallskip

\begin{rmk} \label{ricciflatrmk2}
Corollary~\ref{riccicor} also implies that a \Gs\ whose Riemann curvature tensor is in $\mathrm{Sym}^2(\wtws)$
is also Ricci-flat. However, this can only occur if $R_{ijkl}$ vanishes
identically (that is, the manifold is {\em flat}). To see this, suppose 
$\pi_{14}(\mathrm{Riem}) = 0$. Then for fixed $i,j$, we have that $R_{ijkl} = A^m_{ij}
\ph_{mkl}$ because it is in $\wtws$ as a skew-symmetric matrix in $k,l$. But for
fixed $k,l$, it is also in $\wtws$ as a skew-symmetric matrix in $i,j$, so
$A^m_{ij} = A^{nm} \ph_{nij}$, and hence $R_{ijkl} = A^{nm} \ph_{nij} \ph_{mkl}$.
Since $R_{ijkl} = R_{klij}$, we have $A^{nm} = A^{mn}$ is symmetric. Now by
Corollary~\ref{riccicor}, we have
\begin{eqnarray*}
0 = R_{jk} & = & R_{ijkl} g^{il} = A^{nm} \ph_{nij} \ph_{mkl} g^{il} \\
& = & -A^{nm} ( g_{nm} g_{jk} - g_{nk} g_{jm} - \ps_{njmk} ) \\
& = & - \tr_g(A) g_{jk} + A_{jk} + 0
\end{eqnarray*}
Taking the trace gives $6 \tr_g(A) = 0$, hence $A_{jk} = 0$, and therefore
$R_{ijkl} = 0$.
\end{rmk}

\smallskip

Corollary~\ref{riccicor} allows us to derive an explicit formula for the
Ricci curvature of the metric of any \Gs.

\smallskip

\begin{prop} \label{ricciprop}
Given a \Gs\ $\ph$ with full torsion tensor $T_{lm}$, its associated metric
$g$ has Ricci curvature $R_{jk}$ given by
\begin{equation} \label{riccieq}
\fbox{
  \begin{minipage}{0.7\textwidth}
    \leftline{$\, \, R_{jk} = (\nab{i} T_{jm} - \nab{j} T_{im}) \ph_{nkl} g^{mn} g^{il}$}
\vskip 0.05in
    \leftline{$\qquad \qquad \, \, {}-  T_{jl} g^{li} T_{ik} + \tr_g(T) T_{jk}
- T_{jb} T_{ia} g^{il} g^{ap} \ps_{lpqk} g^{bq}$}
  \end{minipage}
}
\end{equation}
This can also be written in the form
\begin{eqnarray} \label{riccieq2}
R_{jk} & = & \nab{i}( T_{jm} \ph_{nkl} g^{mn} g^{il}) - \nab{j} (T_{im} \ph_{nkl}
g^{mn} g^{il}) \\ \nonumber & & \qquad {} -  T_{jl} g^{li} T_{ik} + \tr_g(T) T_{jk}
+ T_{jb} T_{ia} g^{il} g^{ap} \ps_{lpqk} g^{bq}
\end{eqnarray}
(note the change of sign on the last term.)
\end{prop}
\begin{proof}
We use Corollary~\ref{riccicor} and equations~\eqref{pi7riemanntempeq}
and~\eqref{pi7riemmaneq} to compute:
\begin{eqnarray*}
R_{jk} & = & (\nab{i} T_{jm} - \nab{j} T_{im} - T_{ia} T_{jb} g^{ap} g^{bq} \ph_{pqm}) \, g^{mn}
\ph_{nkl} g^{il} \\ & = & (\nab{i} T_{jm} - \nab{j} T_{im}) \ph_{nkl} g^{mn} g^{il} - 
T_{ia} T_{jb} g^{ap} g^{bq} g^{il} (g_{pk} g_{ql} - g_{pl} g_{qk} - \ps_{pqkl})
\\ & = & (\nab{i} T_{jm} - \nab{j} T_{im}) \ph_{nkl} g^{mn} g^{il} - T_{ik} T_{jl} g^{il}
+ T_{il} T_{jk} g^{il} - T_{jb} T_{ia} g^{il} g^{ap} \ps_{lpqk} g^{bq}
\end{eqnarray*}
using Lemma~\ref{g2identities1lemma}. To obtain~\eqref{riccieq2}, we write
\begin{eqnarray*}
(\nab{i} T_{jm}) \ph_{nkl} g^{mn} g^{il} & = & \nab{i} (T_{jm} \ph_{nkl} g^{mn}
g^{il}) - T_{jm} (\nab{i} \ph_{nkl}) g^{mn} g^{il} \\ & = & 
\nab{i} (T_{jm} \ph_{nkl} g^{mn} g^{il}) - T_{jm} T_{ip} g^{pq} \ps_{qnkl} g^{mn}
g^{il}
\end{eqnarray*}
and similarly for the other term involving $\nab{} T$. The two extra terms add up
and are negatives of the last term in~\eqref{riccieq}. The details are left
to the reader.
\end{proof}

\smallskip

\begin{rmk} \label{riccibryantrmk}
  Robert Bryant derived a formula, equation [4.27] in~\cite{Br3}, for
  the Ricci tensor in terms of the intrinsic torsion in the language
  of frame bundles and canonical \G-connections. Our
  formula~\eqref{riccieq} is an equivalent expression in the language
  of local coordinates. Additionally, equation [4.30] in~\cite{Br3}
  expresseses the Ricci tensor in terms of the four independent
  torsion forms and the representation-theoretic operators `$Q$' and
  `$j$.' Our expression for $R_{jk}$ in terms of $T_{lm}$ can also be
  combined with Theorem~\ref{torsionthm} to obtain an expression for
  $R_{jk}$ in terms of $\tau_0$, $\tau_1$, $\tau_2$, and $\tau_3$.
\end{rmk}

\smallskip

\begin{rmk} \label{ricciCIrmk}
Cleyton and Ivanov~\cite{CI} also derived a formula for the Ricci tensor
for {\em closed} \Gs s in terms of $\dl \ph$.
\end{rmk}

\smallskip

\begin{cor} \label{scalarcor}
Given a \Gs\ $\ph$ with full torsion tensor $T_{lm}$, its associated metric
$g$ has scalar curvature $R$ given by
\begin{equation} \label{scalareq}
\boxed{\, \, R = -12 \, g^{il} \nab{i} (\tau_1)_j + \frac{21}{8} \, \tau_0^2
- {|\tau_3|}^2 + 5\, {|\tau_1|}^2 - \frac{1}{4} \, {|\tau_2|}^2 \, \, }
\end{equation}
where ${|A|}^2 = A_{ij} A_{kl} g^{ik} g^{jl}$ is the matrix norm.
\end{cor}
\begin{proof}
We use~\eqref{riccieq2} and trace $R = g^{jk} R_{jk}$ to obtain $R = $
\begin{equation} \label{scalartempeq}
 -2 \, g^{il} \nab{i} (T_{jm} \ph_{knl} g^{jk} g^{mn}) - T_{jl} T_{ik} g^{li}
g^{jk} + {(\tr_g(T))}^2 + T_{jb} (T_{ia} g^{il} g^{ap} \ps_{lpqk}) g^{bq} g^{jk}
\end{equation}
Now recall that $T_{lm} = \frac{1}{4} \tau_0 \, g_{lm} - (\tau_3)_{lm} + (\tau_1)_{lm}
- \frac{1}{2} (\tau_2)_{lm}$ from Theorem~\ref{torsionthm}. Using the fact that the
four components of $T_{lm}$ are mutually orthogonal, and that $\tau_3$ is symmetric
while $\tau_1$ and $\tau_2$ are skew-symmetric, the second term
in~\eqref{scalartempeq} becomes
\begin{eqnarray*}
& & {}- \left( \frac{1}{4} \tau_0 \, g_{jl} - (\tau_3)_{jl} + (\tau_1)_{jl}
- \frac{1}{2} (\tau_2)_{jl} \right) \left( \frac{1}{4} \tau_0 \, g_{ik} -
(\tau_3)_{ik} + (\tau_1)_{ik} - \frac{1}{2} (\tau_2)_{ik} \right) g^{li} g^{jk} \\
& = & {}- \frac{7}{16} \, \tau_0^2 - {|\tau_3|}^2 + {|\tau_1|}^2 + \frac{1}{4} \,
{|\tau_2|}^2
\end{eqnarray*}
The third term is $\frac{49}{16} \, \tau_0^2$. We use Proposition~\ref{wtwprop}
on the third term to write it as
\begin{eqnarray*}
& & T_{jb} \left( \left( \frac{1}{4} \tau_0 \, g_{ia} - (\tau_3)_{ia} +
(\tau_1)_{ia} - \frac{1}{2} (\tau_2)_{ia} \right) g^{il} g^{ap} \ps_{lpqk} \right)
g^{bq} g^{jk} \\ & = & T_{jb} \left( 0 + 0 - 4 \, (\tau_1)_{qk} - (\tau_2)_{qk}
\right) g^{bq} g^{jk} = 4\, {|\tau_1|}^2 - \frac{1}{2} \, {|\tau_2|}^2
\end{eqnarray*}
Combining these three expressions gives the last four terms of~\eqref{scalareq}.
Finally, for the first term of~\eqref{scalartempeq} we note that by
Proposition~\ref{wtwprop} and~\eqref{wtwvecseq}, we have
\begin{equation*}
T_{jm} \ph_{knl} g^{jk} g^{mn} = (\tau_1)_{jm} \ph_{knl} g^{jk} g^{mn} =
6 \, (\tau_1)_l  
\end{equation*}
and the proof is complete.
\end{proof}

\smallskip

\begin{rmk} \label{scalarbryantrmk}
The equation~\eqref{scalareq} agrees exactly with equation [4.28] of~\cite{Br3}.
It looks different because Bryant is using the differential form norms, while we
use the matrix norms in our equation~\eqref{scalareq}.
\end{rmk}

\smallskip

We close this section by considering a specific subclass of \Gs s, namely those for which
the only nonvanishing component of the torsion tensor is $\tau_0$, the scalar torsion.
These are sometimes called {\em nearly} \G\ manifolds. In this case we have $d \ps = 0$ and 
$d\ph = \tau_0 \ps$, so differentiating the second equation gives $d \tau_0 \wedge \ps = 0$, and
hence $d \tau_0 = 0$, since wedge product with $\ps$ is an isomorphism from $\wons$ to $\wfis$.
Thus $\tau_0$ is constant and hence $\nab{i} T_{jl} = 0$ in this case. Substituting
into~\eqref{riccieq2} gives
\begin{equation*}
R_{jk} = \frac{3}{8} \tau_0^2 \, g_{jk}
\end{equation*}
and we recover the known fact that nearly \G\ manifolds are always Einstein with positive
Einstein constant. The squashed $7$-sphere $S^7$ is an example of such a manifold.

\appendix

\section{Identites involving $\ph$ and $\ps$}
\label{g2identitiessec}

In this section we derive several useful identities involving the
3-form $\ph$ and the 4-form $\ps$ for a manifold with a \Gs. Some of
these have appeared previously in~\cite{CI} and in~\cite{Br3}. Here
we give proofs of them as well as some new identities. Some identities are
described invariantly and others in terms of arbitrary local
coordinates.

\subsection{Basic relations of $\G$-geometry}
\label{basicrelationssec}

The main ingredients for deriving the identities of \G-geometry are
the following:

\smallskip

\begin{lemma} \label{g2knownlemma}
The metric $g$, cross product $\times$, and 3-form $\ph$ satisfy
the following relations:
\begin{eqnarray}
\label{g2maineq} g ( u \times v , w) & = & \ph (u,v,w) \\
\label{crosseq} {(u \times v)}^{\flat} & = & v \hk u \hk \ph \; =
\; \st ( \us \wedge \vs \wedge \ps ) \\ \label{iteratedeq} u \times
\left( v \times w \right) & = & - g(u, v) w + g(u, w) v - {\left( u
\hk v \hk w \hk \ps \right) }^{\sharp}
\end{eqnarray}
where $u, v, w$ are vector fields and $\vs$ denotes the 1-form
which is metric dual to $v$.
\end{lemma}
\begin{proof}
See, for example, \cite{K1} for a proof. Note that in~\cite{K1}
there is a sign error in~\eqref{iteratedeq}. With the sign
convention used in that paper, the last term should have a plus
sign instead of a minus sign. In this paper we use the opposite
orientation convention, and so~\eqref{iteratedeq} has a minus sign
in front of the last term. See also~\cite{K2} for more about sign
conventions and orientations. 
\end{proof}
From Lemma~\ref{g2knownlemma}, we get:

\smallskip

\begin{cor} \label{g2identitycor}
Let $a, b, c, d$ be vector fields. Then
\begin{equation*}
g( a \times b, c \times d) = g( a \wedge b, c \wedge d) - \ps( a,
b, c, d)
\end{equation*}
\end{cor}
\begin{proof}
We compute
\begin{eqnarray*}
g( a \times b, c \times d) & = & \ph( a, b, c \times d) \\ & = & -
\ph( a , c \times d, b) \\ & = & - g( a \times (c \times d) , b)
\\ & = & - g( -g(a,c) d + g(a,d) c - {(a \hk c \hk d \hk
\ps)}^{\sharp}, b) \\ & = & g(a,c) g(b,d) - g(a,d) g(b,c) +
\ps(d,c,a,b)
\end{eqnarray*}
which is what we wanted to prove.
\end{proof}

We will also have occasion to use the following relations, whose
proofs can be found in~\cite{K1}.

\smallskip

\begin{prop} \label{g2relationsprop}
Let $\alpha$ be a $1$-form on $M$, let $w$ be a vector field on $M$,
and $\ws$ be the $1$-form dual to $w$. Then the following relations
hold:
\begin{eqnarray*}
|\ph|^2 = 7 \qquad \qquad \quad & \qquad \qquad & \qquad \qquad \quad 
|\ps|^2 = 7 \\ |\ph \wedge \alpha|^2 = 4 | \alpha |^2 \qquad \quad \,
& \qquad \qquad & \qquad \quad \, \, |\ps \wedge \alpha|^2 = 3
|\alpha |^2 \\ \st ( \ph \wedge \st ( \ph \wedge \alpha ) ) = -4
\alpha \, \, \quad \quad \quad & \qquad \qquad & \st (\ps \wedge \st
(\ps \wedge \alpha) ) = 3 \alpha \\ \ps \wedge \st (\ph \wedge
\alpha) = 0 \! \! \quad \quad \quad \quad \quad & \qquad \qquad &
\quad \, \, \ph \wedge \st (\ps \wedge \alpha) = -2 \ps \wedge \alpha
\\ \st (\ph \wedge \ws) = w \hk \ps \, \quad \quad \quad & \qquad
\qquad & \quad \quad \, \, \, \, \st ( \ps \wedge \ws) = w \hk \ph \\
\ph \wedge ( w \hk \ph) = -2 \st (w \hk \ph) \, & \qquad \qquad &
\quad \quad \ps \wedge (w \hk \ph) = 3 \st \ws \\ \ph \wedge (w \hk
\ps) = -4 \st \ws \quad \quad & \qquad \qquad & \quad \quad \ps
\wedge ( w \hk \ps) = 0
\end{eqnarray*}
\end{prop}

\subsection{Some coordinate-free identities of \G-geometry }
\label{evolutionidentitiessec}

We now derive some identities which are useful for describing the
decomposition of the space of forms in Section~\ref{formssec} and for
computing the evolution of the metric $g$ from the $3$-form $\ph$ in
Section~\ref{metricevolutionsec}.

\smallskip

\begin{lemma} \label{wtwscubedlemma}
Let $u$, $v$, and $w$ be vector fields on $M$. Let $\us$, $\vs$, and
$\ws$ denote their dual $1$-forms with respect to the metric $g$. Then
the following identity holds:
\begin{equation*}
\st \left( (u \hk \ph) \wedge (v \hk \ph) \wedge ( w \hk \ph) \right)
= -2 g(u, v) \ws - 2 g(u, w) \vs - 2 g(v, w) \us
\end{equation*}
\end{lemma}
\begin{proof}
We begin with the general relation between $\ph$, $g$, and the
volume form:
\begin{equation*}
(u \hk \ph ) \wedge (v \hk \ph) \wedge \ph = -6 g(u, v) \vol
\end{equation*} 
Taking the interior product of this equation with $w$ and using
the fact that $w \hk \vol = \st \ws$, we obtain
\begin{eqnarray*}
- 6 g(u, v) \st \ws & = & (w \hk u \hk \ph) \wedge (v \hk \ph)
\wedge \ph + (u \hk \ph) \wedge (w \hk v \hk \ph) \wedge \ph \\ & &
{}+ (u \hk \ph) \wedge (v \hk \ph) \wedge ( w \hk \ph) \\ & = & -2
(u \times w)^{\flat} \wedge \st (v \hk \ph) - 2 (v \times
w)^{\flat} \wedge \st(u \hk \ph) \\ & & {}+ (u \hk \ph) \wedge (v
\hk \ph) \wedge ( w \hk \ph)
\end{eqnarray*}
where we have used equation~\eqref{crosseq} and the relation $(v
\hk \ph) \wedge \ph = -2 \st( v \hk \ph)$ from
Proposition~\ref{g2relationsprop}. We rearrange this equation and
use $\st (v \hk \ph) = \vs \wedge \ps$ to obtain
\begin{eqnarray*}
(u \hk \ph) \wedge (v \hk \ph) \wedge ( w \hk \ph) & = & -6 g(u,v)
\st \ws -2 \vs \wedge (u \times w)^{\flat} \wedge \ps \\ & & {} -2
\us \wedge (v \times w)^{\flat} \wedge \ps
\end{eqnarray*}
We now use~\eqref{crosseq} and~\eqref{iteratedeq}, and take $\st$
of both sides to get
\begin{eqnarray*}
& & \st \left( (u \hk \ph) \wedge (v \hk \ph) \wedge ( w \hk \ph)
\right) \\ & & = -6 g(u,v) \ws -2 (v \times (u \times w) )^{\flat}
-2 (u \times (v \times w) )^{\flat} \\ & & = -6g(u,v) \ws -2 \left(
-g(u,v) \ws + g(v,w) \us - (v \hk u \hk w \hk \ps)^{\sharp}
\right) \\ & & {} -2 \left( -g(u,v) \ws + g(u,w) \vs - (u \hk v
\hk w \hk \ps)^{\sharp} \right) \\ & & = -2 g(u, v) \ws - 2 g(u,
w) \vs - 2 g(v, w) \us
\end{eqnarray*}
and the proof is complete.
\end{proof}

\begin{rmk} \label{wtwscubedrmk}
If we put $u = \ddx{i}$, $v = \ddx{j}$, and $w = \ddx{l}$ this
identity takes the form
\begin{equation*}
\st \left( (\ddx{i} \hk \ph) \wedge (\ddx{j} \hk \ph) \wedge (
\ddx{l} \hk \ph) \right) = -2 \left( g_{ij} g_{lm} + g_{il}
g_{jm} + g_{jl} g_{im} \right) \dx{m}
\end{equation*}
If we take $\st$ of both sides of this equation, and wedge both
sides with an arbitrary $1$-form $\alpha = \alpha_k \dx{k}$, we get
\begin{eqnarray*}
& & \alpha \wedge \left( (\ddx{i} \hk \ph) \wedge
(\ddx{j} \hk \ph) \wedge ( \ddx{l} \hk \ph) \right) = -2 \left(
g_{ij} g_{lm} + g_{il} g_{jm} + g_{jl} g_{im} \right) \alpha_k
\dx{k} \wedge \st \dx{m} \\ & & \frac{1}{8} \alpha_{s_7}
\ph_{is_1s_2} \ph_{js_3s_4} \ph_{ls_5s_6} \dx{s_1} \wedge \ldots
\wedge \dx{s_7} = -2 \left( g_{ij} g_{lm} + g_{il} g_{jm} + g_{jl}
g_{im} \right) \alpha_k g^{km} \vol
\end{eqnarray*}
and hence
\begin{eqnarray*}
& & \frac{1}{8} \sum_{\sigma \in S_7} \, \sgn(\sigma) \,
\ph_{i\sigma(1)\sigma(2)} \ph_{j\sigma(3)\sigma(4)}
\ph_{l\sigma(5)\sigma(6)} \alpha_{\sigma(7)} \, \dx{1} \wedge
\ldots \wedge \dx{7} \\ & & = -2 \left( g_{ij} \alpha_l + g_{il}
\alpha_j + g_{jl} \alpha_i \right) \sqrt{\det(g)} \, \dx{1} \wedge
\ldots \wedge \dx{7}
\end{eqnarray*}
which yields the useful relation
\begin{eqnarray} \nonumber
\sum_{\sigma \in S_7} \, \sgn(\sigma)
\ph_{i\sigma(1)\sigma(2)} \ph_{j\sigma(3)\sigma(4)}
\ph_{l\sigma(5)\sigma(6)} \alpha_{\sigma(7)} & = & -16 \left(
g_{ij} \alpha_l + g_{il} \alpha_j + g_{jl} \alpha_i \right)
\sqrt{\det(g)} \\ \label{wtwscubedeq} & = & \frac{8}{3} \left(
B_{ij} \alpha_l + B_{il} \alpha_j + B_{jl} \alpha_i \right)
\end{eqnarray}
using the fact that $B_{ij} = -6 g_{ij} \sqrt{\det(g)}$.
\end{rmk}

\smallskip

\begin{cor} \label{wtwscubedcor}
Let $u$, $v$, and $w$ be vector fields on $M$. Then
the following holds:
\begin{equation*}
\st \left( (v \hk w \hk \ph) \wedge (u \hk \ph) \wedge \ph \right)
= 2g(u,v) \ws - 2g(u,w) \vs + 2 \st (u \hk v \hk w \hk
\ps)^{\sharp}
\end{equation*}
\end{cor}
\begin{proof}
This can be seen in the proof of Lemma~\ref{wtwscubedlemma}
above, by observing what happens to the appropriate term.
\end{proof}

\begin{rmk} \label{wtwscubedcorrmk}
If we put $v = \ddx{l}$, $w = \ddx{i}$, and $u = \ddx{j}$ this
identity takes the form
\begin{equation*}
\st \left( (\ddx{l} \hk \ddx{i} \hk \ph) \wedge (\ddx{j} \hk \ph)
\wedge \ph \right) = 2 ( g_{lj} g_{im} - g_{ji} g_{lm} +
\ps_{iljm} ) \dx{m}
\end{equation*}
If we take $\st$ of both sides of this equation, and wedge both
sides with an arbitrary $1$-form $\alpha = \alpha_k \dx{k}$, we get
\begin{eqnarray*}
& & \alpha \wedge \left( (\ddx{l} \hk \ddx{i} \hk \ph) \wedge (\ddx{j} \hk \ph)
\wedge \ph \right) =  2 ( g_{lj} g_{im} - g_{ji} g_{lm} +
\ps_{iljm} ) \alpha_k \dx{k} \wedge \st \dx{m} \\ & & \frac{1}{12}
\alpha_{s_1} \ph_{ils_2} \ph_{js_3s_4} \ph_{s_5s_6s_7} \dx{s_1}
\wedge \ldots \wedge \dx{s_7} = 2 ( g_{lj} g_{im} - g_{ji} g_{lm} +
\ps_{iljm} ) \alpha_k g^{km} \vol
\end{eqnarray*}
and hence
\begin{eqnarray*}
& & \frac{1}{12} \sum_{\sigma \in S_7} \, \sgn(\sigma) \,
\alpha_{\sigma(1)} \ph_{il\sigma(2)} \ph_{j\sigma(3)\sigma(4)}
\ph_{\sigma(5)\sigma(6)\sigma(7)} \, \dx{1} \wedge
\ldots \wedge \dx{7} \\ & & = 2 \left( g_{lj} \alpha_i - g_{ji}
\alpha_l + \ps_{iljm} g^{km} \alpha_k \right) \sqrt{\det(g)} \,
\dx{1} \wedge \ldots \wedge \dx{7}
\end{eqnarray*}
which yields the useful relation
\begin{eqnarray} \nonumber
\sum_{\sigma \in S_7} \, \sgn(\sigma) \alpha_{\sigma(1)}
\ph_{il\sigma(2)} \ph_{j\sigma(3)\sigma(4)}
\ph_{\sigma(5)\sigma(6)\sigma(7)} & = & 24
\left( g_{lj} \alpha_i - g_{ji} \alpha_l + \ps_{iljm} g^{km}
\alpha_k \right) \sqrt{\det(g)} \\ \nonumber
 & = & -4 \left( B_{lj} \alpha_i - B_{ji} \alpha_l \right) \\
\label{wtwscubedcoreq} & & {} + 24 \, \ps_{iljm} g^{km}
\alpha_k \sqrt{\det(g)}
\end{eqnarray}
where we have again used $B_{ij} = -6 g_{ij} \sqrt{\det(g)}$.
\end{rmk}

\smallskip

\begin{prop} \label{vfzeroprop}
Let $u$, $v$, and $w$ be vector fields on $M$. Then
the following holds:
\begin{equation*}
(u \hk v \hk \ps) \wedge (w \hk \ph) \wedge \ph = (v \hk \ps)
\wedge (u \hk w \hk \ph) \wedge \ph
\end{equation*}
\end{prop}
\begin{proof}
We start with the (necessarily zero) $8$-form
\begin{equation*}
(v \hk \ps) \wedge (w \hk \ph) \wedge \ph = 0
\end{equation*}
and take the interior product with $u$. This gives
\begin{equation*}
(u \hk v \hk \ps) \wedge ( w \hk \ph) \wedge \ph = (v \hk \ps)
\wedge (u \hk w \hk \ph) \wedge \ph + (v \hk \ps) \wedge (w \hk
\ph) \wedge (u \hk \ph)
\end{equation*}
from which the result follows, using the fact that $(v \hk \ps)
\wedge (w \hk \ph) \wedge (u \hk \ph) = 0$ for any $u, v, w$ which
is proved in~\cite{K1}, Theorem 2.4.7.
\end{proof}

\begin{rmk} \label{vfzerormk}
If we put $u = \ddx{i}$, $w = \ddx{j}$, and $v = \ddx{l}$ this
identity takes the form
\begin{equation*}
(\ddx{i} \hk \ddx{l} \hk \ps) \wedge (\ddx{j} \hk \ph) \wedge \ph =
(\ddx{l} \hk \ps) \wedge (\ddx{i} \hk \ddx{j} \hk \ph) \wedge \ph
\end{equation*}
In coordinates this becomes
\begin{equation*}
\frac{1}{24} \ps_{lis_1s_2} \ph_{js_3s_4} \ph_{s_5s_6s_7} \dx{s_1}
\wedge \ldots \wedge \dx{s_7} = (\ddx{l} \hk \ps) \wedge (\ddx{i}
\hk \ddx{j} \hk \ph) \wedge \ph
\end{equation*}
and hence
\begin{equation} \label{vfzeroeq}
\sum_{\sigma \in S_7} \, \sgn(\sigma) \,
\ps_{li\sigma(1)\sigma(2)} \ph_{j\sigma(3)\sigma(4)}
\ph_{\sigma(5)\sigma(6)\sigma(7)} \qquad
\text{is skew-symmetric in $i,j$.}
\end{equation}
We also remark that the identity $(v \hk \ps)
\wedge (w \hk \ph) \wedge (u \hk \ph) = 0$ in local coordinates
becomes
\begin{equation} \label{vfzeroeq2}
\sum_{\sigma \in S_7} \, \sgn(\sigma) \,
\ph_{i\sigma(1)\sigma(2)} \ph_{j\sigma(3)\sigma(4)}
\ps_{l\sigma(5)\sigma(6)\sigma(7)} = 0 
\end{equation}
\end{rmk}

\smallskip

\begin{prop} \label{topformsprop}
Let $v,w,a,b,c,d$ be vector fields on $M$. The following identities hold:
\begin{eqnarray*}
\as \wedge \bs \wedge \cs \wedge \ps & = & \ph (a, b, c) \vol \\ \as
\wedge \bs \wedge \cs \wedge \ds \wedge \ph & = & \ps(a,b,c,d) \vol \\
\as \wedge \bs \wedge \cs \wedge \ws \wedge (v \hk \ps) & = & \left(
g(v,w) \ph(a,b,c) - g(a,v) \ph(w,b,c) \right. \\ & & {}- \left. g(b,v)
\ps(a,w,c) - g(c,v) \ph(a,b,w) \right) \vol \\ \as \wedge \bs \wedge
\cs \wedge \ds \wedge \ws \wedge ( v \hk \ph) & = & \left( g(v,w)
\ps(a,b,c,d) - g(a,v) \ps(w,b,c,d) \right. \\ & & {}- \left. g(b,v)
\ps(a,w,c,d) - g(c,v) \ps(a,b,w,d) \right. \\ & & {}- \left. g(d,v)
\ps(a,b,c,w) \right) \vol
\end{eqnarray*}
\end{prop}
\begin{proof}
The first two equations follow from repeated application
of~\eqref{iosrelationseq}. To prove the third, start with the
(necessarily zero) $8$-form
\begin{equation*}
\as \wedge \bs \wedge \cs \wedge \ws \wedge \ps = 0
\end{equation*}
and take the interior product with $v$, and rearrange terms. The
fourth equation is proved similarly.
\end{proof}

\smallskip

\begin{prop} \label{vecvecphipsiprop}
Let $a, b, c, d$ be vector fields on $M$. The following relation
holds:
\begin{equation*}
\as \wedge \bs \wedge (c \hk \ph) \wedge (d \hk \ps) = ( 2 g(a
\wedge b, c \wedge d) + \ps(a,b,c,d) ) \vol.
\end{equation*}
\end{prop}
\begin{proof}
We start with the relation $\ps \wedge (\cs \hk \ph) = 3 \st \cs$ from
Proposition~\ref{g2relationsprop}. Taking the interior product with
$d$, using Lemma~\ref{g2knownlemma}, and rearranging, we obtain
\begin{equation*}
(c \hk \ph) \wedge (d \hk \ps) = 3 \st ( \cs \wedge \ds ) - {(c
\times d)}^{\flat} \wedge \ps
\end{equation*}
Taking the wedge product with $\as \wedge \bs$,
\begin{eqnarray*}
\as \wedge \bs \wedge (c \hk \ph) \wedge (d \hk \ps) & = & 3 \, (\as
\wedge \bs ) \wedge \st (\cs \wedge \ds) - \as \wedge \bs \wedge {(c
\times d)}^{\flat} \wedge \ps \\ & = & 3 \, g(a \wedge b, c \wedge d)
\vol - \ph(a, b, c \times d) \vol \\ & = & 3 \, g(a \wedge b, c \wedge
d) \vol - g( a \times b, c\times d) \vol \\ & = & ( 2 \, g(a \wedge b,
c \wedge d) + \ps(a,b,c,d) ) \vol.
\end{eqnarray*}
In the second equality, we used $\us \wedge \vs \wedge \ws \wedge \ps
= \ph(u,v,w) \vol$ from Proposition~\ref{topformsprop}, in the third
equality we used $\ph(u,v,w) = g(u \times v, w)$, and in the final
equality we used Corollary~\ref{g2identitycor}.
\end{proof}

\subsection{Contractions of $\ph$ and $\ps$}
\label{contractionssec}

The next set of identities are various contractions of $\ph$, $\ps$,
and their derivatives, in index notation. They are used repeatedly in
calculations throughout this paper.

In local coordinates $x^1, x^2, \ldots, x^7$, the 3-form
$\ph$ and the dual 4-form $\ps$ are
\begin{eqnarray*}
\ph & = & \frac{1}{6} \ph_{ijk}\, dx^i \wedge dx^j \wedge dx^k \\
\ps & = & \frac{1}{24} \ps_{ijkl} \, dx^i \wedge dx^j \wedge dx^k
\wedge dx^l
\end{eqnarray*}
where $\ph_{ijk}$ and $\ps_{ijkl}$ are totally skew-symmetric in
their indices. The metric is given by $g_{ij} = g( \ddx{i},
\ddx{j})$. The cross product is a $(2,1)$ tensor which we write
as
\begin{equation} \label{crosscoeq}
\ddx{i} \times \ddx{j} = P^k_{ij} \ddx{k}
\end{equation}
where $P^k_{ij} = -P^k_{ji}$. From~\eqref{g2maineq} it follows that
\begin{equation} \label{phipeq}
\ph_{ijk} = g_{kl} P^l_{ij} \qquad \qquad \qquad P^l_{ij} =
g^{kl} \ph_{ijk}
\end{equation}
Setting $u = \ddx{i}$, $v = \ddx{j}$ and $w = \ddx{k}$
in~\eqref{iteratedeq}, we obtain
\begin{eqnarray}
\nonumber \ddx{i} \times \left( \ddx{j} \times \ddx{k} \right) & =
& -g_{ij} \ddx{k} + g_{ik} \ddx{j} + \ps_{ijkl} (dx^l)^{\sharp} \\
\label{iteratedcoeq} P^m_{il} P^l_{jk} \ddx{m} & = & -g_{ij}
\ddx{k} + g_{ik} \ddx{j} + \ps_{ijkl} g^{lm} \ddx{m}
\end{eqnarray}

The first set of identities is the following.

\smallskip

\begin{lemma} \label{g2identities1lemma}
Let the tensors $g$, $\ph$, $\ps$, and $P$ be as given above. Then
the following identities hold:
\begin{eqnarray*}
P^k_{il} P^l_{jk} & = & -6 g_{ij} \\ \ph_{ijk} \ph_{abc} g^{ia}
g^{jb} g^{kc} & = & 42 \\ \ph_{ijk} \ph_{abc} g^{jb} g^{kc} & = & 6
g_{ia} \\ \ph_{ijk} \ph_{abc} g^{kc} & = & g_{ia} g_{jb} - g_{ib}
g_{ja} - \ps_{ijab} 
\end{eqnarray*}
\end{lemma}
\begin{proof}
We prove the fourth equation. The other three follow by
contraction with $g^{ij}$ and using~\eqref{phipeq}. To obtain the
fourth equation, take the inner product of~\eqref{iteratedcoeq}
with $\ddx{n}$ and simplify:
\begin{eqnarray*}
P^m_{il} P^l_{jk} g_{mn} & = & - g_{ij} g_{kn} + g_{ik} g_{jn} +
\ps_{ijkn} \\ \ph_{ila}g^{am} \ph_{jkb} g^{bl} g_{mn} & = &
-g_{ij} g_{kn} + g_{ik} g_{jn} + \ps_{injk} \\ \ph_{iln} \ph_{jkb}
g^{bl} & = & -g_{ij} g_{kn} + g_{ik} g_{jn} + \ps_{injk}
\end{eqnarray*}
which is what we wanted to show, since $\ph_{iln} = -
\ph_{inl}$. We also note that the second identity above is just a
restatement in terms of indices of the fact that the pointwise
norm ${|\ph|}^2$ of
$\ph$ is 7.
\end{proof}

The next set of identities involves contractions of $\ph$ with
$\ps$.

\smallskip

\begin{lemma} \label{g2identities2lemma}
Let the tensors $g$, $\ph$, and $\ps$ be as given above. Then
the following identities hold:
\begin{eqnarray*}
\ph_{ijk} \ps_{abcd} g^{ib} g^{jc} g^{kd} & = & 0 \\
\ph_{ijk} \ps_{abcd} g^{jc} g^{kd} & = & - 4 \ph_{iab} \\ \ph_{ijk}
\ps_{abcd} g^{kd} & = & g_{ia} \ph_{jbc} + g_{ib} \ph_{ajc} +
g_{ic} \ph_{abj} \\ & & {} - g_{aj} \ph_{ibc} - g_{bj} \ph_{aic} -
g_{cj} \ph_{abi}
\end{eqnarray*}
\end{lemma}
\begin{proof}
Again, the first two follow from the third. To prove the third, we
take the inner product of~\eqref{crosscoeq}
with~\eqref{iteratedeq} and use Corollary~\ref{g2identitycor}:
\begin{equation*}
g\left( \ddx{a} \times \ddx{b}, \ddx{i} \times \left(\ddx{j} \times
\ddx{k} \right) \right) = g_{ai} \ph_{jkb} - g_{ib} \ph_{jka}
- \ps_{abil} P^l_{jk}
\end{equation*}
But this also equals
\begin{eqnarray*}
& = & g \left( P^l_{ab} \ddx{l} , -g_{ij} \ddx{k} + g_{ik} \ddx{j}
+ \ps_{ijkn} g^{nm} \ddx{m} \right) \\ & = & -g_{ij} g_{lk}
P^l_{ab} + g_{ik} g_{lj} P^l_{ab} + P^l_{ab} g^{nm} g_{lm}
\ps_{ijkn} \\ & = & -g_{ij} \ph_{abk} + g_{ik} \ph_{abj} +
P^l_{ab} \ps_{ijkl}
\end{eqnarray*}
Combining the two expressions and reaaranging, we obtain
\begin{equation*}
g_{ia} \ph_{jkb} - g_{ib} \ph_{jka} + g_{ij} \ph_{abk} - g_{ik}
\ph_{abj} - \ph_{jkc} \ps_{abil} g^{cl} - \ph_{abc} \ps_{ijkl}
g^{cl} = 0
\end{equation*}
Denote the above expression by $A_{ijkab}$. Then it is tedious but
straightforward to check that
\begin{equation*}
A_{ijkab} + A_{ajkbi} + A_{bijka} - A_{kijab} - A_{jkabi} = 0
\end{equation*}
yields the desired identity.
\end{proof}

Finally we can contract $\ps$ with itself.

\smallskip

\begin{lemma} \label{g2identities3lemma}
Let the tensors $g$, $\ph$, and $\ps$ be as given above. Then
the following identities hold:
\begin{eqnarray*}
\ps_{ijkl} \ps_{abcd} g^{ia} g^{jb} g^{kc} g^{ld} & = & 168 \\
\ps_{ijkl} \ps_{abcd} g^{jb} g^{kc} g^{ld} & = & 24 g_{ia} \\
\ps_{ijkl} \ps_{abcd} g^{kc} g^{ld} & = & 4 g_{ia} g_{jb} - 4
g_{ib} g_{ja} - 2 \ps_{ijab} \\ \ps_{ijkl} \ps_{abcd} g^{ld} & = &
-\ph_{ajk} \ph_{ibc} - \ph_{iak} \ph_{jbc} - \ph_{ija} \ph_{kbc} \\
& & {} + g_{ia} g_{jb} g_{kc} + g_{ib} g_{jc} g_{ka} + g_{ic}
g_{ja} g_{kb} \\ & & {} - g_{ia} g_{jc} g_{kb} - g_{ib} g_{ja}
g_{kc} - g_{ic} g_{jb} g_{ka} \\ & & {} -g_{ia} \ps_{jkbc} -
g_{ja} \ps_{kibc} - g_{ka} \ps_{ijbc} \\ & & {} + g_{ab}
\ps_{ijkc} - g_{ac} \ps_{ijkb}
\end{eqnarray*}
\end{lemma}
\begin{proof}
As usual, all the identities follow from the last one. (However,
this time establishing the third from the fourth also requires
using Lemma~\ref{g2identities1lemma}.) To prove the last identity,
we take the inner product of~\eqref{iteratedcoeq} with itself. The
calculation involves the use of Corollary~\ref{g2identitycor}
twice as well as Lemmas~\ref{g2identities1lemma}
and~\ref{g2identities2lemma}. The details are tedious and not
enlightening, hence they are left to the reader. As in the case of
$\ph$, the first identity is a restatement of the fact that
${|\ps|}^2 = 7$ pointwise.
\end{proof}

\begin{rmk} \label{psipsiidentityrmk}
The particular expression for $\ps_{ijkl} \ps_{abcd} g^{ld}$ above
is not manifestly skew-symmetric in $a,b,c$. This is due to the
particular way in which it was derived. Since it of course is
skew-symmetric in $a,b,c$, one can skew-symmetrize to obtain the
more natural expression $ \ps_{ijkl} \ps_{abcd} g^{ld} =$
\begin{eqnarray*}
& & g_{ia} g_{jb} g_{kc} + g_{ib} g_{jc} g_{ka} + g_{ic} g_{ja} g_{kb}
- g_{ia} g_{jc} g_{kb} - g_{ib} g_{ja} g_{kc} - g_{ic} g_{jb} g_{ka}
\\ & & {}-\frac{1}{3} \left( \ph_{ibc} \ph_{ajk} + \ph_{aic} \ph_{bjk}
+ \ph_{abi} \ph_{cjk} \right) - \frac{1}{3} \left( \ph_{jbc} \ph_{iak}
+ \ph_{ajc} \ph_{ibk} + \ph_{abj} \ph_{ick} \right) \\ & &
{}-\frac{1}{3} \left( \ph_{kbc} \ph_{ija} + \ph_{akc} \ph_{ijb} +
\ph_{abk} \ph_{ijc} \right) - \frac{1}{3} \left( g_{ia} \ps_{jkbc} +
g_{ib} \ps_{jkca} + g_{ic} \ps_{jkab} \right) \\ & & {} -\frac{1}{3}
\left( g_{ja} \ps_{kibc} + g_{jb} \ps_{kica} + g_{jc} \ps_{kiab}
\right) - \frac{1}{3} \left( g_{ka} \ps_{ijbc} + g_{kb} \ps_{ijca} +
g_{kc} \ps_{ijab} \right)
\end{eqnarray*}
but since this is much more unwieldly and contains non-integer
coefficients, we prefer to use the more manageable expression given in
Lemma~\ref{g2identities3lemma}.
\end{rmk}

\smallskip

Next we consider contractions involving the covariant derivatives of
$\ph$ and $\ps$.

\smallskip

\begin{prop} \label{g2derivativeidentitiesprop}
Let $g$, $\ph$, and $\ps$ be as before. The following identities hold:
\begin{eqnarray*}
( \nab{m} \ph_{ijk} ) \ph_{abc} g^{ia} g^{jb} g^{kc} & = & 0 \\ (
\nab{m} \ps_{ijkl} ) \ps_{abcd} g^{ia} g^{jb} g^{kc} g^{ld} & = & 0 \\
(\nab{m} \ph_{ijk}) \ps_{abcd} g^{ib} g^{jc} g^{kd} & = & - \ph_{ijk}
( \nab{m} \ps_{abcd} ) g^{ib} g^{jc} g^{kd}
\end{eqnarray*}
\begin{eqnarray*}
(\nab{m} \ph_{ijk}) \ph_{abc} g^{jb} g^{kc} & = & - \ph_{ijk} (\nab{m}
\ph_{abc} ) g^{jb} g^{kc} \\ (\nab{m} \ps_{ijkl}) \ps_{abcd} g^{jb}
g^{kc} g^{ld} & = & - \ps_{ijkl} (\nab{m} \ps_{abcd} ) g^{jb} g^{kc}
g^{ld} \\ (\nab{m} \ph_{ijk} ) \ps_{abcd} g^{jc} g^{kd} & = & - \ph_{ijk} 
(\nab{m} \ps_{abcd} ) g^{jc} g^{kd} - 4 \, \nab{m} \ph_{iab}
\end{eqnarray*}
and finally also
\begin{equation} \label{nablapsieq}
\nab{m} \ps_{abcd} = - (\nab{m} \ph_{abj} ) \ph_{cdk} g^{jk} -
\ph_{abj} (\nab{m} \ph_{cdk}) g^{jk}
\end{equation}
\end{prop}
\begin{proof}
The method for establishing all of these identities should be
clear. Simply take the covariant derivative $\nab{m}$ of the various
equations in Lemmas~\ref{g2identities1lemma},
\ref{g2identities2lemma}, and~\ref{g2identities3lemma}, and in some
cases relabel indices to obtain the desired result. We omit the
details.
\end{proof}

\smallskip

We also have the following important relation.
\begin{prop} \label{nabphirelnabpsiprop}
The following relation holds between $\nab{} \phi$ and $\nab{} \psi$:
\begin{equation*}
(\nab{m} \ps_{ijkl} ) \ps_{abcd} g^{jb} g^{kc} g^{ld} = 3 \,
(\nab{m} \ph_{ijk} ) \ph_{abc} g^{jb} g^{kc}
\end{equation*}
\end{prop}
\begin{proof}
We substitute~\eqref{nablapsieq} into the left hand side above, and
use Lemma~\ref{g2identities2lemma} to obtain:
\begin{eqnarray*}
& & -\left( (\nab{m} \ph_{ijp}) \ph_{klq} g^{pq} + \ph_{ijp} (\nab{m}
\ph_{klq}) g^{pq} \right) \ps_{abcd} g^{jb} g^{kc} g^{ld} \\ &
= & -(\nab{m} \ph_{ijp} ) g^{pq} g^{jb} (\ph_{qkl} \ps_{abcd}
g^{kc} g^{ld}) - (\nab{m} \ph_{klq} ) g^{pq} g^{kc} g^{ld} (
\ph_{pij} \ps_{cdab} g^{jb} ) \\ & = & {} - (\nab{m} \ph_{ijp} )
g^{pq} g^{jb} (-4 \, \ph_{qab}) \\ & & {} - (\nab{m} \ph_{klq}
) g^{pq} g^{kc} g^{ld} ( g_{pc} \ph_{ida} + g_{pd} \ph_{cia} +
g_{pa} \ph_{cdi} - g_{ic} \ph_{pda} - g_{id} \ph_{cpa} - g_{ia}
\ph_{cdp} ) \\ & = & 4 \, (\nab{m} \ph_{ijp} ) \ph_{abq} g^{pq}
g^{jb} + 0 + 0 \\ & & {}- (\nab{m} \ph_{klq}) ( \delta^q_a g^{kc}
g^{ld} \ph_{cdi} - \delta^k_i g^{pq} g^{ld} \ph_{pda} - \delta^l_i
g^{pq} g^{kc} \ph_{cpa} - g_{ia} g^{pq} g^{kc} g^{ld} \ph_{cdp})
\\ & = & 4 \, (\nab{m} \ph_{ijk} ) \ph_{abc} g^{jb} g^{kc} -
(\nab{m} \ph_{kla} ) \ph_{cdi} g^{kc} g^{ld} + (\nab{m}
\ph_{ilq} ) \ph_{pda} g^{pq} g^{ld} \\ & & {} + (\nab{m}
\ph_{kiq} ) \ph_{cpa} g^{pq} g^{kc} + (\nab{a} \ph_{klq} )
\ph_{cdp} g^{pq} g^{kc} g^{ld} g_{ia}
\end{eqnarray*}
Using Proposition~\ref{g2derivativeidentitiesprop} on the second and
final terms, the final term vanishes and the remaining terms all combine
(after relabelling some indices) to yield
\begin{equation*}
3 \, (\nab{m} \ph_{ijk} ) \ph_{abc} g^{jb} g^{kc}
\end{equation*}
and the proof is complete.
\end{proof}

\section{Review of Basic Flow Formulas} \label{reviewsec}

In this section we briefly review several formulas that are used
frequently when studying geometric flows. Let $(M, g)$ be an
oriented Riemannian manifold. We choose local coordinates $x^1,
\ldots, x^n$ and denote the Riemannian metric by $g_{ij} =
g(\ddx{i}, \ddx{j})$. The associated volume form is then
$\vol = \sqrt{\det(g)} \, \dx{1} \wedge \ldots \wedge \dx{n}$.
We use $g^{ij}$ to denote the inverse matrix to $g_{ij}$, and
$g^{ij}$ is the induced metric on $1$-forms: $g^{ij} = g( \dx{i},
\dx{j})$.

\smallskip

\begin{lemma} \label{invmetriclemma}
Suppose $g_{ij}$ depends smoothly on some parameter $t$. Then
\begin{equation*}
\ddt g^{ij} = - g^{ik} \left( \ddt g_{kl} \right) g^{lj} 
\end{equation*}
\end{lemma}
\begin{proof}
We differentiate the relation $g_{kl} g^{lj} = \delta^k_j$, to
obtain
\begin{eqnarray*}
& & \left( \ddt g_{kl} \right) g^{lj} + g_{kl} \left( \ddt g^{lj}
\right) = 0 \\ & & g^{ik} g_{kl} \left( \ddt g^{lj} \right) = \ddt
g^{ij} = - g^{ik} \left( \ddt g_{kl} \right) g^{lj}
\end{eqnarray*}
as claimed.
\end{proof}

\begin{lemma} \label{volumeflowlemma}
Suppose $g_{ij}$ depends smoothly on some parameter $t$. Then
\begin{eqnarray*}
\ddt \det(g) & = & \left( \ddt g_{ij} \right) g^{ij} \det(g) \\
\ddt \vol & = & \frac{1}{2} \left( \ddt g_{ij} \right) g^{ij}
\vol
\end{eqnarray*}
\end{lemma}
\begin{proof}
The second formula follows easily from the first, using the local
coordinate expression for the volume form. Cramer's rule from
linear algebra says
\begin{equation*}
g_{ik} G^{kj} = \det(g) \delta^j_i 
\end{equation*}
where $G^{ij}$ is the classical adjoint of $g_{ij}$, the transpose
of the matrix of cofactors. Note that $g^{ij} = \frac{1}{\det(g)}
G^{ij}$. The determinant $\det(g)$ is a {\em linear} function
$F(g_1, \ldots, g_n)$ of the columns $g_1, \ldots, g_n$ of the
matrix $g_{ij}$. Therefore
\begin{equation*}
\ddt \det(g) = F(\ddt g_1, g_2, \ldots, g_n) + F(g_1, \ddt g_2,
\ldots, g_n) + \ldots + F(g_1, g_2, \ldots, \ddt g_n)
\end{equation*}
where, for example, $F(\ddt g_1, g_2, \ldots, g_n)$ is the
determinant of the matrix $g_{ij}$ with the first column $g_{i1}$
replaced with $\ddt g_{i1}$. Now we expand the first of these
determinants along the first column, the second along the second
column, and so on. We obtain
\begin{eqnarray*}
\ddt \det(g) & = & \left( \ddt g_{i1} \right) G^{i1} + \left( \ddt
g_{i2} \right) G^{i2} + \ldots + \left( \ddt g_{in} \right) G^{in}
\\ & = & \left( \ddt g_{ij} \right) G^{ij} = \left( \ddt g_{ij}
\right) g^{ij} \det(g)
\end{eqnarray*}
which completes the proof.
\end{proof}

\end{document}